\documentclass[a4paper, twoside, 11pt]{article}  
\usepackage[latin1]{inputenc}
\usepackage[T1]{fontenc}
\usepackage{amssymb, amsmath, theorem, amsfonts, amscd, geometry,graphics,epsf,epsfig,psfrag,multicol,path,
changebar, array, delarray, enumerate, mathrsfs, setspace, fancyhdr, amssymb}

\usepackage[all]{xy}

\makeindex

\numberwithin{equation}{section}

\newcommand\qed{\hfill$\square$}

\newtheorem{The}{Theorem}[section]
\newtheorem{Pro}[The]{Proposition}
\newtheorem{Lem}[The]{Lemma}
\newtheorem{Def}[The]{Definition}
\newtheorem{Cor}[The]{Corollary}

\DeclareMathOperator{\End}{End}
\DeclareMathOperator{\Hom}{Hom}
\DeclareMathOperator{\Sym}{Sym}

\DeclareMathOperator{\Lie}{Lie}
\DeclareMathOperator{\Spec}{Spec}
\DeclareMathOperator{\Spf}{Spf}

\DeclareMathOperator{\Nilp}{Nilp}

\DeclareMathOperator{\Gal}{Gal}
\DeclareMathOperator{\id}{id}

\DeclareMathOperator{\diag}{diag}
\DeclareMathOperator{\dist}{dist}
\DeclareMathOperator{\univ}{univ}
\DeclareMathOperator{\PGL}{PGL}

\DeclareMathOperator{\Sch}{Sch_{\text{ln}}}
\DeclareMathOperator{\mult}{mult}
\DeclareMathOperator{\vol}{vol}
\DeclareMathOperator{\GL}{GL}
\DeclareMathOperator{\Sp}{Sp}
\DeclareMathOperator{\tr}{tr}
\DeclareMathOperator{\Ad}{Ad}
\DeclareMathOperator{\U}{U}

\newcommand\N{\mathbb{N}}
\newcommand\Q{\mathbb{Q}}
\newcommand\R{\mathbb{R}}
\newcommand\C{\mathbb{C}}
\newcommand\A{\mathbb{A}}
\newcommand\F{\mathbb{F}}
\newcommand\Z{\mathbb{Z}}
\newcommand\Fp{\mathbb{F}_p}
\newcommand\Fpq{\mathbb{F}_{p^2}}

\newcommand\Zp{\mathbb{Z}_p}
\newcommand\Zl{\mathbb{Z}_l}
\newcommand\Zlp{\mathbb{Z}_{(p)}}
\newcommand\Zpq{\mathbb{Z}_{p^2}}
\newcommand\Qp{\mathbb{Q}_p}
\newcommand\Qpq{\mathbb{Q}_{p^2}}
\newcommand\MHB{\mathcal{M}^{HB}}
\newcommand\M{\mathcal{M}}

\begin{document}
\title{Intersections of arithmetic Hirzebruch-Zagier cycles}
\author{Ulrich Terstiege}
\date{}
\maketitle

\thispagestyle{empty}
\section*{Introduction}

 In this paper we will  establish a close connection between the intersection multiplicity of three arithmetic Hirzebruch-Zagier cycles and the Fourier coefficients of the derivative  of a certain Siegel-Eisenstein series at its center of symmetry. 
Our main result proves a conjecture of Kudla and Rapoport.

 Kudla has proposed a general program which relates intersection multiplicities of special cycles on (arithmetic models of) Shimura varieties to Fourier coefficients (of derivatives) of Eisenstein series, see his paper \cite{Ku1}, his ICM-talk \cite{Ku2}, his Bourbaki-talk \cite{Ku3} and his CDM-talk \cite{Ku4}. Besides the paper \cite{Ku1}, which is based on the paper \cite{GK} by Gross and Keating, other examples of results in this direction may be found in the  papers \cite{KR1} - \cite{KR4} by Kudla and Rapoport, in the monograph \cite{KRY} by Kudla, Rapoport and Yang, in 
\cite{AR},  which is mainly an exposition of \cite{GK}, in the papers \cite{H1}, \cite{H2} and \cite{H3} by Howard and in  other papers.
The present paper contributes to this program in a situation of degenerate intersections.
\newline

We now describe our results in detail.

We fix a prime number $p \neq 2.$
Arithmetic Hirzebruch-Zagier cycles are defined in the paper \cite{KR2} by Kudla and  Rapoport  as cycles over a certain moduli scheme $\M$ of abelian schemes over $\Zlp.$
 The generic fibre of $\M$ is the canonical model over $\Q$ of a Hilbert-Blumenthal surface.

 Let us make this more precise.
Let $V$ be a quadratic space over $\Q$ with signature $(2,2)$. Let $C(V)$ be the corresponding Clifford algebra, and let $C^+(V)$ be its even part. Then $C^+(V)$ is of the form $B_0\otimes_{\Q}k$, where $B_0$ is an indefinite quaternion algebra over $\Q$, and $k,$ the center of $C^{+}(V),$ is a  real quadratic extension of $\Q$.

We suppose that there exists a self-dual $\Zlp$-lattice $\Lambda \subset V$(satisfying an additional technical condition, see section 1), and we fix such a lattice. Then $p$  is unramified  in $k.$ We assume that $p$ is inert in $k,$ in particular, $k$ is a field (see \cite{KR2}, section 11 for the case of a split prime).  Let $\mathcal{O}_C$ be the Clifford algebra of $\Lambda$ and let $\mathcal{O}_k$ be the ring of $\Zlp$-integers in $k$. Finally, we consider the  algebraic group $G$ over $\Q$ with 
\[
G(R)= \{ g \in (C^{+}(V)\otimes_{\Q}R)^{\times}; \   \nu(g)\in R^{\times}\}
\]
for any $\Q$-algebra $R$. Here 
 $ \nu$ denotes the spinor norm. We denote by $\A$ the  adeles of $\Q$, by   $\A_f$ the finite adeles of $\Q$ and by   $\A_f^p$ the finite adeles of $\Q$ with trivial $p$-adic component. Further $(\Sch / \Zlp)$ denotes the category of locally noetherian schemes over $\Spec \Zlp.$

If  $K^p \subset G(\A_f^p)$ is a sufficiently small compact open subgroup, there is a moduli scheme  $\M$ which represents the following functor on $(\Sch / \Zlp)$. It associates to a locally noetherian $\Zlp$-scheme $S$ the set of isomorphism classes of  tuples $(A, \lambda, \iota, \overline{\eta^p})$, where $A$ is an abelian scheme over $S$ up to prime to $p$ isogeny, $\lambda$ is a $\Zlp^{\times}$-class of $p$-principal polarizations on $A$,  further $\iota: \mathcal{O}_C \otimes \mathcal{O}_k \longrightarrow \End(A)\otimes \Zlp$ is a homomorphism (satisfying a compatibility condition with the Rosati involution induced by $\lambda$), and $\overline{\eta^p}$ is a $K^p$-level structure. See section 1 for the precise conditions.
Then $\M$ is quasi-projective and smooth of relative dimension $2$ over $\Spec \Zlp$.

To define arithmetic Hirzebruch-Zagier cycles, we recall that a 
 \emph{special endomorphism} of a tuple $(A, \lambda, \iota)$ over $S\in (\Sch/\Zlp)$ (as above) is an element $j \in \End(A)\otimes \Zlp$ such that $j=j^*$ for the Rosati involution and satisfying a compatibility condition with $\iota$, see section 1.
 If $S  $ is connected, 
then for any $(A, \lambda, \iota, \overline{\eta^p})$ over $S$ as above the  $\Zlp$-module of special endomorphisms of $(A, \lambda, \iota)$  becomes a quadratic $\Zlp$-module via the quadratic form $Q$ given by $j^2=Q(j)\cdot \id$.

Now let  $\omega \subset V(\A_f^p)$ be compact, open and stable under the $K^p$-action, and let $t \in \Q^{\times}$. We consider the functor on $(\Sch / \Zlp)$ which associates to  $S \in (\Sch / \Zlp)$  the set of isomorphism classes of $5$-tuples $(A, \lambda, \iota, \overline{\eta^p}, j)$, where $(A, \lambda, \iota, \overline{\eta^p})$ is as before, and where $j$ is a special endomorphism of $(A, \lambda, \iota)$  such that 
 $j^2=t$
 and  satisfying a compatibility condition with $\overline{\eta^p}$ and $\omega$.
 Again, if $K^p$ is sufficiently  small, this functor is representable by a scheme $Z(t, \omega)$ which maps  by a finite and unramified morphism (forgetting $j$) to $\mathcal{M}$. An \emph{arithmetic Hirzebruch-Zagier cycle} or \emph{special cycle} is a  scheme of the form $Z(t, \omega)$.

 We now fix three  special cycles $Z_1=Z(t_1, \omega_1)$, $ Z_2=Z(t_2, \omega_2)$ and  $Z_3=Z(t_3, \omega_3)$.   Let 
 \[
 Z=Z_1 \times_{\mathcal{M}} Z_2 \times_{\mathcal{M}}Z_3.
 \]
Locally, a special cycle is isomorphic to a (relative) divisor in $\M$ (or it is empty). This suggests to intersect three special cycles.
If $S \in (\Sch / \Zlp)$ is connected and $\xi =(A,\lambda, \iota, \overline{\eta^p}, j_1,j_2,j_3) \in Z(S)$, then following \cite{KR2} we associate to  $\xi$ its \emph{fundamental matrix} $T_{\xi}$ whose entry at $(a,b)$ is  given by $\frac{1}{2}(Q(j_a+j_b)-Q(j_a)-Q(j_b))$. It is a symmetric $ 3 \times 3$ matrix with entries in $\Q$.
 The map $\xi \mapsto T_{\xi}$ is locally constant, and hence we can write 
 \[
 Z=\coprod_T Z_T,
 \]
 where $Z_T$ is the union of the connected components of $Z$ which have fundamental matrix $T$.
 We assume that $Z$ is not empty, and we fix  $T$ such that $Z_T$ is not empty. Note that then all entries of $T$ lie in $\Zlp$ and that the diagonal entries of $T$ are $t_1,t_2,t_3.$ We define the intersection multiplicity 
 \[
 \chi_T(Z_1,Z_2,Z_3)= \chi(Z_T, \mathcal{O}_{Z_1}\otimes^{\mathbb{L}}\mathcal{O}_{Z_2}\otimes^{\mathbb{L}}\mathcal{O}_{Z_3}).
 \]
 Here, $\otimes^{\mathbb{L}}$ is the derived tensor product of $\mathcal{O}_{\mathcal{M}}$-modules and $\chi$ is the Euler-Poincaré  characteristic.

 If $T$ is singular, then the generic fibre of $Z_T$ might be nonempty. 
 Hence we additionally assume that  $T$ is nonsingular. Then $Z_T$ lies over the supersingular locus of $\M$ and its support is proper over $\Fp$ so that  $\chi_T(Z_1,Z_2,Z_3)$ is finite.    One of the main results of \cite{KR2} states that $Z_T$ is a discrete set of points if and only if $T$ (as matrix over $\Zlp$) is not divisible by $p$ and that it is  otherwise of dimension $1.$ 

 It is the aim of this paper to give an explicit expression for $ \chi_T(Z_1,Z_2,Z_3)$ and to express it in terms of certain Eisenstein series. The case that $Z_T$ is a discrete set of points was treated in \cite{KR2}, so we will focus on the case of a degenerate intersection, i.e., the case that $T$ is divisible by $p.$
 
 Let $P\subset \Sp_6$ be the Siegel parabolic subgroup. It is given as follows.
$P=MN,$ where
\[
\begin{split}
& M=\left\{m(a)=
\begin{pmatrix}
  a& 0 \\ 
 0 & ^t{a^{-1}}
\end{pmatrix} \mid a \in \GL(3) \right\} \\
 \text{and} \\ 
& N= \left\{n(b)=
\begin{pmatrix}
  1& b \\ 
 0 & 1
\end{pmatrix} \mid b \in \Sym_3  \right\}.
\end{split}
\]
Let $\tilde{K}=\tilde{K}_{\infty}\tilde{K}_f= \prod_v \tilde{K}_v$,
where
\[
 \tilde{K}_v = 
 \begin{cases}
\Sp_6(\Zl),\ & \text{if $v=l < \infty ;$ } \\
\\
\left\{\begin{pmatrix}
  a& b \\ 
 -b & a
\end{pmatrix} \mid a+ib \in \U_3(\R) ) \right\},
 \ & \text{if $v = \infty.$ } 
 \end{cases}
 \]
Then we have the Iwasawa decomposition
$$
\Sp_6(\A)=P(\A)\tilde{K}.
$$
To  $\omega= \omega_1\times \omega_2 \times \omega_3$ we will associate in section 6 a section $\Phi$ of the induced representation $I_3(s, \chi_V)$ of $\Sp_6(\A)$ (see \cite{Ku1}).
We consider the following Eisenstein series on $\Sp_6(\A)$ (as in \cite{KR2}, see also \cite{RW})  
 $$
E(g,s,\Phi)=\sum_{\gamma \in P(\Q) \setminus \Sp_6 (\Q) } \Phi(\gamma g,s).
$$
There is a Fourier decomposition (see \cite{Ku1})
$$  
E(g,s,\Phi)=\sum_{U \in \Sym_3(\Q)}E_U(g,s,\Phi).
$$
This yields 
a Fourier decomposition of its derivative 
$$  
E^{'}(g,s,\Phi)=\sum_{U \in \Sym_3(\Q)}E_U^{'}(g,s,\Phi).
$$
Recall our fixed nonsingular $T$ such that $Z_T$ is not empty.
The main theorem of this paper is the following. 
\begin{The}\label{Haupt}  Let $h \in \Sp_6(\R)$ and suppose that $\omega_1 \times \omega_2 \times \omega_3$ is locally centrally symmetric.
 There is  the following relation between the  derivative in $s=0$ of the  $T$-th Fourier coefficient of  $E(g,s,\Phi)$ and the  intersection multiplicity $\chi_T(Z_1,Z_2,Z_3):$
 \[
 \begin{split}
  E^{'}_{T}(h,0,\Phi)&= -\frac{1}{2}\log(p) \cdot \kappa \cdot \chi_T(Z_1,Z_2,Z_3) \cdot W_T^2(h),
 \end{split}
 \]
 where $\kappa$ is some  volume constant given in section 6.
 \end{The}
To explain this notation, we write $h=m(a)n(b)k \in M(\R)N(\R)\tilde{K}_{\infty}=\Sp_6(\R).$ Let $\tau=b+i \cdot a  ^{t}a.$ Then the Whittaker funktion $W_T^2(h)$ (see \cite{KR2}, \cite{Ku1}) is of the form  $$W_T^2(h)=c(a,k) e^{2\pi i \cdot \tr (T\tau)},$$ 
where $c(a,k)$ is a constant depending on $a$ and $k$.

The theorem was proved in \cite{KR2} in case that $Z_T$ is a discrete set of points or, equivalently, that $T$ is not divisible by $p$. 
\newline

We now   sketch the further content of this paper and  some of the techniques used to prove Theorem \ref{Haupt}.

Since our fixed $T$ is assumed to be nonsingular, $Z_T$ lies over the supersingular locus of the special fibre of $\mathcal{M}$ (\cite{KR2}, Corollary 3.8).
 Let $\F=\overline{\F}_p$  be a fixed algebraic closure of $\Fp,$ and let us fix  $\xi \in Z(\F)$ such that $T=T_{\xi}$.
Let $\xi =(A,\lambda, \iota, \overline{\eta^p}, j_1,j_2,j_3)$ and let $A(p^{\infty})$ be the corresponding $p$-divisible group. It is equipped with an action by $\mathcal{O}_C\otimes \mathcal{O}_k \otimes \Zp \cong M_4(\Zpq),$ which allows us to write $A(p^{\infty})=\mathcal{A}^4$.
Then $\mathcal{A}$ is a principally polarized supersingular  formal 
$p$-divisible group  over $\F$ of height 4 and dimension 2 
 which is equipped with a $\Zpq$-action satisfying some compatibility condition, see section 2. 

  Let $W=W(\F)$ be the ring of Witt vectors of $\F$.  For any $W$-scheme $S$ we write  $\overline{S}=S\times_{\Spec W}\Spec \F.$ We consider the following functor ${\mathcal{M}}^{HB}$ on the category 
$\begin{rm}Nilp \end{rm}$ of $W$-schemes $S$ such that $p$ is locally nilpotent in
$\mathcal{O}_S$. It associates to a scheme $S \in \Nilp$ 
 the set of isomorphism classes of pairs $(X, \varrho)$, consisting of a principally polarized $p$-divisible group $X$   over $S$ which is equipped with an $\Zpq$-action and a quasi-isogeny of height zero
$
\varrho: \mathcal{A}\times_{\Spec \F}\overline{S} \rightarrow X\times_S \overline{S},
$ such that certain compatibility conditions are satisfied, see section 2.

By Theorem 3.25 of  \cite{RZ}, the functor ${\mathcal{M}}^{HB}$ is representable by a formal scheme over $\Spf W$ which we
also denote by  ${\mathcal{M}}^{HB}$.
Let  $\widehat{\mathcal{M}_W/\mathcal{M}_W^{ss}}$ be the completion of $\mathcal{M}_W=\mathcal{M}\times_{\Spec \Zlp} \Spec W $ along its supersingular locus $\mathcal{M}_W^{ss}$. 

 Then  $\widehat{\mathcal{M}_W/\mathcal{M}_W^{ss}}$  and $\MHB$ are closely connected via  the  uniformization theorem
(Theorem 6.23 of loc. cit.), see (\ref{unif}).

The underlying reduced subscheme of $\MHB$ is a union of projective lines over $\F$  which  are in bijection with the vertices of the building $\mathscr{B}:=\mathscr{B}(\PGL_2(\Qpq))$. 
Two such projective lines $\mathbb{P}_{[\Lambda]}$ and $\mathbb{P}_{[\Lambda^{'}]}$ intersect in  a point if and only if the vertices $[\Lambda]$ and $[\Lambda^{'}]$ are neighbours in the building.

 As in \cite{KR2}, we define 
 $$
V^{'}=\{j \in  \End(A) \otimes \Q;  j=j^* \text{ and }  \iota(c \otimes a) \circ j= j\circ \iota(c \otimes a^{\sigma})\}.
$$
It is a quadratic space over $\Q$ of dimension $4.$ As above, the quadratic form is given by $j^2=Q(j)\cdot \id.$  We may regard $V^{'}$ as a $\Q$-subspace of $\End^0( \mathcal{A})= \End( \mathcal{A})\otimes \Q$.
We will also define a quadratic $\Qp$-space $V^{'}_p$ of special endomorphisms of $\mathcal{A}$ which contains $V^{'}$ as a dense $\Q$-subspace, cf. section 2.

For any $j\in V^{'} $ (in fact also for $j\in V_p^{'} $ ) we define the \emph{formal special cycle} $Z(j)$ to be the closed formal subscheme of $\MHB$ consisting of all points $(X, \varrho)$ such that $\varrho \circ j \circ \varrho^{-1}$ lifts to an endomorphism of $X$.

Recall our fixed three special endomorphisms $j_1,j_2,j_3$ associated to $\xi$. We define the intersection multiplicity $(Z(j_1),Z(j_2),Z(j_3))$ of the associated formal special cycles again to be the Euler-Poincar\'e characteristic of the derived tensor product of the structure sheaves of the $Z(j_i)$.
One shows (Proposition \ref{locglob}) 
 that 
 the global intersection multiplicity $\chi_T(Z_1,Z_2,Z_3)$ is a multiple of  $(Z(j_1),Z(j_2),Z(j_3))$ by an explicit integer. Thus in order to calculate $\chi_T(Z_1,Z_2,Z_3)$, it is enough to calculate  $(Z(j_1),Z(j_2),Z(j_3))$. 
\newline

To this end we first investigate the structure of formal special cycles. 
Let us for the moment fix $j \in V^{'}$ such that $j^2\neq 0$ and $a := \nu_p(j^2)\geq 0$, where $\nu_p$ is the $p$-adic valuation. 
First we show that $Z(j)$ is a relative divisor in $\MHB.$
We give a complete description of the special fibre of  $Z(j).$ To state the result, we recall from \cite{KR2} that $j$ gives rise to a morphism  $\beta$ of the building
Then $\nu_p(\beta^2)=a-1$. For $[\Lambda]\in \mathscr{B}$ denote by $d_{[\Lambda]}$ the distance of $[\Lambda]$ to $\mathscr{B}^{\beta}$, the fixed point locus of $\beta$. Then $\mathscr{B}^{\beta}$ is one single edge if $a$ is even, and it is a subbuilding of the  form $\mathscr{B}(\PGL_2(\Qp))$ if $a$ is odd.

 \begin{The}

The special fibre $Z(j)_p$ of $Z(j)$ is of the form
$$
Z(j)_p=\sum_{[\Lambda]; \  d_{[\Lambda]} \leq \frac{a-1}{2}}(1+p+\ ...\ +p^{\frac{a-1}{2} -d_{[\Lambda]}})\mathbb{P}_{[\Lambda]} \ + \ p^{a/2}\cdot s.
$$
  Here, $s$ is a divisor in the special fibre of ${\MHB}$  which is $\neq 0$ if and only if $a$ is even. In this case 
  $s$ is reduced and irreducible, and it meets the supersingular locus only in the unique supersingular point of  $Z(j/p^{a/2})$.
  
\end{The}
This theorem is proved using the theory of displays of Zink (\cite{Z1}, \cite{Z2}).

We also associate to any formal special cycle $Z(j)$ (such that $j^2 \in \Zp \setminus \{0\}$) a divisor $D(j),$ the "difference divisor",  whose ideal  is locally  generated by the quotient of a generator of the ideal of $Z(j)$ and  a generator of the ideal of $Z(j/p)$. We show that $D(j)$ is always a regular formal scheme. The proof of this uses the fact that for $j_1 \in V^{'}_p$ such that $\nu_p(j_1^2)=1,$ the formal special cycle $Z(j_1)$ is isomorphic to the formal model of the Drinfeld upper half plane (for $\Qp$ and base changed with $\Spf W$), and the fact that for $j \perp j_1$ the structure of $Z(j)\cap Z(j_1)$ is known from \cite{T}. This "testing by the Drinfeld space" is one of the basic techniques in our paper.
\newline

The next step is to investigate the structure of intersections of formal special cycles. First we show the statements of the following
\begin{The} 

\begin{enumerate}\label{th03}
\item Let $j, j^{'} \in V^{'}$ be linearly independent. Suppose further that for any  $x\in \MHB(\F)$  the equation in $x$ (that is a generator of the ideal in  $\mathcal{O}_{\MHB, x}$)  of at least one of the corresponding two formal special cycles  is not divisible by $p$. Then $\mathcal{O}_{Z(j)}\otimes^{\mathbb{L}}\mathcal{O}_{Z(j^{'})} $  is represented in the derived category by $ \mathcal{O}_{Z(j)}\otimes \mathcal{O}_{Z(j^{'}).}$

\item  The derived tensor product $ \mathcal{O}_{Z(j_1)}\otimes^{\mathbb{L}}\mathcal{O}_{Z(j_2)}\otimes^{\mathbb{L}}\mathcal{O}_{Z(j_3)}$ depends only on the $\Zlp$-span $\bf{j}$ of $j_1,j_2,j_3$ in $V^{'}$.
\item  If $T$ (as matrix over $\Zlp$) is not divisible by $p,$ then
$\mathcal{O}_{Z(j_1)}\otimes^{\mathbb{L}}\mathcal{O}_{Z(j_2)}\otimes^{\mathbb{L}} \mathcal{O}_{Z(j_3)}$ is represented in the derived category by $\mathcal{O}_{Z(j_1)}\otimes \mathcal{O}_{Z(j_2)}\otimes \mathcal{O}_{Z(j_3)}.$
\item Suppose that at least one of  $\nu_p(j_1^2), \nu_p(j_2^2), \nu_p(j_3^2)$ is odd. Then 
$$
(Z(j_1),Z(j_2),Z(j_3))=\sum_{l,m,n}(D(j_1/p^l),D(j_2/p^m),D(j_3/p^n)),
$$
where the sum is taken over all possible triples $(l,m,n)$ (i.e. setting $a_i=\nu_p(j_i^2)$, we have $l \leq [a_1/2],\  m \leq [a_2/2], \ n \leq [a_3/2]$, where $[ \ ]$ denotes Gauss brackets).
\end{enumerate}
\end{The}

Point 1 is the key observation for the proof of point 2. Since $p\neq 2,$ point 2 allows us to assume that $j_1,j_2,j_3$ are pairwise perpendicular to each other, i.e., $T$ is a diagonal matrix  $T=\diag(\varepsilon_1p^{a_1},\varepsilon_2p^{a_2},\varepsilon_3p^{a_3}),$ where  $\varepsilon_i \in \Zlp^{\times} $ for all $i$ and $a_1\leq a_2 \leq a_3.$ In this case the assumption of point 4 is fulfilled, which allows us by induction to restrict ourselves to the calculation of the intersection multiplicity $(D(j_1),D(j_2),D(j_3))$. 
Point 3 shows that (if $T$ is not divisible by $p,$) the length (over $\Zlp$) of  the local ring of a point in $Z_T$ resp. the length of the artinian $W$-scheme  $Z(j_1)\cap Z(j_2)\cap Z(j_3)$ is the same as the intersection multiplicity $(Z(j_1),Z(j_2),Z(j_3))$ in our sense (recall that in the situation of point 3, $Z_T$ is a discrete set of points and  the length of its local rings  is calculated in \cite{KR2}).

Given two perpendicular special endomorphisms $y_1,y_2 \in V^{'},$ we investigate the multiplicities of the several projective lines $\mathbb{P}_{[\Lambda]}$ in the intersections $D(j_1)\cap D(j_2)$ as a divisor in (say) $D(j_1)$ (Propositions \ref{unglmult}, \ref{glmult}). Further we investigate (the existence of) horizontal components of $D(j_1)\cap D(j_2)$, that is of components which do not have support in the special fibre (Propositions \ref{horiz}, \ref{horiz0}).
Using this we calculate the intersection multiplicity $(D(j_1),D(j_2),D(j_3))$ (resp. in some cases $(D(j_1),Z(j_2),D(j_3))$). By induction, this leads to the following

\begin{The}\label{expformi} Suppose that $T$ is $\GL_3(\Zlp)$-equivalent to $\diag(\varepsilon_1p^{a_1},\varepsilon_2p^{a_2},\varepsilon_3p^{a_3})$, where  $\varepsilon_i \in \Zlp^{\times} $ for all $i$ and $a_1\leq a_2 \leq a_3.$ Then
there is the following explicit expression for the intersection multiplicity  $(Z(j_1),Z(j_2),Z(j_3))$.
\[
\begin{split}
(Z(j_1),Z(j_2),Z(j_3))=& -\sum_{i=0}^{a_1}\sum_{j=0}^{\frac{a_1+a_2-\sigma}{2}-i}p^{i+j}(-1)^i (i+2j) \\
& - \eta \sum_{i=0}^{a_1}\sum_{j=0}^{\frac{a_1+a_2-\sigma}{2}-i}p^{\frac{a_1+a_2-\sigma}{2}-j}(-1)^{a_3+\sigma+i}(a_3+\sigma+i+2j) \\
& -  \tilde{\xi}^2 p^{\frac{a_1+a_2-\sigma}{2}+1}\sum_{i=0}^{a_1}\sum_{j=0}^{a_3-a_2+2\sigma-4}\tilde{\xi}^j
(-1)^{a_2-\sigma+i+j} (a_2-\sigma+2+i+j).
\end{split}
\]
\end{The}
See section 5 for the definition of the invariants $\eta  \in \{ \pm 1\}, \sigma \in \{ 1,2\}$ and $\tilde{\xi} \in \{ \pm 1, 0\}$ of $T$. As mentioned before, this theorem is proved by induction on $a_1+a_2+a_3.$ The induction start is given by the treatment of the cases $a_1=0$ in \cite{KR2} (together with point 3 of Theorem \ref{th03})  and $a_1=1$ in \cite{T}.
\newline

We may express $(Z(j_1),Z(j_2),Z(j_3))$ in terms of some representation densities.
To state the result, recall that
 for $S\in \begin{rm}  Sym \end{rm}_m(\mathbb{Z}_p)$ and 
$U\in \begin{rm}  Sym \end{rm}_n(\mathbb{Z}_p)$
with $\begin{rm}  det \end{rm}(S)\neq 0$ and $\begin{rm}  det \end{rm}(U)\neq 0$,
 the representation density is defined as 
\[
\alpha_p(S,U)= 
\operatorname*{lim}_{t\rightarrow\infty} p^{-tn(2m-n-1)/2} \mid \{x \in M_{m,n}(\mathbb{Z}/p^t\mathbb{Z});
\ S[x]-U \in p^t\begin{rm}  Sym \end{rm}_n(\mathbb{Z}_p)\}\mid.
\]
\begin{The}
 Let  $S=\begin{rm}  diag \end{rm}(1,-1,1,-\Delta),$ where $\Delta \in \Zp^{\times}$ is not a square. There is the following relation between intersection multiplicities and representation densities:
\[
(Z(j_1),Z(j_2),Z(j_3))= -\frac{p^4}{(p^2+1)(p^2-1)}\alpha_p^{'}(S,T).
\]
\end{The}
(See section 6 for an explanation of the derivative $\alpha_p^{'}(S,T).$)
This theorem was proved before in case $a_1=0$ in \cite{KR2} and in case $a_1=1$ in \cite{T}.
Using this theorem and the connection 
of $E^{'}_{T}(h,0,\Phi)$ and $\alpha_p^{'}(S,T)$ (see the proof of Theorem \ref{Ziel}),
Theorem \ref{Haupt} follows.
\newline

\emph{Acknowledgements.} I want to thank all people who helped me to write this paper. In particular I thank M. Rapoport for the suggestion of the problem and for his interest and advice. Thanks are also due to Prof. Kudla for helpful comments on Eisenstein series
and to U. G\"ortz for helpful discussions.

 
\newpage

\tableofcontents
\newpage

\section{Arithmetic Hirzebruch-Zagier cycles}

In this section we review the notion of arithmetic  Hirzebruch-Zagier cycles and introduce the intersection problem we want to consider. We closely follow the paper \cite{KR2} by Kudla and Rapoport to which we refer for more details. 

Let $V$ be a vector space  over $\Q$  of dimension $4$ which is equipped with a quadratic form $q$ such that the signature of $V$ with respect to $q$ is $(2,2)$. Denote by $C^{+}(V)$ the even part of the Clifford algebra $C(V)$. 
Let $\beta $ be the main involution of $C(V).$ Denote  the center of $C^{+}(V)$ by $k.$ It is a  real quadratic extension of $\Q.$
We consider the algebraic group $G$ over $\Q$ with 
\[
G(R)= \{ g \in (C^{+}(V)\otimes_{\Q}R)^{\times}; \   \nu(g)=g\cdot g^{\beta} \in R^{\times}\}
\]
for any $\Q$-algebra $R$.
Let   $K$ be a  compact open subgroup  of $G(\A_f).$ 
 One associates to $G$   and   $K$   a Shimura variety $Sh(G, D)_K$ with $\C$-valued points
\[
Sh(G,D)_K(\C)=G(\Q)\setminus[D\times G(\A_f)/K],
\]
where $D$ is the space of oriented negative $2$-planes in $V(\R).$ 
If $K$ is sufficiently small, then there is a canonical model $M_K$ over $\Q$ of  $Sh(G, D)_K$ which
 represents a moduli functor over $\Q$ whose points consist of polarized abelian schemes with endomorphism structure and $K$-level structure. See \cite{KR2} for the precise statement and  more details.
\newline

Now let $p\neq 2$ be a prime number. 
 We are interested in a model of $M_K$ over $\Zlp$. 
 
  Let $v \in V$ such that $q(v)>0.$ 
Then $C^{+}(V)$ can be written in  the form  $B_0\otimes_{\Q}k$, where $B_0 $ is the fixed algebra of $\Ad(v)$ in  $C^{+}(V).$ Further $B_0$ is an indefinite quaternion algebra over $\Q.$
We fix an element $\tau \in B_0^{\times}$ such that  $\tau^{\beta}=-\tau$
and $\tau^2<0.$ Then 
\begin{equation}
x \mapsto x^{*}=\tau x^{\beta}\tau^{-1}
\end{equation}
is a positive involution of $C(V).$  
 We assume that there exists a self-dual  $\Zlp$-lattice $\Lambda \subset V(\Q)$. Then  $p$ is unramified in $k.$  We further assume that $\Lambda= \tau\Lambda \tau^{-1}$, and 
 we fix such a lattice $\Lambda$. 
Denote by $\mathcal{O}_C$ the Clifford algebra of $\Lambda$, and denote by $\mathcal{O}_k$ the ring of $\Zlp$-integers in $k$. Then  $\mathcal{O}_C$ is invariant under $*.$  We assume that $p$ is inert in $k,$ in particular $k$ is a field. 
  We denote by $U$ the underlying $\Q$-vector space of $C(V).$ It is equipped with a nondegenerate alternating form given by $\langle x,y\rangle=\tr^0(y^{\beta}\tau x),$ where $\tr^0$ denotes the reduced trace on $C(V).$ 
 Further  $U$ inherits an action $i$ by $C(V)\otimes_{\Q} k$ via $i(c \otimes a)x=cxa.$ 
Finally we fix a compact open subgroup $K^p \subset G(\A_f^p)$.

We consider the moduli problem over $\Zlp$ which associates to a locally noetherian scheme $S$ over $\Zlp$ the set of isomorphism classes of $4$-tuples $(A, \lambda, \iota, \overline{\eta^p})$, where
\begin{enumerate}[(1)] 
\item $A$ is an abelian scheme over $S$, up to prime to $p$ isogeny,
\item $\lambda: A {\longrightarrow} \hat{A}$
is a $\Zlp^{\times}$-class of $p$-principal polarizations on $A,$
\item $\iota: \mathcal{O}_C \otimes \mathcal{O}_k \longrightarrow \End(A)\otimes \Zlp$ is a homomorphism satisfying
\[
\iota(c \otimes  a)^{*}=\iota(c^{*}\otimes a)
\] for the Rosati involution with respect to $\lambda$ on $\End^0(A)$ resp. the involution $*$ on $C(V)$ introduced above,
\item $\overline{\eta^p}$ is a $K^p$-equivalence class of  $ \mathcal{O}_C \otimes \mathcal{O}_k$-equivariant isomorphisms
\[
\left(\prod_{l \neq p} T_l(A)\right) \otimes \Q \stackrel{\sim}{\longrightarrow}U(\A_f^p).
\]
 We require that the elements of $\overline{\eta^p}$ preserve the symplectic forms on $\left( \prod_{l \neq p} T_l(A)\right) \otimes \Q $ (coming from $\lambda$), resp. on $U(\A_f^p)$ up to a scalar in $(\A_f^p)^{\times}$. 
 \end{enumerate}
We further require that the action of $\iota(c \otimes a)$ on $\Lie(A)$ satisfies the determinant condition, see \cite{KR2} resp. \cite{Ko}.  This implies in particular that the relative  dimension of $A$ over $S$ is $8$.
 \newline
 
 If $K^p$ is sufficiently small, this moduli problem is represented by a quasi-projective scheme $\mathcal{M}$ over $\Spec \Zlp$ which is smooth of relative dimension $2$ over $\Spec \Zlp$. If further $K=K_p \cdot K^p$, where $K_p \subset G(\Qp)$ is the stabilizer of $\Lambda \otimes \Zp$, then
 \[M_K \cong \mathcal{M}\times_{\Spec \Zlp}\Spec \Q.\] 
From now on, we suppose that $K$ has this form.

 \begin{Def}
 \emph{Let $(A, \lambda, \iota)$ (over a locally noetherian $\Zlp$-scheme $S$) as above. A} special endomorphism \emph{of $(A, \lambda, \iota)$  is an  element $j\in \End(A)\otimes \Zlp$ such that $j=j^*$ for the Rosati involution and  $\iota(c \otimes a) \circ j= j\circ \iota(c \otimes a^{\sigma})$
 for all $c \in \mathcal{O}_C $ and $a \in \mathcal{O}_k$. (Here $<\sigma>= \Gal(k/ \Q)$).}
 \end{Def}
 Suppose $S\in (\Sch / \Zlp)$ (the category of locally noetherian schemes over $\Spec \Zlp$) is connected. Then for any $(A, \lambda, \iota, \overline{\eta^p})$ over $S$ as above the  $\Zlp$-module of special endomorphisms of $(A, \lambda, \iota)$  becomes a quadratic $\Zlp$-module via the quadratic form $Q$ given by $j^2=Q(j)\cdot \id$.
\newline

Now we come to the notion of  \emph{special cycles} or \emph{arithmetic Hirzebruch-Zagier cycles}. 
 Let $\omega \subset V(\A_f^p)$ be compact, open and stable under the $K^p$-action, and let $t \in \Q^{\times}$. We consider the functor on $(\Sch / \Zlp)$ which associates to  $S \in (\Sch / \Zlp)$  the set of isomorphism classes of $5$-tuples $(A, \lambda, \iota, \overline{\eta^p}, j)$, where $(A, \lambda, \iota, \overline{\eta^p})$ is as before and $j$ is a special endomorphism of $(A, \lambda, \iota)$  such that 
 \begin{enumerate}
\item $Q(j)=t$
\item  \text{For $\eta \in \overline{\eta^p} $ we have} $\eta j {\eta}^{-1} \in \omega$.
 \end{enumerate}
 The second condition means that for any (equivalent: for one)  $\eta \in \overline{\eta^p}$ the endomorphism $\eta j {\eta}^{-1}$ of $U(\A_f^p)$ is given by right multiplication by an element of $\omega.$ 
 Again, if $K^p$ is small enough (which we always assume), this functor is representable by a scheme $Z(t, \omega)$ which maps  by a finite and unramified morphism (forgetting $j$) to $\mathcal{M}$. 
 
 \begin{Def}
\emph{A} special cycle \emph{or} arithmetic Hirzebruch-Zagier cycle \emph{over} $\mathcal{M}$ \emph{is a scheme of the form $Z(t, \omega)$ as described above.}
 \end{Def}
 
 We now fix three  special cycles $Z_1=Z(t_1, \omega_1)$, $ Z_2=Z(t_2, \omega_2)$ and  $Z_3=Z(t_3, \omega_3)$. Let 
 \[
 Z=Z_1 \times_{\mathcal{M}} Z_2 \times_{\mathcal{M}}Z_3.
 \]
 
 If $S \in (\Sch / \Zlp)$ is connected and $\xi =(A,\lambda, \iota, \overline{\eta^p}, j_1,j_2,j_3) \in Z(S),$ then we associate to $\xi$ its \emph{fundamental matrix} 
 $T_{\xi}.$ By definition, its entry at $(a,b)$ is  given by $\frac{1}{2}(Q(j_a+j_b)-Q(j_a)-Q(j_b))$. Thus it is a symmetric $ 3 \times 3$ matrix with entries in $\Q$.
Since the map $\xi \mapsto T_{\xi}$ is locally constant, we can write 
 \[
 Z=\coprod_T Z_T,
 \]
 where $Z_T$ is the union of the connected components of $Z$ which have fundamental matrix $T$.  We assume that $Z$ is not empty, and we fix  $T$ such that $Z_T$ is not empty. 
 Note that then all entries of $T$ lie in $\Zlp.$ We define the intersection multiplicity 
 \[
 \chi_T(Z_1,Z_2,Z_3)= \chi(Z_T, \mathcal{O}_{Z_1}\otimes^{\mathbb{L}}\mathcal{O}_{Z_2}\otimes^{\mathbb{L}}\mathcal{O}_{Z_3})
 \]
 Here, $\otimes^{\mathbb{L}}$ is the derived tensor product of $\mathcal{O}_{\mathcal{M}}$-modules and $\chi$ is the Euler-Poincaré  characteristic (defined analogous to \cite{KR1}, section 4:  Recall that this means the following. Let $\pi :  \M \rightarrow \Spec \Zlp$. Then for a   sheaf of modules  $\mathcal{F}$ on $\M$ one defines $\chi(\mathcal{F})= \sum_i(-1)^i \lg_{\Zlp} (R^i\pi_*\mathcal{F}).$ This  is finite if $\mathcal{F}$ is coherent and has  support in  the special fibre $\M_p$ of $\M$ and this support is proper over $\Fp$.  Further, $\chi$ is additive in short exact sequences.  For a (bounded) complex of sheaves of  modules $\mathcal{F}^{\bullet}$ on $\M$ one defines $\chi(\mathcal{F}^{\bullet})=\sum_i (-1)^{i}\chi(\mathcal{F}^i)$.) The first argument $Z_T$ of $\chi$ indicates  that we only consider the value of $\chi$ which comes from the part of $ \mathcal{O}_{Z_1}\otimes^{\mathbb{L}}\mathcal{O}_{Z_2}\otimes^{\mathbb{L}}\mathcal{O}_{Z_3}$ which has support in 
  $Z_T$.

 If $T$ is singular, then the generic fibre of $Z_T$ might be nonempty. 
   Hence we additionally assume that our fixed $T$ is nonsingular. Then (as follows from \cite{KR2}, Propossition 3.8 and sections 5 and 8), the support of $Z_T$ lies in the supersingular locus  of $\M$ and is proper over $\F_p.$ Using this one sees that then $\chi_T(Z_1,Z_2,Z_3)$ is finite. It is the aim of this paper to give an explicit expression for $ \chi_T(Z_1,Z_2,Z_3)$ in terms of certain Eisenstein series.


\section{Formal special cycles and the local intersection problem}

In this section the intersection problem is reformulated in terms of an intersection problem of  formal special cycles on a formal moduli space of $p$-divisible groups. We will solve our intersection problem in the subsequent sections in this reformulation.

Since our fixed $T$ is assumed to be nonsingular, $Z_T$ lies over the supersingular locus of the special fibre of $\mathcal{M}$ (\cite{KR2}, Corollary 3.8).
 Let $\F=\overline{\F}_p$   be a fixed algebraic closure of $\Fp,$ and let us fix  $\xi \in Z(\F)$ such that $T=T_{\xi}$. We write $\xi =(A,\lambda, \iota, \overline{\eta^p}, j_1,j_2,j_3).$
 As in \cite{KR2}, we define 
 $$
V^{'}=\{j \in \End^0(A);  j=j^* \text{ and }  \iota(c \otimes a) \circ j= j\circ \iota(c \otimes a^{\sigma})\ \ \forall \ c\otimes a \in \mathcal{O}_C\otimes \mathcal{O}_k\}.
$$
 Note that $V^{'}$ depends only on the isogeny class of $(A,\lambda, \iota).$ Further, as shown in \cite{KR2}, its dimension as $\Q$-vector space is $4$, and it is equipped with the quadratic form $Q$ given by $j^2=Q(j)\cdot id$.

 Let $A(p^{\infty})$ be the  $p$-divisible group corresponding to $A.$ It is equipped with an action by $\mathcal{O}_C\otimes \mathcal{O}_k \otimes \Zp \cong M_4(\Zpq)$ which allows us to write $A(p^{\infty})=\mathcal{A}^4$.
Then $\mathcal{A}$ is a supersingular  formal 
$p$-divisible group  over $\F$ of height 4 and dimension 2 which is equipped with an action 
\[
\iota:   \mathbb{Z}_{p^2}\rightarrow \begin{rm}End\end{rm}(\mathcal{A}),
\]
 such that $\mathcal{A}$ is special with respect to $\iota$ (see below for a definition of the term \emph{special}). Further $\mathcal{A}$ is equipped 
  with a principal polarization

\begin{equation*}
\begin{CD}
\lambda: \mathcal{A}  \stackrel{\sim}{\longrightarrow}  \hat{\mathcal{A}},
\end{CD}
\end{equation*}
such that for the Rosati involution $\iota(a)^{*}=\iota(a)$ for all $a \in \Zpq$.

  Let $W=W(\F)$ be the ring of Witt vectors of $\F$ and denote its Frobenius also by $\sigma$.  For any $W$-scheme $S$ we write  $\overline{S}=S\times_{\Spec W}\Spec \F.$ We consider the following functor ${\mathcal{M}}^{HB}$ on the category 
$\begin{rm}Nilp \end{rm}$ of $W$-schemes $S$ such that $p$ is locally nilpotent in
$\mathcal{O}_S$. It associates to a scheme $S \in \Nilp$ 
 the set of isomorphism classes of the following data.
\begin{enumerate}[(1)]
\item A $p$-divisible group $X$ over $S$, with an action
\[
\iota: \mathbb{Z}_{p^2}\rightarrow \begin{rm}End\end{rm}(X),
\]
such that $X$ is special with respect to this $\mathbb{Z}_{p^2}$-action.
\item 
A quasi-isogeny of height zero
\[
\varrho: \mathcal{A}\times_{\Spec \F}\overline{S} \rightarrow X\times_S \overline{S},
\]
which commutes with the action of $\mathbb{Z}_{p^2}$   such that the following condition holds.
Let   ${\lambda}_{\overline{S}}: \mathcal{A}_{\overline{S}} \rightarrow
 \hat{\mathcal{A}}_{\overline{S}}$ be the map induced by $\lambda$.
Then we require the existence of an isomorphism $\tilde{\lambda}
:X \rightarrow \hat{X}$
such that for the induced map 
  $\tilde{\lambda}_{\overline{S}}:X_{\overline{S}} \rightarrow \hat{X}_{\overline{S}}$ we have 
  the relation 
${\lambda}_{\overline{S}}=\hat{\varrho}\circ \tilde{\lambda}_{\overline{S}}\circ \varrho$.

\end{enumerate}
Here, a $p$-divisible group $X$ over $S$ with $\mathbb{Z}_{p^2}$-action is said to be {\em special} if the induced $\mathbb{Z}_{p^2}\otimes \mathcal{O}_S$-module $\Lie X$ is, locally on $S$, free of rank $1$.

 Let $\mathcal{M}^{\bullet}$ be the same functor on Nilp as $\MHB$ but  without the condition on the height of the quasi-isogeny, and such that the relation ${\lambda}_{\overline{S}}=\hat{\varrho}\circ \tilde{\lambda}_{\overline{S}}\circ \varrho$ is required  to hold only up to a constant in $\Qp^{\times}$. Then by  
 Theorem 3.25 of \cite{RZ}, the functor $\mathcal{M}^{\bullet}$ is representable by a formal scheme which we
also denote by $\mathcal{M}^{\bullet}$. From this it follows that $\mathcal{M}^{HB}$  also also representable by a formal scheme which we
also denote by $\mathcal{M}^{HB}$, and $\mathcal{M}^{\bullet} = \MHB \times \Z.$
 The formal schemes $\M^{\bullet}$ and $\MHB$ are formally locally of finite type over $\Spf W $ and are formally smooth  over $\Spf W $.

Denote by $\widehat{\mathcal{M}_W/\mathcal{M}_W^{ss}}$ the completion of $\mathcal{M}_W=\mathcal{M}\times_{\Spec \Zlp} \Spec W $ along its supersingular locus $\mathcal{M}_W^{ss}$. Then   by Theorem 6.23 of loc. cit., we have (using that our  $K^p$ is sufficiently  small) an isomorphism of formal schemes over $\Spf W$
\begin{equation}\label{unif}
\widehat{\mathcal{M}_W/\mathcal{M}_W^{ss}} \cong G^{'}(\Q)\setminus(\mathcal{M}^{HB}\times \Z \times G(\A_f^p)/K^p),  
\end{equation}
where $G^{'}$ is the inner form of $G$ described in \cite{KR2}, section 4, i.e., $G^{'}$ is defined as $G$ for $B=C^+(V)= B_0 \otimes_{\Q}k$, but for the quaternion algebra $B^{'}$ over $k$ which is ramified at the two archimedean primes and is isomorphic to $B$ at all finite primes.  Then $G^{'}(\Q)$ can be identified with the group of quasi-isogenies of the abelian variety $A,$ which respect $\iota$ and (up to scalar) $\lambda.$ Further $g \in G^{'}(\Q)$ acts on $\mathcal{M}^{\bullet} $ (and therefore on  $\MHB \times \Z$) by sending $(X, \varrho)$ to $(X, \varrho g^{-1}).$
See \cite{RZ} for more details.
\newline
Denote the isocrystal of $\mathcal{A}$  by $N$. It is equipped with $\sigma$- resp. $\sigma^{-1}$-linear operators $F$ and $V$ such that $FV=VF=p.$  From the      
polarization $\lambda$ we get a perfect symplectic form $\langle  ,\rangle  $ on the Dieudonné module of $\mathcal{A}$ and hence also on $N$. 
As in \cite{KR2}, section 5, we define in this context the space of special endomorphisms of $\mathcal{A}$
\[
V_p^{'}=\{j \in \End(N,F); j\iota(a)=\iota(a^{\sigma})j \ \  \forall a \in \Zpq  \ \text{ and } \ j^*=j\},
\]
where  $*$ denotes 
the adjoint with
respect to the alternating form $\langle  ,\rangle  .$ 
We may and will regard $V^{'}$ as the  $\Q$-subspace of $V^{'}_p$ consisting of the special endomorphisms which come from the abelian variety $A$.  Then $V^{'} $ is dense in $V_p^{'}$. 
As shown in loc. cit., $V_p^{'}$ is a $4$-dimensional vector space over
$\mathbb{Q}_p$ with  a quadratic form, which we also denote by  $Q$ and which is  again given by 
\[
j^2=Q(j)\cdot \id.
\]
Instead of $Q(j)$ we will also write $j^2.$ Note that $V^{'}$ is positive definite (\cite{KR2}, Proposition 3.5), that $V^{'}(\Qp)\cong V^{'}_p$ and that  $V^{'}(\A_f^p) \cong V(\A_f^p)$ via  some (fixed) $\eta \in \overline{\eta^p}$ (where $\overline{\eta^p}$ belongs to $\xi$).  We will always identify $V^{'}(\A_f^p)$ and $V(\A_f^p).$
 In particular we may regard $V^{'}$ as being contained in $V(\A_f^p).$

Given special endomorphisms $y_1,...,y_n\in V_p^{'}$, we again define their {\em fundamental matrix} to be the symmetric $n\times n$-matrix whose entry at $(a,b)$ is given by $\frac{1}{2}(Q(y_a+y_b)-Q(y_a)-Q(y_b))$.
\begin{Def}
\emph{Let $j \in V_p^{'}$ and regard it as an element of $\End^0(\mathcal{A})$. Then the} formal special cycle $Z(j)$ \emph{associated to $j$ is the closed formal subscheme of ${\mathcal{M}}^{HB}$ consisting of
all points $(X, \varrho)$ such that $\varrho \circ j \circ \varrho^{-1}$ lifts to an endomorphism of $X$.}
\end{Def}
The fact that $Z(j)$ is a closed formal subscheme of ${\mathcal{M}}^{HB}$ follows from \cite{RZ}, Proposition 2.9.

Our next aim is a description of the underlying reduced scheme ${\MHB}^{red}$ of $\MHB$. 
Let $X$ be the set of 
 $W$-lattices $L$  in the isocrystal $N$ such that  for the dual lattice $L^{\perp}$ we have  $L^{\perp}=cL$  (with respect to $\langle, \rangle$) for some constant $c$,  and which are stable under $F$, under $V$  and under the $\Zpq$-action. The latter gives a grading $L=L_0\oplus L_1 \subset N=N_0\oplus N_1$ (see also \cite{KR2}, section 4). We call the index $i\in \Z/2$ critical for $L$ if $F^2L_i=pL_i$. 
 By \cite{KR2}, Lemma 4.2, for each such lattice at least one index is critical. We call $L$ superspecial if both indices are critical.
We may assume that $\xi$ is superspecial, which is means that  the corresponding Dieudonné  module $M \subset N$ of $\mathcal{A}$ is superspecial. (It follows for example from \cite{KR2}, Lemma 8.12 or from the beginning of the next section that there is always a superspecial point in $Z_T.$)
 
 The $\F$-valued points of $\MHB$ correspond to those lattices in $X$ which have the same volume as $M$ (by the condition that the quasi-isogenies have height $0$). We denote the set of such lattices by $X^{HB}$.

Let now $\eta^0_i=M_i^{p^{-1}F^2}$ for $i \in \Z/2$.  If $L \in X^{HB}$ has critical index $0,$ then let $\eta_0^L=L_0^{p^{-1}F^2}$. It has the same volume as $\eta_0^0$. 
We assign to $L$ the point on the projective line $\mathbb{P}(L_0/pL_0)= \mathbb{P}( \eta_0^L \otimes_{\Zpq}\F)$ corresponding to the line $FL_1/pL_0 \subset L_0/pL_0= \eta_0^L \otimes_{\Zpq}\F$.

 If $L \in X^{HB}$ has critical index $1,$ then let $\eta_1^L=L_1^{p^{-1}F^2}$. It has the same volume as $\eta_1^0$.  We assign to $L$ the point on $\mathbb{P}(L_1/pL_1)= \mathbb{P}( \eta_1^L \otimes_{\Zpq}\F)$ corresponding to the line $FL_0/pL_1 \subset L_1/pL_1 = \eta_1^L \otimes_{\Zpq}\F$. If $e_1,e_2$ is a basis of $\eta_0^0,$ then $F\eta_1^L$ has the same volume as the $\Zpq$-lattice in $\eta^0_0\otimes_{\Zp}\Qp$ spanned by $e_1, pe_2$.

In both cases, any point on $\mathbb{P}(L_i/pL_i)$ is recovered by some lattice in $X^{HB}$, see also \cite{KR2}, Lemma 4.3.  
If we identify $\eta^0_0\otimes_{\Zp}\Qp$ with $\Qpq^2$, we see that the projective lines occuring in the above construction for the lattices in $X^{HB}$ are in bijection with the vertices of the building $\mathscr{B}:=\mathscr{B}(\PGL_2(\Qpq))$.

For each  superspecial lattice we get two points, each lying on a different projective line. 
 If we glue  for  any superspecial lattice the corresponding two  projective lines along these ($\Fpq$-rational) points, then we see that (at least on the level of $\F$-valued points) ${\MHB}^{red}$ is a union of projective lines dual to  $\mathscr{B}$: For each vertex $[\Lambda] \in \mathscr{B}$ we have a projective line $\mathbb{P}_{[\Lambda]}$ (of the form $\mathbb{P}( \eta_1^L \otimes_{\Zpq}\F)$ or $\mathbb{P}( \eta_0^L \otimes_{\Zpq}\F)$ for a suitable $L$ as above) over $\F$ and with a $\Fpq$-rational structure. Further for each edge, the  projective lines corresponding to the vertices of the edge intersect in a $\Fpq$-rational point. All $\Fpq$-rational points are such intersection points and these correspond precisely to  the superspecial points (comp. \cite{KR2}, Lemma 4.3).
 
  To see that this description is also true on the level of schemes, we show that the $\F$-valued points on each projective line are in fact  the $\F$-valued points of a closed subscheme of ${\MHB}^{red}$  which equals the projective line as a scheme. To this end one chooses three pairwise perpendicular special endomorphisms $y_1,y_2, y_3 \in V^{'}_p$ such that $y_1^2=p$ and such that  the scheme $(Z(y_1)\cap Z(y_2)\cap Z(y_3))^{red}$ is (at least on the level of points) precisely the projective line. (It is easy to see that such $y_1,y_2,y_3$ exist using the description of $V^{'}$ given in the beginning of the next section and Lemma \ref{hillem} or using [KR2], Lemma 8.12.) By the identification of $Z(y_1)$ and the Drinfeld upper half plane (whose underlying reduced subscheme is known to be a union of projective lines) given in \cite{T}, section 4, it follows that this equality also holds on the level of schemes. (See also \cite{KR2}, (4.10).)
\newline

We now come back to our intersection problem. 
Let $S_1, S_2, S_3$ be closed subschemes  of $\MHB$. Then we define the (possibly infinite) intersection multiplicity
\[
(S_1,S_2,S_3):= \chi(\mathcal{M}^{HB},\mathcal{O}_{S_1}\otimes^{\mathbb{L}}\mathcal{O}_{S_2}
\otimes^{\mathbb{L}}\mathcal{O}_{S_3}),
\]
 Here, $\otimes^{\mathbb{L}}$ is the derived tensor product of $\mathcal{O}_{\MHB}$-modules and $\chi$ is the Euler-Poincaré  characteristic. 
 Again, this is finite if $S_1 \cap S_2 \cap S_3$ has support in the supersingular locus and this support is proper over $\F.$

Recall our fixed $\xi \in Z_T(\F)$. 
Writing again $\xi =(A,\lambda, \iota, \eta, j_1,j_2,j_3)$, the special endomorphisms $j_1,j_2,j_3$ induce special endomorphisms of $\mathcal{A}$ which we also call $j_1,j_2,j_3$.

Our next aim is to express the intersection multiplicity $\chi_T(Z_1,Z_2,Z_3)$ in terms of the intersection multiplicity 
$(Z(j_1),Z(j_2),Z(j_3))$, the calculation of which we will call the local intersection problem.
Denote by $Z^{'}$ the center of $G^{'}$ and denote by $Z^{'}(\Q)^0$ the set of elements  $z \in Z^{'}(\Q)$ such that $\nu_p(\det(z))=0$ (where $\nu_p$ denotes the $p$-adic valuation). Let $\underline{x} =(j_1,j_2,j_3)\in  V^{' 3}$. Further let 
\[ I(\underline{x}, \omega_1\times \omega_2 \times \omega_3)= \{  gK^p \in G(\A_f^p)/K^p \ ; \ \ g^{-1}\underline{x}g \in  \omega_1\times \omega_2 \times \omega_3 \}.\]

\begin{Pro}\label{locglob}
The intersection multiplicity $\chi_T(Z_1,Z_2,Z_3)$ can be expressed in the following way,
\[
\chi_T(Z_1,Z_2,Z_3)=2\cdot (Z(j_1),Z(j_2),Z(j_3)) \cdot \mid Z^{'}(\Q)^0 \setminus  I(\underline{x}, \omega_1\times \omega_2 \times \omega_3)  \mid.
\]

\end{Pro}

The proof is analogous to the proof of the corresponding Theorem 8.5  in \cite{KR1}. As in loc. cit., it is based on the following incindence relations.

Let $Z^{'}=Z(t^{'}, \omega^{'})$ be a special cycle. 
For a  special endomorphism $j \in V^{'}_p$ define the closed subscheme $Z^{\bullet}(j)$ of $\mathcal{M}^{\bullet}$ by the obvious analogon to $Z(j)$. Further denote by $V^{'}_{t^{'}}$ the set of elements $j$ in $V^{'}$ with $j^2=t^{'}$. Denote by  $\widehat{Z^{'}_W}$ the module over  $\widehat{\mathcal{M}_W/\mathcal{M}_W^{ss}}$ obtained from $Z^{'}$ by completing $Z^{'} \times_{\Spec \Zlp}\Spec W$ along $\mathcal{M}_W^{ss}.$ 
Then the uniformization morphism (\ref{unif}) induces an  inclusion
$$
\widehat{Z^{'}_W} \hookrightarrow G^{'}(\Q)\setminus (V^{'}_{t^{'}} \times \mathcal{M}^{\bullet}\times G(\A_{f}^p)/K^p),
$$
and a point $(j, (X, \varrho), gK^p)$ lies in the image if and only if
$ g^{-1}jg \in \omega^{'}  \ \ \text{and} \ \ (X,\varrho)\in Z^{\bullet}(j).  $ For several special cycles $Z^{'}_i$ we get an analogous description for the fibre product of the  $\widehat{{Z_i^{'}}_W}$. 
See the discussion in \cite{KR1}, section $8$ for more details.
\newline

The proposition shows that it is enough to handle the local intersection problem and this is what we will do in the subsequent sections.
 

\section{On the structure of formal special cycles}

Our first object of study in this section is the underlying reduced subscheme of a formal special cycle. Here the results follow from \cite{KR2}, sections 5 and 8. We denote by $\nu_p$ the $p$-adic valuation.

 Let $j \in V_p^{'}$ such that $j^2 \neq 0$ and $\nu_p(j^2) \geq0$. Let $Y= N_0^{p^{-1}F^2}$. Then $\beta:=F^{-1}j\mid Y$ is an $\sigma$-linear endomorphism and induces an endomorphism of $\mathscr{B}.$   If $\nu_p(j^2)=0,$ then the underlying reduced subscheme of $Z(j)$ is a single superspecial  point. It is the superspecial  point corresponding the unique edge $\mathscr{B}^{\beta}$ in $\mathscr{B}$ which is fixed by $\beta$. If $\nu_p(j^2)\geq 1,$ then the underlying reduced subscheme of $Z(j)$ is  a connected union of projective lines. The corresponding locus in the building is given by the set 
 \[
 \mathcal{T}(\beta)=\{x \in \mathscr{B}\ ; \ \dist(x, \mathscr{B}^{\beta})\leq \frac{1}{2}\nu_p(\det \beta)\}.
 \]
Here, $\mathscr{B}^{\beta}$ is the fixed point set of $\beta$ in $\mathscr{B}$, and $\dist$ denotes the distance in the building.
The fixed point set  $\mathscr{B}^{\beta}$ is the midpoint of an edge if $\nu_p(j^2)$ is even and it is a subbuilding of the form $\mathscr{B}({\PGL}_2(\Q_p))$ if $\nu_p(j^2)$ is odd.
See \cite{KR2}, sections 5 and 8  for more details. 
\begin{Def}
\emph{ For a special endomorphism} $j \in V^{'}_p$ \emph{such that $j^2\neq 0$ denote by} $Core(j)$ \emph{ the set of superspecial points in $\MHB$ such that the corresponding midpoints of edges in $\mathscr{B}$ belong to $\mathscr{B}^{\beta}$, where $\beta:=F^{-1}j\mid Y$.}
\end{Def}
In the sequel, we say that a special endomorphism $j \in V^{'}_p$ is even, resp. odd, if $\nu_p(j^2)$ is even resp. odd. Thus, if $j$ is even, then $Core(j)$ consists of a single superspecial point. If $j$ is odd, then $Core(j)$ is an infinite set.
\begin{Pro} \label{div}
Let $j \in V_p^{'}$ such that $j^2 \neq 0$ and $\nu_p(j^2) \geq 0$. Then $Z(j)$ is a relative divisor over $\Spf W $.
\end{Pro}
\emph{Proof.} The fact that $Z(j)$ is a divisor is proved in 
 \cite{T}, Proposition 4.5. The fact that the local equation of $Z(j)$ is nowhere divisible by $p$ follows from \cite{KR3} Lemma 3.6 and the fact that not the whole supersingular locus of $\mathcal{M}^{HB}$ (which is connected) belongs to $Z(j)$  (see above).
 \qed
\newline

 Let $L \in X^{HB}$ be a superspecial lattice. 
The space  $N_L \subset V_p^{'}$ of special endomorphisms, which map $L$ into itself is a quadratic $\Zp$-modul. (It is the set of all special endomorphisms $j \in V^{'}_p$ such that the $p$-divisible group  corresponding to $L$ belongs to $Z(j)$.) We write $L=L_0\oplus L_1$, and then (by the compatibility condition with $\iota$) for $y\in N_L$ we can write 
$$
y= \begin{pmatrix}
  & y_1 \\ 
 y_0 &  
\end{pmatrix} ; \ \ y_i \in \Hom(L_i, L_{i+1})
$$
A \emph{standard basis} of $L$ is a basis $e_1,e_2,e_3,e_4$ of $L$ such that $e_1,e_2$ is a basis of $L_0$ and such that $e_3=Fe_1$ and $e_4=p^{-1}Fe_2$ and such that for the symplectic form coming from the principal polarization we have $\langle e_1, e_2 \rangle =1$.
After choosing a standard basis of $L$  (see \cite{KR2},  Lemma 5.1 for the existence of a standard basis) we have by  \cite{KR2}, Corollary 5.2
\[
N_L \cong \{ x=
\begin{pmatrix}
  a& b \\ 
 c & pa^{\sigma} 
\end{pmatrix} ; \ \ a,b,c \in \Zpq,\  b^{\sigma}=-b, \ c^{\sigma}=-c\},
\]  
and the quadratic form is given by
\[
Q(x)=paa^{\sigma}-bc.
\]
In this isomorphism, $x$ is the matrix of $y_0$ in the standard basis. We write $\Zpq=\Zp[\delta]/(\delta^2-\Delta)$, where $\Delta \in \Zp^{\times}$ is not a square in $\Zp$.

A $\Zp$-basis of $N_L$ is given by
\[
 s_1=
\begin{pmatrix}
 & -\delta \\ 
 \delta^{-1} 
\end{pmatrix}, \ \
 s_2=
\begin{pmatrix} & \delta \\ 
 \delta^{-1} 
\end{pmatrix}, \ \
 s_3=
\begin{pmatrix}
 1   \\ 
 & p
\end{pmatrix}, \ \
 s_4=
\begin{pmatrix}
 \delta  \\ 
 & -\delta p &
\end{pmatrix}.
\]
The matrix of the bilinear form induced by $Q$  with respect to this basis is:
\[ 
S^{'}=\diag(1,-1,p,-\Delta p).
\]

\begin{Lem}\label{drinsenk}
Let $y_1, y_2 \in V_p^{'}$ such that $y_1^2,y_2^2\neq 0$ and
 let $x \in (Z(y_1)\cap Z(y_2))(\F)$. 
 Then there exists a special endomorphism $y \in V^{'}_p $ such that $y \perp y_1, y_2$ and such that $\nu_p(y^2)=1$ and $x \in Z(y)(\F)$.
\end{Lem}
{\em Proof.} 
We may assume that $y_1$ and $y_2$ are linearly independent.

Suppose that  $x\in  \mathbb{P}_{[\Lambda]} $. 
Let $z\in  \mathbb{P}_{[\Lambda]} $ be superspecial such that $z\notin Core(y_1)\cap Core(y_2) $. It is enough to find $y \in V^{'}_p $ such that $\nu_p(y^2)=1$ and $y \perp y_1,y_2$ and such that $ \mathbb{P}_{[\Lambda]} \subset Z(y)$.
We consider the basis $s_1, ..., s_4$ described above for the quadratic $\Zp$-module of special endomorphisms in $V^{'}_p$ which are endomorphisms of the $p$-divisible group belonging to $z$.
Write $y_1= \sum a_i s_i$ and $y_2=\sum b_i s_i$, where $a_i, b_i \in \Qp$. If $a_4=b_4=0,$ then we can choose $y=s_4$. 
Hence we may assume that $a_4 \neq 0.$ 
Let  $\tilde{y_2}=y_2-(b_4/a_4)y_1=\sum c_i s_i.$ Then  $\tilde{y_2}\neq 0$ and $c_4 =0.$ 
We consider the special endomorphisms
$$
r_1=a_4c_3\cdot s_2 + \frac{1}{p}a_4 c_2 \cdot s_3 + \frac{-1}{\Delta p}(a_2c_3-a_3c_2)\cdot s_4
$$
and
$$
 r_2= -a_4c_3 \cdot s_1 + \frac{1}{p}a_4 c_1 \cdot s_3 + \frac{-1}{\Delta p}(a_1c_3 - a_3 c_1)\cdot s_4. 
$$
One checks that $r_1$ and $r_2$ are both perpendicular to $y_1$ and to $\tilde{y_2}$ and hence also to $y_2$.

Suppose that the coefficients of $s_3$ and $s_4$ in $r_1$ and in $r_2$ do simultaneously vanish. Since $a_4 \neq 0$ it follows that $c_1 =  c_2 =0$. Hence $c_3 \neq 0$ and hence $a_1=a_2=0$.
Hence one of the special endomorphisms $(p+1)s_1+s_2$ and $(p-1)s_1+s_2$  fulfills the requirements of $y$  (for both we have $\nu_p(y^2)=1$  and  $y \perp y_1, y_2,$ and $\mathbb{P}_ {[\Lambda]} $ belongs to the  formal special cycle of one of them).

Suppose now that we have chosen $i$ such that   the coefficients of $s_3$ and $s_4$ in  $r_i$  do not both vanish. 
We choose the integer $m$ minimal such that $y=p^m r_i$ lies in the $\Zp$-span of $s_1,...,s_4$.
Then $p^m a_4 c_3\equiv  0 \mod p$ since otherwise $y^2 \in \Zp^{\times}$, hence $\{z\}=Core(y)$, hence $z\in Core(y_1) \cap Core(y_2)$ which contradicts our assumption on $z$.
Hence the coefficient of $s_3$ or of $s_4$ in $y$ is a unit and hence $\nu_p(y^2)=1$ and $y$ fulfills all requirements above. (Note that both projective lines passing through $z$ belong to $Z(y)$ (see also Lemma \ref{hillem} below), in particular  $\mathbb{P}_{[\Lambda]} \subset Z(y) $.) 
\qed
\newline
 
For $j \in V_p^{'}$ as in Proposition \ref{div} we define $D(j)=Z(j)-Z(j/p)$. This is meant in the following sense. If $Z(j)$ is (locally) given by the equation $f=0$ and $Z(j/p)$ is (also locally) given by $g=0,$ then $D(j)$ is the relative effective divisor locally given by
 $f \cdot g^{-1} =0$. Further, if $\nu_p(j^2)=0$ or $=1$, then $D(j)=Z(j)$.
\begin{Lem} \label{diperg}
Let $y,y_1\in V_p^{'}$ such that  $y_1 \perp y$ and $y^2 \in \Zp \setminus \{0\}$ and $\nu_p(y_1^2)=1$. Let $y^2=\varepsilon p^a$ and $y_1^2=\varepsilon_1 p,$ where $\varepsilon, \varepsilon_1\in \Zp^{\times}.$
\begin{enumerate}
\item If $a \geq 2,$ then $D(y)\cap Z(y_1)$ is a reduced union of projective lines. The same holds if $a=1$ and the image of $-\varepsilon\varepsilon_1$ in $\Fp$ is not a square.
\item If $a=0$, then $Z(y)\cap Z(y_1)$ is a reduced, irreducible horizontal divisor in $Z(y_1)$ (i.e. has no component with support in the special fibre of $Z(y_1)$) which intersects the special fibre of $Z(y_1)$ only in the point of  $Core(y)$. Further, the intersection multiplicity of this horizontal component with each of the two projective lines in the special fibre of $Z(y_1)$ passing through the point of  $Core(y)$ equals $1$. 
\item   If $a=0$ and $y_0\in V_p^{'}$ satisfies $y_0 \perp y, y_1$ and $\nu_p(y_0^2)=1$, then 
$(Z(y_0),Z(y_1),Z(y))=2$.
\item If  $a=0$ and $y_0\in V_p^{'}$ satisfies $y_0 \perp y, y_1$  and $\nu_p(y_0^2)=0$, then 
 $(Z(y_0),Z(y_1),Z(y))=1$. 
\end{enumerate}
\end{Lem}
\emph{Proof.} All claims follow from \cite{T}. More precisely, in the situation of 1., by the equations for antispecial cycles in given in  \cite{T}, section 2, and the connection of antispecial cycles and formal special cycles in our sense discussed in \cite{T}, section 4, 
$D(y)\cap Z(y_1)$ is indeed the reduced union of the projective lines  $\mathbb{P}_{[\Lambda]}$ which belong to the special fibres of $Z(y)$  and $Z(y_1)$.

If $a=0$, then the assertion of 2.  follows by comparing the intersection $Z(y)\cap Z(y_1)$ with the intersection $Z(py)\cap Z(y_1)$. More precisely, again by combining as above the results of sections 2 and 4 of \cite{T},  $Z(py)\cap Z(y_1)$ consists of two (reduced) projective lines in the special fibre and a horizontal component which has all the properties which are claimed for the horizontal component of $Z(y)\cap Z(y_1)$ in the statement of the lemma. Since no projective lines in the special fibre belong to $Z(y)$, it follows that the horizontal components of   $Z(y)\cap Z(y_1)$ and of  $Z(py)\cap Z(y_1)$ are equal.  

Using this, 
in case $\nu_p(y_0^2)=1$, the claim that $$(Z(y_0),Z(y_1),Z(y))=((Z(y_0)\cap Z(y_1)), (Z(y) \cap Z(y_1)))=2$$ (intersection multiplicity in $Z(y_1)$, see \cite{T}, Proposition 4.7 for the first equality) follows from the fact that (by the same reasoning as above) $Z(y_0)\cap Z(y_1)$ is the reduced union of the projective lines  $\mathbb{P}_{[\Lambda]}$ which belong to the special fibres of $Z(y_0)$  and  $Z(y_1)$.

In case $\nu_p(y_0^2)=0$ one has to prove that (as intersection multiplicity in  $Z(y_1)$) on has $((Z(y_0)\cap Z(y_1)), (Z(y) \cap Z(y_1)))=1.$ In other words, this means that the horizontal components of   $(Z(py_0)\cap Z(y_1))$      and of $(Z(py) \cap Z(y_1))$ intersect in $Z(y_1)$ with multiplicity $1$. This in turn follows from the comparison of intersections of (horizontal components of) special cycles in the sense of \cite{KR1} and antispecial cycles resp. formal special cycles in our sense given in \cite{T}, sections 3, 4. 
\qed

\begin{Pro}\label{reg}
Let $j \in V_p^{'}$  be as in Proposition 6, i.e. $j^2\neq 0$ and $\nu_p(j^2)\geq 0$. Then $D(j)$ is a regular formal scheme.
\end{Pro}
\emph{Proof.}
If $\nu_p(j^2)=0,$  then we find $j_0,j_1 \in V^{'}$ such that $Core(j)\subset (Z(j_0) \cap Z(j_1))(\F)$ and  all three endomorphisms are perpendicular to each other and such that $\nu_p(j_0^2)=0$ and  $\nu_p(j_1^2)=1$.
Then by Lemma \ref{diperg} the length (over $W$) of the local ring (which is a local Artin ring) of $Z(j)\cap Z(j_0) \cap Z(j_1)$ is $1$, hence this ring is 
regular. But then $Z(j)$ is also regular in the point of  $Core(j).$

If $\nu_p(j^2)=1,$ then $Z(j)$ is isomorphic to the Drinfeld upper half plane, see \cite{T}, section 4 for the precise statement, which is regular.  

Suppose now that $\nu_p(j^2)\geq 2$. Let $x\in Z(j)(\F)$. We have to show that the local ring $\mathcal{O}_{D(j),x}$ is regular. 

{\bf First case} $x\notin  Z(j/p)(\F)$.

We choose $j_1 \in V_p^{'}$ such that $j_1 \perp j$ and $\nu_p(j_1^2)=1$ and such that $x \in Z(j_1)(\F)$, see Lemma \ref{drinsenk}. Then $Z(j_1)$ is isomorphic to the Drinfeld upper half plane (see \cite{T}, section 4). By Lemma \ref{diperg} we know that locally around $x$, the intersection $Z(j_1) \cap D(j)$ is isomorphic to a projective line over $\F$, hence this intersection is regular in $x$. Hence $D(j)$ is regular in $x$ as well.

{\bf Second case} $x\in  Z(j/p)(\F)$.

If $x$ is not superspecial, then we choose $j_1$ as in case 1 and the reasoning is also the same as above. Suppose now that $x$ is superspecial, $\{x\} = \mathbb{P}_{[\Lambda_0]} \cap \mathbb{P}_{[\Lambda_1]} $, and
 $\mathbb{P}_{[\Lambda_0]} \subset Z(j/p) $.
 We choose $j_1$ linearly independent of $j$ such that $\nu_p(j_1^2)=1$ and  $\mathbb{P}_{[\Lambda_0]} \subset Z(j_1)$ but $\mathbb{P}_{[\Lambda_1]} \not\subset Z(j_1)$.
 We consider the $\Zp$-submodule $\bf{j}$ of $V_p^{'}$ spanned by $j/p$ and $j_1$. Then there is  $j_2 \in \bf{j}$ such that $j_1 \perp j_2$ and $\{j_1, j_2\}$ is a basis of $\bf{j}$. (This can be seen as follows: Consider an arbitrary  quadratic $\Zp$-module $(M,q)$ which is free of finite rank. Then for any element $m\in M$ such that $\nu_p(q(m))$ is minimal, there is an orthogonal  basis of $M$ which contains $m$ as an element.  Thus, if we apply this to $\bf{j}$, and argueing by contradiction,  we see that if our claim were wrong  there would exist some $j^{'} \in \bf{j}$ such that $\nu_p({j^{'}}^2)=0$. But this contradicts the fact that $\mathbb{P}_{[\Lambda_0]} \subset Z(j^{''})$ for all $j^{''}\in \bf{j}$).
 We write $j/p = aj_1 + bj_2$. Then $D(j) \cap Z(j_1)= D(pbj_2)\cap Z(j_1)$ which is regular in $x$ by the same reasoning as above. Again it follows that $D(j)$ is regular in $x$. \qed
\newline

Our next aim is a description of the special fibre of a formal special cycle $Z(j)$, where $j\in V^{'}$ and $j^2\neq 0$.
In the sequel we use the following notations. 
For a vertex $[\Lambda]\in \mathscr{B}$ denote by $d_{[\Lambda]}$ the distance of $[\Lambda]$ to $\mathscr{B}^{\beta}$, where $\beta:=F^{-1}j \mid Y$ as above.
For any formal scheme $X$ over $\Spf W$ we denote by $X_p$ the special fibre of $X$. 
\begin{The}\label{spezf1}
Let $j\in V^{'}$ such that $j^2\neq 0$ and $a := \nu_p(j^2)\geq 0$.
Then the divisor  $Z(j)_p$ in ${\MHB}_p$ is of the form
$$
Z(j)_p=\sum_{[\Lambda]; \  d_{[\Lambda]} \leq \frac{a-1}{2}} (1+p+\ ...\ +p^{\frac{a-1}{2} -d_{[\Lambda]}})\mathbb{P}_{[\Lambda]} \ + \ p^{a/2}\cdot s.
$$
  Here, $s$ is a divisor in ${\MHB}_p$  which is $\neq 0$ if and only if $j$ is even. In this case 
  $s$ meets the supersingular locus only in the point of $Core(j)$.
\end{The}
We will see later (Lemma \ref{horizl}) that (in case that $j$ is even) the intersection multiplicity of $s$ with each of the two projective lines to which the point in $Core(j)$ belongs is $1$. In particular, $s$ is reduced and irreducible.
\newline

\emph{Proof.} 
We know by Proposition \ref{div} that $Z(j)_p$ is a divisor in ${\MHB}_p.$ 
It follows from the description of the underlying reduced scheme of $Z(j)$ given above that the part of $Z(j)_p$ having support in the supersingular locus is of the form $$Z(j)_p^{ss}=\sum_{[\Lambda]; \  d_{[\Lambda]} \leq \frac{a-1}{2}} m_{[\Lambda]}\mathbb{P}_{[\Lambda]}$$ (without knowing the values of the multiplicities $m_{[\Lambda]}$ of the irreducible components $\mathbb{P}_{[\Lambda]}$).
 Further it follows from \cite{KR2}, section 10 that the rest of $Z(j)_p$ is of the form $p^{a/2}\cdot s$ where $s$ is as above. (In the corresponding statement Corollary 10.3 in  loc. cit. it is erroneously asserted that the multiplicity of $s$ is $p^a$ instead of $p^{a/2}$.).

Thus it remains to find the values of $m_{[\Lambda]}$. It suffices to show that 
\begin{enumerate}[(i)]
\item If $d_{[\Lambda]}=\frac{a-1}{2},$ then $m_{[\Lambda]}=1$.
\item If $\mathbb{P}_{[\Lambda]}$ intersects  $\mathbb{P}_{[\Lambda]^{'}}$ and $d_{[\Lambda]} <  d_{[\Lambda^{'}]}\leq \frac{a-1}{2},$ then $m_{[\Lambda]}+m_{[\Lambda^{'}]}=2\cdot(1+p+\ ...\ +p^{\frac{a-1}{2} -d_{[\Lambda{'}]}})+p^{\frac{a-1}{2} -d_{[\Lambda]}}$.
\end{enumerate}

For any point $x\in Z(j)_p(\F)$ we use the following notation. Let $I_x$ be the ideal of $Z(j)_p$ in the local ring  $\mathcal{O}_{\MHB_p,x}$ and let $\mathfrak{m}_x$ be  the maximal ideal of  $\mathcal{O}_{\MHB_p,x}$. Denote  by $n_x$ the maximal integer $n$ with the property that $I_x \subset \mathfrak{m}_x^n$. Denote by $X$ (instead of $X_x$) the $p$-divisible group over $\F$  belonging to the point $x$.

Suppose now we are in the situation of (i) and let  $x \in \mathbb{P}_{[\Lambda]}$ be a superspecial point which does not belong to another  $\mathbb{P}_{[\Lambda^{'}]} \subset Z(j)_p$. Let $f=0$ be the equation of $\mathbb{P}_{[\Lambda]}$ in the local ring
 $\mathcal{O}_{\M_p,x}$.  It follows that $I_x=(f^{m_{[\Lambda]}})$. Since $f$ belongs to a regular parameter system of $\mathfrak{m}$, it follows that $n_x=m_{[\Lambda]}$.

Suppose now that  we are in the situation of (ii) and let  $x$ be the intersection point of $\mathbb{P}_{[\Lambda]}$ and  $\mathbb{P}_{[\Lambda^{'}]}$. Then it follows as above that  $n_x= m_{[\Lambda]}+m_{[\Lambda^{'}]}$.
Further one sees easily that  $j$ as an endomorphism of $X$ is divisible by $p^{\frac{a-1}{2} -d_{[\Lambda]}}$ but not by  $p^{\frac{a-1}{2} -d_{[\Lambda]}+1}$. (Note that $ d_{[\Lambda]}=\frac{a-1}{2}$ means that  $\mathbb{P}_{[\Lambda]}$ belongs to $Z(j)$ but not to $Z(j/p)$.)

Using the notation just introduced it follows that it is enough to show the following:
\begin{enumerate}[(1)]
\item Suppose $x$ is a superspecial point belonging to precisely  one $ \mathbb{P}_{[\Lambda]} \subset Z(j)_p$.  Then $n_x=1$.
\item Suppose we are in the situation of (ii) and $x$ is as above the intersection point of $\mathbb{P}_{[\Lambda]}$ and  $\mathbb{P}_{[\Lambda^{'}]}$.  Suppose further that $j$ is as an endomorphism of $X$ divisible by $p^r$ but not by  $p^{r+1}$. Then $n_x=2\cdot(1+p+\ ...\ +p^{r-1})+p^r.$
\end{enumerate}

For this we are going to use the theory of displays (see \cite{Z2}) in the manner of \cite{KR3}, section 8. 
Denote by $(M,F,V)$ the Dieudonn\'e module to the $p-$divisible group $X$ corresponding to the point $x$ in (1) or in (2). Then the display for $X$ is given by the data $(M, VM, F, V^{-1})$. We find a basis $f_1,...,f_4$ of $M$ for which $F$  and $V$ have matrices
\[
F=
\begin{pmatrix}
&&p& \\
&&&1 \\
1&&& \\
&p&& \\
\end{pmatrix} \sigma, \ \ \ \ \ \ V=
\begin{pmatrix}
&&p& \\
&&&1 \\
1&&& \\
&p&& \\
\end{pmatrix} \sigma^{-1}, Ê 
\] and such that for the alternating  form  coming from the polarization we have $\langle f_1,f_2 \rangle=\langle f_3,f_4 \rangle=1,$ see \cite{KR2}, section 5.
Let $e_1=f_1,\ e_2=f_4,\ e_3=f_2,\ e_4=f_3 $. Denote by  $T$ the $W$-span of $e_1,e_2$ and by $L$ the $W$-span of $e_3,e_4$. Then
\[
M=L\oplus T, \ \ \ VM=L\oplus pT.
\]
Define the matrix $(\alpha_{ij})$ by 
\[
Fe_j=\sum_i \alpha_{ij}e_i \text{ for } j=1,2,  \\
 \]
 \[
 V^{-1}e_j=\sum_i \alpha_{ij}e_i \text{ for } j=3,4.
\]
Hence, 
\[
(\alpha_{ij})=
\begin{pmatrix}
&&&1 \\
&&1& \\
&1&& \\
1&&& \\
\end{pmatrix}.
\]
It follows (see [Z1, p.48]) that the universal deformation of $X$ over $\F[\![ t_1,t_2,t_3,t_4]\!]$ corresponds to the display $(L\oplus T)\otimes W(\F[\![t_1,t_2,t_3,t_4 ]\!]) $ with matrix $(\alpha_{ij})^{\univ}$ (wrt. the basis $e_1,...,e_4$ and with entries in $W(\F[\![t_1,t_2,t_3,t_4 ]\!])$) given by 
\[
(\alpha_{ij})^{\univ}=
\begin{pmatrix}
1&&[t_1]&[t_2] \\
&1&[t_3]&[t_4] \\
&&1& \\
&&&1 \\
\end{pmatrix} \cdot
\begin{pmatrix}
&&&1 \\
&&1& \\
&1&& \\
1&&& \\
\end{pmatrix}=
\begin{pmatrix}
[t_2]&[t_1]&&1 \\
[t_4]&[t_3]&1& \\
&1&& \\
1&&& \\
\end{pmatrix}.
\]
Here the $[t_i]$ denote the Teichm\"uller representatives of the $t_i$.

 Now let $A^{'}=W[\![t_1,t_2,t_3,t_4]\!]$,  let  $R^{'}=\F[\![t_1,t_2,t_3,t_4]\!]$. Denote by 
 $t$ resp. $u$ the image of $t_1$ resp. of $t_4$ in $A^{'}/(t_2,t_3)$  so that   $A^{'}/(t_2,t_3) = W[\![t,u]\!]=:A$.
 Let $R=\F[\![t,u]\!]$.  We extend the Frobenius  $\sigma $ on $W$ to $A^{'}$ resp. $A$ putting $\sigma(t_i)=t_i^p$ resp. $\sigma(t)=t^p$ and $\sigma(u)=u^p$. 
 Further for any $n \in \N$ denote by $\mathfrak{a}_n$ resp. $\mathfrak{r}_n$ the ideal in $A$ resp. in $R$ generated by the monomials $t^iu^{n-i}$, $i=0,...,n$, and let $A_n=A/\mathfrak{a}_n$ and $R_n=R/\mathfrak{r}_n$.  Then $A^{'}$ is a frame for $R^{'}$ resp. $A$ is a frame for $R$ resp. $A_n$ is a frame for $R_n$.  

For an  $A^{'}$ - $R^{'}$-window $(M^{'}, M^{'}_1, \Phi^{'},\Phi^{'}_1)$ let  $M_1^{'^{\sigma}}=A^{'}\otimes_{A^{'}, \sigma}M_1^{'}$ and denote by $\Psi^{'}:  M_1^{'^{\sigma}} \rightarrow M^{'}$ the linearization of $\Phi^{'}_1$. It is an isomorphism of  $A^{'}$-modules.
 Denote by  $\alpha^{'}:M_1^{'} \rightarrow M_1^{'^{\sigma}}$ the composition of the inclusion map $M_1^{'}  \hookrightarrow M^{'} $ followed by $\Psi^{'^{-1}}$. In this way, the category of pairs $(M_1^{'}, \alpha^{'})$ consisting of a free $A^{'}$-module of finite rank and an $A^{'}$-linear injective homomorphism $\alpha^{'}:M_1^{'} \rightarrow M_1^{'^{\sigma}}$ such that Coker  $\alpha^{'}$ is a free $R^{'}$-module becomes equivalent to the category of formal $p$-divisible groups over $R^{'}.$  (Note that the so called nilpotence condition is fulfilled automatically here, comp. also \cite{KR3}, section 8.) 
 A corresponding description holds for the category of formal $p$-divisible groups over $R$ resp. over $R_n$.
 We consider the $A^{'}$ - $R^{'}$ window $(M^{'},M^{'}_1, \Phi^{'},\Phi^{'}_1)$ given by 
\[ 
M^{'}=M\otimes A^{'}, \ \ \ \ M_1^{'}=VM\otimes A^{'}, \ \ \ \ 
\Phi^{'}=
\begin{pmatrix}
t_2&t_1&&p \\
t_4&t_3&p& \\
&1&& \\
1&&& \\
\end{pmatrix} \sigma, \ \ \ \ \Phi^{'}_1=\frac{1}{p}\cdot \Phi^{'},
\]
the matrix of $\Phi^{'}$ being described in the basis $e_1, e_2,e_3,e_4$. The corresponding display is the universal display described above (easy to see using the procedure dscribed on p.2 of \cite{Z2}). Hence $(M^{'},M^{'}_1, \Phi^{'},\Phi^{'}_1)$ is the universal window.
The corresponding matrix of $\alpha^{'}$ wrt. the basis $pe_1,pe_2,e_3,e_4$ of $ M_1^{'}$ resp. the basis $p(1\otimes e_1),p(1\otimes e_2),1\otimes e_3,1\otimes e_4$ of $ M_1^{'^{\sigma}}$ is given by 
\[\alpha^{'}=
\begin{pmatrix}
&&&1 \\
&&1& \\
&p&-t_3&-t_4 \\
p&&-t_1&-t_2 \\
\end{pmatrix}.
\]

The $p$-divisible group $X$ is equipped with a principal polarization $\lambda$  and with a $\Zpq$-action $\iota$.
\begin{Lem} The ideal in $A^{'}$  resp. $R^{'}$ describing the universal deformation of $(X, \lambda, \iota)$ is (in both rings) given by $(t_2,t_3)$.
 The universal object $(M_1,\alpha)$ over $A$ (corresponding to the completion of  $\mathcal{O}_{\MHB_p,x}$) is given by the image of $(M_1^{'}, \alpha^{'})$  over $A$ (i.e. is obtained from  $(M_1^{'}, \alpha^{'})$ by tensoring over $A^{'}$ with $A$).  In the basis $pf_1,f_2,f_3,pf_4$ of $M_1$ resp. the basis $p(1\otimes f_1),1\otimes f_2 ,1\otimes f_3, p(1\otimes f_4)$ of $M_1^{\sigma}$ the map $\alpha $ is then given by the matrix
\[\alpha=
\begin{pmatrix}
&&1& \\
&&-u&p \\
p&-t&& \\
&1&&& \\
\end{pmatrix}=\begin{pmatrix}
&U\\
T&\\
 \end{pmatrix}.
\]

Here
\[  U=
\begin{pmatrix}
 1& \\
-u&p \\ \end{pmatrix} 
 \text{  resp. } \ T=
  \begin{pmatrix}
p&-t \\
&1& \\
 \end{pmatrix}. 
\]

\end{Lem}
 \emph{Proof.}
The principal polarization $\lambda$ of $X$  corresponds to the perfect symplectic pairing on $M$ given by $\langle e_1, e_3 \rangle=1, \langle e_2, e_4 \rangle=-1$ and $\langle e_i, e_j\rangle=0$ for all other pairs $i  \leq j$ (see above). One checks easily that it lifts to a perfect symplectic bilinear form $\langle, \rangle$ of the display over $W(R^{'}/(t_2,t_3))$ induced (by base change) by the universal one above, such that  $\langle , \rangle$ satisfies  the condition $\langle V^{-1}(x), V^{-1}(y) \rangle^V = \langle x, y \rangle $ (where $V$ is the Verschiebung of  $W(R^{'}/(t_2,t_3))$).
Hence $\lambda$ lifts over $R^{'}/(t_2,t_3),$ comp \cite{G}, p. 231.
Further the $\Zpq$-operation $\iota$ also lifts over $R^{'}/(t_2,t_3)$, since on $ M_1^{'}\otimes A^{'}/(t_2,t_3)$ we have  a $\Z/2$-grading
 lifting the $\Z/2$-grading on $VM$ and such that the map induced by $\alpha^{'}$ is homogenous of degree $1$. More precisely, the $\Z/2$-grading is given by $(M_1^{'}\otimes A^{'}/(t_2,t_3))^0= $ span of $\{ pe_1,e_3\}$ and $(M_1^{'}\otimes A^{'}/(t_2,t_3))^1=$ span of $\{ e_2,pe_4\} $. Hence $(X,\lambda, \iota)$ lifts over $R^{'}/(t_2,t_3)$. Since the completion of  $\mathcal{O}_{\MHB_p,x}$ is isomorphic to $R^{'}/(t_2,t_3)$, it follows that the ideal $(t_2,t_3) $ in 
$R^{'}$ is the ideal describing the deformation of $(X, \lambda, \iota)$. 
 Hence the universal object $(M_1,\alpha)$ over $A$ corresponding to the completion of  $\mathcal{O}_{\MHB_p,x}$ is given by the image of $(M_1^{'}, \alpha^{'})$  over $A$.  It follows that in the basis $pf_1,f_2,f_3,pf_4$ of $M_1$ resp. the basis $p(1\otimes f_1),1\otimes f_2 ,1\otimes f_3, p(1\otimes f_4)$ of $M_1^{\sigma}$ the map $\alpha $ is then given by the matrix
\[\alpha=
\begin{pmatrix}
&&1& \\
&&-u&p \\
p&-t&& \\
&1&&& \\
\end{pmatrix}=\begin{pmatrix}
&U\\
T&\\
 \end{pmatrix}
\]

as asserted.
\qed
\newline

Next we want to determine the maximal $n$ such that $j$ lifts over $R_n$. The  $\Z/2$-grading of the Dieudonn\'e module $M$ of $X$ is given by $M^0=$ span of $ \{f_1,f_2\} $ and $M^1=$  span of $\{ f_3,f_4 \}$. The matrix of $j$ in the basis $f_1,..., f_4$ is of the form 
\[j=
\begin{pmatrix}
&&pa^{\sigma}&-b\\
 && -c&a \\
  a& b &&\\ 
 c & pa^{\sigma} &&\\
\end{pmatrix} ; \ \ a,b,c \in \Zpq,\  b^{\sigma}=-b, \ c^{\sigma}=-c,
\]
see \cite{KR2}, section 5. Write $j=p^r \bar{j}$, where  
\[\bar{j}=
\begin{pmatrix}
&&pa_0^{\sigma}&-b_0\\
 && -c_0&a_0 \\
  a_0& b_0 &&\\ 
 c_0 & pa_0^{\sigma} &&\\
\end{pmatrix} ; \ \ a_0,b_0,c_0 \in \Zpq \text{ are not all divisible by } p.
\]
\begin{Lem} \label{hilflem1}
If $b_0\equiv 0 \mod p$ and $c_0 \equiv 0 \mod p,$ then
$x$ is the intersection point of $\mathbb{P}_{[\Lambda]}$ with another projective line  $\mathbb{P}_{[\Lambda^{'}]} \subset Z(j)$  and  $d_{[\Lambda^{'}]}=d_{[\Lambda]}=0$. 
\end{Lem}
\emph{Proof.}
Since in this case $\nu_p(\bar{j}^2)=1$ we can identify $Z(\bar{j})$ with the Drinfeld upper half plane (see \cite{T}, chapter 4) and therefore the statement  just means that under this identification $x$ is a superspecial point in the Drinfeld upper half plane (in the sense of \cite{KR1}).
 Denoting the Dieudonn\'e module of $x$ with its $\Z/2$-grading as above by $M=M^0 \oplus M^1$, we have to show that $\bar{j}M^0=VM^0$ and $\bar{j}M^1=VM^1$.
Now $\bar{j}M^0$ is the $W$-span of $a_0f_3+c_0f_4$ and $b_0f_3+pa_0f_4$. Since $a_0$ is a unit in $W$ and $b_0$ and $c_0$ are divisible by $p$ this is the same as the $W$-span of $f_3$ and  $b_0f_3+pa_0f_4$, and this is the same as the $W$-span of $f_3$ and $pf_4$ which is  $VM^0$. The equality $\bar{j}M^1=VM^1$ is proved in the same way. \qed
\newline

It follows that neither in the situation of (1) nor in the situation of (2) the assumption made in the lemma is fulfilled. 
Write now $(M_1(n), \alpha(n))$ for $M_1\otimes A_n$ and the map $\alpha(n): M_1(n)  \rightarrow M_1(n)^{\sigma}$ induced by $\alpha$.  Write $j(1)$ for the endomorphism of $M_1(1)$ induced by $j$. In the basis $pf_1,f_2,f_3,pf_4$ it is given by the matrix 
\[j(1)=
\begin{pmatrix}
&&a^{\sigma}&-b\\
 && -c&pa \\
  pa& b &&\\ 
 c & a^{\sigma} &&\\
\end{pmatrix} = \begin{pmatrix}
&j_1(1)\\
j_0(1)&\\
 \end{pmatrix}.
\]
Here
\[  j_0(1)=
\begin{pmatrix}
 pa& b \\ 
 c & a^{\sigma} \\
 \end{pmatrix} 
 \text{  resp. } \ j_1(1)=
  \begin{pmatrix}
a^{\sigma}&-b\\
 -c&pa \\
 \end{pmatrix}. 
\]
In order to lift $j$ over $R_n$ we need an endomorphism $j(n)$ of $M_1(n)$ lifting $j(1)$ such that  the following diagram commutes:
\[
\xymatrix{ M_1(n) \ar[d]_{j(n)} \ar[r]^{\alpha(n)}&  M_1(n)^{\sigma} 
\ar[d]^{\sigma(j(n))}  \\M_1(n)\ar[r]_{\alpha(n)} & M_1(n)^{\sigma}  .}
\]
In other words we are looking for liftings $j_0(n)$ of $j_0(1)$ and  $j_1(n)$ of $j_1(1)$ over $A_n$ such that \[ Uj_0(n)= \sigma (j_1(n))T \  \text{ and } \   Tj_1(n)= \sigma (j_0(n))U. \tag{*}\]
Suppose $n=p^l$, where $l\geq1$, and suppose we have found liftings  $j_0(p^{l-1})$ and  $j_1(p^{l-1})$  satifying (*). For any choice of liftings   $j_0(p^{l})$ and  $j_1(p^{l})$  of $j_0(p^{l-1})$ and  $j_1(p^{l-1})$ 
the matrices $\sigma(j_i(p^{l}))$ are equal  to $\sigma(j_i(p^{l-1})$ interpreted as matrix over $A_{p^l}$. Hence there are liftings  $j_0(p^{l})$ and  $j_1(p^{l})$  satifying (*) if and only if 
\[ U^{-1}\sigma (j_1(p^{l-1}))T \  \text{ and } \  T^{-1} \sigma (j_0(p^{l-1}))U\] are integral and in this case 
\[j_0(p^l)=U^{-1}\sigma (j_1(p^{l-1}))T \  \text{ and } \ j_1(p^l)= T^{-1} \sigma (j_0(p^{l-1}))U.\] 
Define now inductively matrices $X(l)$ and $Y(l)$ over $A_{p^l}\otimes_{\Z}\Q$ as follows: $X(0)=j_0(1)$ and $Y(0)=j_1(1)$ and 
\[ X(l+1)= U^{-1}\sigma (Y(l))T  \  \text{ and } \  Y(l+1)=T^{-1} \sigma (X(l))U.\] (Again $\sigma (Y(l))$ and $\sigma (X(l))$ are well defined over $A_{p^{l+1}}\otimes_{\Z}\Q$.)

The following lemma can easily be proved by induction.
\begin{Lem} \label{hilflem2}
Let $l \geq 1$. 
The matrix $X(l)$ is of the form
\[ X(l)=
\frac{1}{p^l} \left( \begin{pmatrix}
0&0 \\
0&  (tu)^{1+p+...+p^{l-2}}\cdot (   (-1)^l \sigma^{l+1}(a) \cdot (tu)^{p^{l-1}}   +    b \cdot u^{p^{l-1}}-   c \cdot t^{p^{l-1}} ) \\
\end{pmatrix}+p\cdot A(l)\right), 
\] where $A(l)$ has entries in $A_{p^l}$.  Similarly,  $Y(l)$ is of the form
\[ Y(l)=
\frac{1}{p^l} \left( \begin{pmatrix}
 (tu)^{1+p+...+p^{l-2}} \cdot (   (-1)^l   \sigma^{l+1}(a) \cdot (tu)^{p^{l-1}}   +     b \cdot u^{p^{l-1}}- c \cdot t^{p^{l-1}} )&0 \\
0&0 \\
\end{pmatrix}+p\cdot B(l)\right), 
\] where $B(l)$ has entries in $A_{p^l}$. \qed
\end{Lem}

Lemma \ref{hilflem2} shows together with the discussion above that $j$ lifts over $R_{p^r}$ but not over $R_{p^{r+1}}$. More precisely, if we are in the situation of (1) or (2), then $b$ and $c$ are not both divisible by $p^{r+1}$ (Lemma \ref{hilflem1}) and hence we see from Lemma  \ref{hilflem2} (applied in case $l=r+1$) that $j$ lifts over $R_{2\cdot(1+p+...+p^{r-1})+p^r}$ but not over  $R_{2\cdot(1+p+...+p^{r-1})+p^r+1}$. Hence $n_x=2\cdot(1+p+...+p^{r-1})+p^r$ as claimed in (1) resp. in (2). \qed  
\newline

From Theorem \ref{spezf1} we immediately obtain the following corollary (using the same notations and assumptions).

\begin{Cor}\label{spezf}
If $a=0$, then  $D(j)_p$ is of the form $D(j)_p=s$. If $a>0,$ then $D(j)_p$ is of the form  
$$
D(j)_p=\sum_{[\Lambda]; \  d_{[\Lambda]} \leq \frac{a-1}{2}} p^{\frac{a-1}{2} -d_{[\Lambda]}}\mathbb{P}_{[\Lambda]} \ + \ p^{a/2-1}(p-1)\cdot s.
$$
\end{Cor} \qed

\begin{Lem} \label{pmultlem}
Let  $j \in V^{'}$ such that $j^2 \in \Zlp \setminus \{0\}$,  suppose that 
let
  $x \in Z(j)(\F)$  and let $D\subset \mathcal{M}^{HB}$ be  a divisor which is regular in $x$, such that locally around $x$ we have  $D_p \subset Z(j)_p$. 
  Let  $R= \mathcal{O}_{D,x}$ and let $(f)$ be the ideal of $Z(j)\cap D$ in $R$ (i.e., $\mathcal{O}_{Z(j)\cap D,x}=R/(f)$).
   Then  the ideal of $Z(pj)\cap D$ in $R$ is of the form $(h\cdot p \cdot f)$, where $h$ is coprime to $p$ in $R$. 

\end{Lem}

\emph{Proof.} 
We may assume that $f \neq 0.$ 

The canonical map
\[
\varrho_0: \Spf R/(f) \rightarrow D
\]
factors  through $Z(j)\cap D$. 

\begin{bf} Claim: \end{bf} 
\emph{The canonical map}
\[
\varrho_{1}:\Spf R/(pf) \rightarrow D
\]
\emph{factors through  $Z(pj)\cap D$.} 

By our assumption  $D_p \subset Z(j)_p$ (locally around $x$) we know that $f$ is divisible by $p$.
We have $R/(f)=(R/(pf))/I$, where $I=(f)/(pf)$. Since $f$ is divisible by $p$, the ideal  $I$ carries a nilpotent  $pd$-structure.  Hence we may apply   Grothendieck-Messing theory for the pair $R/(pf)$, $R/(f)$.
We denote by  $M$  the value of the crystal of the  $R/(pf)$-valued point $\varrho_1$ in $R/(pf)$, and by $\overline{M}$ the value of the crystal of the $R/(f)$-valued point $\varrho_0$ in $R/(f)$. Then $\overline{M}=M\otimes R/(f)$. Denote by  ${\mathcal{F}}\hookrightarrow M$ the Hodge filtration of  $\varrho_{1}$ and by $\overline{\mathcal{F}}\hookrightarrow \overline{M}$ the Hodge filtration of $\varrho_{0}$. Then the Hodge filtration of $\varrho_{1}$ lifts the  Hodge filtration of $\varrho_{0}$. We have $j\overline{\mathcal{F}} \subset \overline{\mathcal{F}}$, and we have to show that $pj{\mathcal{F}} \subset {\mathcal{F}}$.  As in the proof of Proposition 4.5 in \cite{T}  we find a basis $e_1,e_2,e_3,e_4$ of $M$, such that, 
if $\overline{e_i}$ denote the images of the $e_i$ in $\overline{M}$, then a basis of $\overline{\mathcal{F}}$ is given by  $\overline{e_2}, \overline{e_3}$, and such that, for suitable  $m_1,m_4 \in I,$ a basis of $\mathcal{F}$  is given by $f_0=e_2+m_1e_1$ and $f_1=e_3+m_4e_4$. 
Since $j\overline{\mathcal{F}}\subset \overline{\mathcal{F}}$ it follows that $jf_i \in \langle e_2, e_3 \rangle + IM$ for $i=1,2$. Hence $pjf_i \in \langle pe_2, pe_3 \rangle + pIM=\langle pe_2, pe_3 \rangle = \langle pf_0, pf_1 \rangle \subset \mathcal{F}$.
This proves the claim. 
\newline
Let $\mathcal{O}_{Z(pj)\cap D,x}=R/(g)$. Since $D$ is regular in $x$ we know that $R$ is a UFD. We write $g=h\cdot s \cdot f$, where $h$ is prime to $p$ and $V(s)$ has support in the special fibre of $D$. 
We must show that  $s=p$. We already know that $s$ is divisible by $p$. 
Suppose that $q$ is prime divisor of $p$ in $R$ such that $g$ divisible by $qpf$. 
Let $C=R/(qpf)$, let $\overline{C}=R/(qf)$, and let $\overline{\overline{C}}=R/(f)$. This is a chain of nilpotent $pd$-extensions.

We consider the canonical $C$-, resp. $\overline{C}$ -, resp. $\overline{\overline{C}}$- valued points 
\[
\varrho_C: \Spf C \rightarrow D, \ \ \ 
\varrho_{\overline{C}}: \Spf \overline{C} \rightarrow D, \ \ \ 
\varrho_{\overline{\overline{C}}}: \Spf \overline{\overline{C}} \rightarrow D.
\]

\begin{bf} Claim: \end{bf} {\em$\varrho_{\overline{C}}$ factors through $Z(j)$.}

We denote by    $M$ the value of the crystal of the  $\varrho_C$-valued point  in $C$,  analogously we define   $\overline{M}$ and  $\overline{\overline{M}}$. 
Further let ${\mathcal{F}}\hookrightarrow M$ by the Hodge filtration of $\varrho_{C}$ and analogously we define $\overline{\mathcal{F}}\hookrightarrow \overline{M}$ and $\overline{\overline{\mathcal{F}}}\hookrightarrow \overline{\overline{M}}$. As in the proof of Proposition 4.5 in \cite{T}  we find a basis $e_1,e_2,e_3,e_4$ of $M$, such that, 
if $\overline{e_i}$ denote the images of the $e_i$ in $\overline{M}$, then a basis of $\overline{\mathcal{F}}$ is given by  $\overline{e_2}, \overline{e_3}$, and such that for suitable  $m_1,m_4 \in (qf)/(qpf)$ a basis of  $\mathcal{F}$ is given by $f_0=e_2+m_1e_1$ and $f_1=e_3+m_4e_4$.
The  endomorphism $j$ induces an endomorphism of the $C$-module $M$, and we have $pj\mathcal{F}\subset\mathcal{F}$.  To show the above claim we must show that  $j\overline{\mathcal{F}}\subset \overline{\mathcal{F}}$. 
Since $pjf_0\in \mathcal{F}$ and by the explicit form of $f_0$ and $f_1$ we can write 
 $jf_0 = u +xe_4$ for a suitable  $u\in \mathcal{F}$ and some $x\in C$ with the property that $px=0$.
 Let  $\hat{x}\in R$ be a lift of $x$. Then $p\hat{x}\in (qpf)$, hence $\hat{x}=yqf$ for a suitable $y\in R$. 
  hence $x \in (qf)/(qpf)$,  hence $j\overline{e_2}\in \overline{\mathcal{F}}$. In the same way one shows that $j\overline{e_3}\in \overline{\mathcal{F}}$. From this the claim follows.

The claim implies the assertion, since it implies that $f$ is divisible by $qf$, a contradiction.

 \qed

 We reformulate the lemma in a more geometrical manner in the following proposition.
 \begin{Pro} \label{pmult}
 In the situation of Lemma \ref{pmultlem} and under the assumption that $f\neq 0,$ the divisor $D(pj)\cap D$ in $D$ is locally around $x$ of the form
 $$D(pj)\cap D=D_p+h,$$ where $h$ is an effective divisor which is horizontal, i.e. it has no component with support in $D_p$.
 \end{Pro} \qed 

\section{On the structure of intersections of formal special cycles} 
By the proof of \cite{T}, Proposition 4.5, for any $j \in V_p^{'}$,  and any  $x\in \MHB(\F)$ the ideal of $Z(j)$  in $\mathcal{O}_{\MHB, x}$ is generated by one element (we do not use the assumptions $j^2 \neq 0$ and $\nu_p(j^2)\geq 0$  made in Proposition \ref{div}).
\begin{Lem} \label{Lweg}
Let $j, j^{'} \in V^{'}$ be linearly independent. Suppose further that for any  $x\in \MHB(\F)$  the equation in $x$ (that is a generator of the ideal in  $\mathcal{O}_{\MHB, x}$)  of at least one of the corresponding two formal special cycles  is not divisible by $p$. Then $$\mathcal{O}_{Z(j)}\otimes^{\mathbb{L}}\mathcal{O}_{Z(j^{'})}= \mathcal{O}_{Z(j)}\otimes \mathcal{O}_{Z(j^{'}).}$$ More precisely, the object on the left hand side is represented in the derived category by the object on the right hand side.
The same formula holds if $Z(j)$ or $Z(j^{'})$ or both are replaced by $D(j)$ resp. $D(j^{'})$.
\end{Lem}
\emph{Proof.} Let $x\in \MHB(\F)$ and let $R=\mathcal{O}_{\MHB, x}$. Let the ideals of $Z(j)$ and $Z(j^{'})$ in $R$ be generated by  $f$ and $f^{'}$.  Suppose that $f$ is not divisible by $p$. We consider the  exact sequence
\begin{equation*}
\begin{CD}
0 @ >>> R @>f\cdot>>R @>>>
R/(f)@>>> 0.
\end{CD}
\end{equation*}
Tensoring this with $R/(f^{'})$ we see that $$
\mathcal{TOR}_1(\mathcal{O}_{Z(j)}, \mathcal{O}_{Z(j^{'})})_x=\ker(R/(f^{'})\stackrel{f \cdot}{\longrightarrow}R/(f^{'})).
$$ To show that this vanishes
 we have to show that $f$ and $f^{'}$ have no common divisor in the regular ring $R$. Let $g= \gcd(f,f^{'})$. We have to show that $g$ is a unit. It follows from our assumption in the lemma that $g$ is coprime to $p$. 

Let now $T_x$ be the matrix of the bilinear form associated to $Q$ for the basis $j, j^{'}$ of $\bf{j}$ $= \Zlp$-span of $j,j^{'}$ in $V^{'}$ and let $t=j^2$ and $t^{'}={j^{'}}^2$. Since $j, j^{'} \in V^{'}$, we find special cycles $Z(t, \omega)$ and $Z(t^{'}, \omega^{'})$ and an  $\F$-valued point $y$ of $Z(t, \omega)\times_{\M}Z(t^{'}, \omega^{'})$ such that the corresponding $p$-divisible group resp. its corresponding special endomorphisms are the $p$-divisible group of $x$ resp. $j,j^{'}$ (see also the incidence relation under Proposition \ref{locglob}).  It follows that for the completions $\hat{\mathcal{O}}_{Z(t, \omega) \times_{\M} Z(t^{'}, \omega^{'}),y}\otimes_{\Zp}W\cong \hat{R}/(f,f^{'})$ (we denote the images of $f,f^{'}$ and $g$ in $\hat{R}$ by the same letters). Suppose that $g$ is not a unit. Then  $R/(g)\otimes_W \Q $ has dimension $1$. Thus $\mathcal{O}_{Z(t, \omega) \times_{\M} Z(t^{'}, \omega^{'}),y}\otimes_{\Zlp} \Q$ has dimension $1$.
  It follows that there is some $\C$-valued point  $\xi$ of the generic fibre of $  Z(t, \omega) \times_{\M} Z(t^{'}, \omega^{'})$ such that the corresponding $2 \times 2$ fundamental matrix in $\xi$ (defined as above) equals $T_x$. But then the proof of  Proposition 1.4 of \cite{KR2} shows that $T_x$ is positive definite.  Then again by the same proposition, the generic fibre of the part of  $Z(t, \omega) \times_{\M} Z(t^{'}, \omega^{'})$ where the fundamental matrix equals $T_x$ has dimension $0$ which is a contradiction to the fact that $\dim(\mathcal{O}_{Z(t, \omega) \times_{\M} Z(t^{'}, \omega^{'}),y}\otimes_{\Zlp} \Q)=1$. The proof in case that $Z(j)$ or $Z(j^{'})$ or both are replaced by $D(j)$ resp. $D(j^{'})$ is the same. \qed

\begin{Pro} \label{diag}
Let $y_1,y_2,y_3\in V^{'}$ such that 
  the corresponding fundamental matrix is nonsingular and has entries in $\Zlp$. Then the derived tensor product $ \mathcal{O}_{Z(y_1)}\otimes^{\mathbb{L}}\mathcal{O}_{Z(y_2)}\otimes^{\mathbb{L}}\mathcal{O}_{Z(y_3)}$ depends only on the $\Zlp$-span $\bf{y}$ of $y_1,y_2,y_3$ in $V^{'}$.
\end{Pro}

\emph{Proof.} First we observe that for any basis $y_1^{'}, y_2^{'}, y_3^{'}$ of $\bf{y}$ the derived tensor product $ \mathcal{O}_{Z(y_1^{'})}\otimes^{\mathbb{L}}\mathcal{O}_{Z(y_2^{'})}\otimes^{\mathbb{L}}\mathcal{O}_{Z(y_3^{'})}$ is invariant under any permutation of the $y_i^{'}$.

{\bf Claim} $ \mathcal{O}_{Z(y_1^{'})}\otimes^{\mathbb{L}}\mathcal{O}_{Z(y_2^{'})}\otimes^{\mathbb{L}}\mathcal{O}_{Z(y_3^{'})}$  \emph{ is invariant if one $y_i$ is replaced by $\varepsilon y_i + zy_l$, where $\varepsilon \in \Zlp^{\times}, \ z \in \Zlp$ and $l \neq i$. }

We want to apply the lemma setting $j=y_i, j^{'}=y_l$. We must show that in no $\F$-valued point of $\MHB$ the equation of $Z(y_i)$ and $Z(y_l)$ is divisible by $p$. Suppose that there is some $\F$-valued point of $\MHB$ in which both equations are divisible by $p$. Then the same is true for any linear combination of $y_i$ and $y_l$. Further it follows from Proposition \ref{div} that $y_i^2=0$ and $y_l^2=0$ and hence again the same is true for any linear combination of $y_i$ and $y_l$. Thus it follows (by diagonalizing) that  the fundamental matrix of $y_i$  and $y_l$ is the zero matrix. But then it follows that the fundamental matrix of $y_1,y_2,y_3$ is singular which contradicts our assumption. 
Thus we may apply  the lemma, by which $\mathcal{O}_{Z(y_i^{'})}\otimes^{\mathbb{L}}\mathcal{O}_{Z(y_l^{'})}=\mathcal{O}_{Z(y_i^{'})}\otimes \mathcal{O}_{Z(y_l^{'})}$ $=\mathcal{O}_{Z(y_i^{'}) \cap Z(y_l^{'})}$ which only depends on the $\Zlp$-span of $y_i, y_l$ in $V^{'}$.
From this the claim follows.

Since we can transform the basis $y_1,y_2,y_3$ by a suitable sequence of  permutations and operations as in the claim into any other basis of $\bf{y}$, the claim of the proposition follows.
\qed
\newline

In particular, the proposition shows that for the calculation of $(Z(j_1),Z(j_2), Z(j_3))$ we may assume  our three fixed special endomorphisms   $j_1,j_2,j_3$ are perpendicular to each other, i.e. we may assume that $T$ is a diagonal matrix.  
\newline
Suppose for the moment that our fixed $T$ is not divisible by $p$. Then, as shown in \cite{KR2}, section 6,  the intersection $Z(j_1)\cap Z(j_2) \cap Z(j_3)$ is the formal spectrum of a local Artin ring whose length (over $W$) is calculated in loc. cit. We want to show that this length equals the intersection multiplicity in the sense defined above. To this end we show the following
\begin{Pro}
If $T$ (as matrix over $\Zlp$) is not divisible by $p,$ then
$$\mathcal{O}_{Z(j_1)}\otimes^{\mathbb{L}}\mathcal{O}_{Z(j_2)}\otimes^{\mathbb{L}} \mathcal{O}_{Z(j_3)}=\mathcal{O}_{Z(j_1)}\otimes \mathcal{O}_{Z(j_2)}\otimes \mathcal{O}_{Z(j_3)}$$
(meaning that the right hand side represents the left hand side in the derived category). In particular,
$$(Z(j_1),Z(j_2),Z(j_3))=\lg_W(\mathcal{O}_{Z(j_1)\cap Z(j_2) \cap Z(j_3),x}),$$ where $x$ is the unique $\F$-valued point of $Z(j_1)\cap Z(j_2) \cap Z(j_3)$.
\end{Pro}
\emph{Proof.} By Proposition \ref{diag} we may assume that $j_1,j_2,j_3$ are pairwise perpendicular to each other. We further assume $\nu_p(j_1^2) \leq \nu_p(j_2^2) \leq  \nu_p(j_3^2)$. Then $\nu_p(j_1^2)=0$.
By Lemma \ref{Lweg} we have $\mathcal{O}_{Z(j_1)}\otimes^{\mathbb{L}}\mathcal{O}_{Z(j_2)}=\mathcal{O}_{Z(j_1)}\otimes \mathcal{O}_{Z(j_2)}$ in the above sense.  Let the ideal of $Z(j_i)$ in  $ \mathcal{O}_{\MHB,x}$ be generated by $f_i$ and let $R=\mathcal{O}_{\MHB,x}/(f_1)$. Denote the images of $f_2$ resp. $f_3$ in $R$ by $\overline{f}_2$ resp. $\overline{f}_3$. The same reasoning as in the proof of Lemma \ref{Lweg} shows that the claim follows if we know that $\overline{f}_2$ and $\overline{f}_3$ are coprime in $R$ (which is a UFD by Proposition \ref{reg}). But this in turn follows from the fact that $\dim(R/(\overline{f}_2,\overline{f}_3))=0$. 
\qed
\newline

Our next aim is to prove the following multilinearity property of the intersection product $(Z(j_1),Z(j_2),Z(j_3))$.

\begin{Pro}\label{multlin}
Suppose that at least one of  $j_1,j_2,j_3$ is odd. Then 
$$
(Z(j_1),Z(j_2),Z(j_3))=\sum_{l,m,n}(D(j_1/p^l),D(j_2/p^m),D(j_3/p^n)),
$$
where the sum is taken over all possible triples $(l,m,n)$ (i.e. setting $a_i=\nu_p(j_i^2)$, we have $l \leq [a_1/2],\  m \leq [a_2/2], \ n \leq [a_3/2]$, where $[ \ ]$ denotes Gauss brackets).
\end{Pro}
Note that in case that $j_1,j_2,j_3$ are pairwise perpendicular to each other, the hypothesis of the proposition is fulfilled as we see from \cite{Ku1}, 1.16, see also (\ref{iden}).
Note also that using the fact that $Z(j_1)\cap Z(j_2) \cap Z(j_3)$ has support in the supersingular locus (follows from \cite{KR2}, Proposition 3.8) and  the fact that this support ist proper over $\F$ (follows for example from the beginning of section 3 and Lemma \ref{hillem} or from \cite{KR2}, section 8) one sees that all of these intersection multiplicities  are finite. The same holds for the other   intersection multiplicities in the subsequent proof.)

\emph{Proof.} 
Suppose $j_1$ is odd and $\nu_p(j_1^2)=2r+1$. Denote by $\mathcal{I}$ the ideal sheaf of $D:=D(j_1/p^r)(=Z(j_1/p^r))$ in $\mathcal{O}_{\MHB}$. Denote by $\mathcal{J}$ the ideal sheaf of $\Delta:=Z(j_1)-D(j_1/p^r)$ (the notation in the latter expression is meant in the same sense as the notation $D(j)=Z(j)-Z(j/p)$). Our first aim is to show that  $$(Z(j_1),Z(j_2),Z(j_3))=(\Delta,Z(j_2), Z(j_3))+(D,Z(j_2),Z(j_3)).$$ 

We consider the canonical short exact sequences 
$$
\begin{CD}
0 @ >>> \mathcal{J}/(\mathcal{J}\cdot \mathcal{I}) @>>> \mathcal{O}_{\MHB}/(\mathcal{J}\cdot \mathcal{I}) @>>>\mathcal{O}_{\MHB}/\mathcal{J}@>>> 0
\end{CD}
$$
and 
$$
\begin{CD}
0 @ >>> (\mathcal{J}+\mathcal{I})/ \mathcal{I}  @>>> \mathcal{O}_{\MHB}/ \mathcal{I} @>>>\mathcal{O}_{\MHB}/(\mathcal{J}+\mathcal{I})@>>> 0.
\end{CD}
$$
\begin{Lem}
The inclusion $\mathcal{J}\cdot \mathcal{I} \hookrightarrow \mathcal{J}\cap \mathcal{I}$ is an equality.
\end{Lem}
\emph{Proof.} We can check this locally. Let $x \in \MHB (\F)$ and let $R= \mathcal{O}_{ \MHB, x}$. For any $l\leq r$ let $(f_l)$ be the ideal of $D(j/p^l)$ in $R$. By Proposition \ref{reg}, each $f_l$ is a prime element in $R$. It follows from Corollary \ref{spezf}  that the $f_l$ are pairwise distinct (more precisely, Corollary \ref{spezf} shows   that even the special fibres of the several $\Spf R/(f_l)$ are pairwise distinct). Since $\mathcal{I}_x$ equals $(f_r)$ and $\mathcal{J}_x$ equals $\prod_{l\neq r}(f_l)$,  the claim of the lemma follows.
\qed
\newline

By the first exact sequence have $$(Z(j_1),Z(j_2),Z(j_3))=(\Delta,Z(j_2), Z(j_3))+ \chi(\mathcal{J}/(\mathcal{J}\cdot \mathcal{I})\otimes^{\mathbb{L}}\mathcal{O}_{Z(j_2)}\otimes^{\mathbb{L}}\mathcal{O}_{Z(j_3)}).$$
 
Now using the lemma we see that $ (\mathcal{J}+\mathcal{I})/ \mathcal{I}=\mathcal{J}/(\mathcal{J}\cap \mathcal{I})=\mathcal{J}/(\mathcal{J}\cdot \mathcal{I})$. This shows together with the second exact sequence that
\begin{equation*}
 \begin{split} \chi(\mathcal{J}/(\mathcal{J}\cdot \mathcal{I})\otimes^{\mathbb{L}}\mathcal{O}_{Z(j_2)}\otimes^{\mathbb{L}}\mathcal{O}_{Z(j_3)})=  (D,Z(j_2),Z(j_3)) 
-(D\cap \Delta,Z(j_2),Z(j_3)). 
\end{split} 
\end{equation*}

Thus in order to show that $$(Z(j_1),Z(j_2),Z(j_3))=(\Delta,Z(j_2), Z(j_3))+(D,Z(j_2),Z(j_3))$$ it remains to show the 

{\bf Claim} $(D\cap \Delta,Z(j_2),Z(j_3))=0$. 

Using  Lemma \ref{Lweg} we see that
$$
\chi(\mathcal{O}_{D\cap \Delta} \otimes^{\mathbb{L}}\mathcal{O}_{Z(j_2)}\otimes^{\mathbb{L}}\mathcal{O}_{Z(j_3)})=\chi(\mathcal{O}_{D\cap \Delta} \otimes^{\mathbb{L}}\mathcal{O}_{Z(j_2) \cap Z(j_3)}).
$$
First we observe that by the proof of the above lemma, $\Delta \cap D$ is as a divisor in $D$.
The same reasoning as in the proof  of Proposition \ref{unglmult} below (second case in the proof) shows that $\Delta \cap D$ is as a divisor in $D$ of the form $rD_p+h$, where 
$h$ is a divisor of the form $\sum_x h_x$ where the sum runs over a discrete set of $\F$-valued points of $D$  and $h_x$ is a horizontal divisor (meaning that its equation is coprime to $p$) meeting the underlying reduced subscheme of $D$ only in $x$. (More precisely, $D$ plays the role of $D(y_1)$  in the situation of  the proof  of Proposition \ref{unglmult}
 and (for any $l<r$), $D(y_1/p^l)$ plays the role of $D(y_2)$ in that proof. Note that $\Delta = \sum_{l<r}D(y_1/p^l).$)  Now we regard $rD_p$ and $h$ as (formal) closed subschemes of $\MHB.$
Suppose that $h=0.$ 
Then $(D\cap \Delta,Z(j_2),Z(j_3))=(rD_p, Z(j_2)\cap Z(j_3)).$ This in turn equals $0$ as  follows from the exact sequence
$$
\begin{CD}
0 @ >>> \mathcal{O}_D @> p^{r}\cdot>> \mathcal{O}_D @>>>\mathcal{O}_D/(p^{r})@>>> 0.
\end{CD}
$$
Thus the claim is proved if we can show that $h=0.$  Suppose that there is $x \in (D\cap \Delta)(\F)$ such that $h_x \neq 0.$
Using Lemma \ref{drinsenk} we choose a special endomorphism $j$ such that $\nu_p(j^2)=1$ and  $x \in Z(j)(\F)$ and $j\perp j_1.$ We write $j_1^2=\varepsilon_1p^{2r+1}$ and $j^2=\varepsilon p.$ It is easy to see that we can choose $j$ such that the image of $-\varepsilon \varepsilon_1$ in $\Fp$ is a square. Then it follows from the results of \cite{T} (in the same way as the proof of Lemma \ref{diperg}) that $Z(j) \cap \Delta$ has support in the special fibre, that this support is proper over $\F$ and that $Z(j)\cap D$ has precisely two horizontal components which meet the special fibre in the same projective line. Further their intersection points with this projective line are both  neither intersection points of two projective lines which belong to $Z(j)$ nor intersection points of two projective lines which belong to $D$. Further the intersection multiplicity of each of these horizontal components with the projective line which they meet is $1.$
We consider the intersection multiplicity $(D, \Delta, Z(j)).$ It can be written as $((Z(j) \cap D),(Z(j)\cap \Delta))$ (intersection multiplicity in $Z(j)$, compare also \cite{T}, section 4). Let $(Z(j) \cap D)_v$  resp. $(Z(j) \cap D)_h$ denote the vertical resp. horizontal part of $(Z(j) \cap D).$  Now  Lemma \ref{mitte0} below shows that $((Z(j) \cap D)_v,(Z(j)\cap \Delta))=0$ (compare the technique of calculating intersection multiplicities used in the propositions of section 5). Further  $((Z(j) \cap D)_h,(Z(j)\cap \Delta))=2r$ and thus $((Z(j) \cap D),(Z(j)\cap \Delta))=0+2r=2r$ (see also below for the additivity of intersection multiplicities in $Z(j)$ used here). On the other hand, $(D, \Delta, Z(j))$ can  be written as $((D \cap \Delta), (D \cap Z(j)) )$ (intersection multiplicity in $D$, again defined to be the Euler-Poincaré characteristic of the derived tensor product of the corresponding structure sheaves).
Repeating the reasoning before the claim with $Z(j_1)$ replaced by $\Delta \cap D$ and $D$ replaced by $rD_p$ and $\Delta$ replaced by $h$, we see the following identity of intersection multiplicities in $D$ 
$$
 2r= \chi(\mathcal{O}_{D\cap \Delta} \otimes^{\mathbb{L}} \mathcal{O}_{Z(j) \cap D} )= (h, Z(j) \cap D) 
  + (rD_p, Z(j) \cap D)  - (h \cap rD_p, Z(j) \cap D)
 $$Now $ (h \cap rD_p, Z(j) \cap D)=0$, since the structure sheaf of $h \cap rD_p$ is a skyscraper sheaf (see also \cite{KR1}, Lemma 4.1). (Analogously one sees the additivity of intersection multiplicities in $Z(j)$ used above.) 
As above we see that $(rD_p, (Z(j) \cap D)_v)=0$ (applying  the above exact sequence, we use that the support of $D\cap (Z(j) \cap D)_v=(Z(j) \cap D)_v$ lies in the special fibre and is proper over $\F$.) Thus  $(rD_p, (Z(j) \cap D))=2r.$
It follows that $(h, Z(j) \cap D)=0$ and hence $h_x=0.$ Since $x$ was arbitrary, it follows that $h=0.$
This confirms the claim.

Now repeating this reasoning we see that $$(Z(j_1),Z(j_2),Z(j_3))=\sum_{l}(D(j_1/p^l),Z(j_2),Z(j_3))$$ (note that in the remaining steps the reasoning in the last part of the proof  becomes simpler since (in contrast to to above case $l=r$) for $l< r$ there are no  horizontal components  in $D(j/p^l) \cap Z(j)$ if $j$ is as above.)
The remaining multilinearity in the other variables follows from the multilinearity of intersections of divisors (such that the the support of their intersection lies in the special fibre and is proper over $\F$) in $D(j_1/p^{l})$, see the beginning of section 5.
\qed
\newline

The proposition and its proof  show that (if at least one of $j_1,j_2,j_3$ is odd) it is enough to calculate the intersection multiplicity $(D(j_i),D(j_l),D(j_k))$ or (if at least one of $j_i, j_l$ is odd) $(D(j_i),Z(j_k),D(j_l))$, where $\{i,k,l\}=\{1,2,3\}$. To calculate such intersections, we use the following observation (already made and used in the proof of the last proposition). 

By Lemma \ref{Lweg}, we have $\mathcal{O}_{D(j_i)}\otimes^{\mathbb{L}}\mathcal{O}_{D(j_l)}=\mathcal{O}_{D(j_i)}\otimes \mathcal{O}_{D(j_l)}=\mathcal{O}_{D(j_i)\cap D(j_l)}$. (Here, each $D(j_l)$ or $D(j_i)$ or both may be replaced by $Z(j_l)$ resp. $Z(j_i)$.) Thus
$$
\chi(\mathcal{O}_{D(j_i)}\otimes^{\mathbb{L}}\mathcal{O}_{D(j_l)}\otimes^{\mathbb{L}}\mathcal{O}_{D(j_l)})=\chi(\mathcal{O}_{D(j_i) \cap D(j_l)}\otimes^{\mathbb{L}}\mathcal{O}_{D(j_k)})=\chi(\mathcal{O}_{D(j_i) \cap D(j_l)}\otimes^{\mathbb{L}}_{\mathcal{O}_{D(j_i)}}\mathcal{O}_{D(j_i) \cap D(j_k)})).
$$ The same is true if $D(j_l)$ or $D(j_k)$ is replaced by $Z(j_l)$ resp. $Z(j_k)$. If the three special endomorphisms are pairwise perpendicular to each other (which we may assume), then
 by 1.16 of \cite{Ku1} (see also (\ref{iden})), at least one of the three special endomorphisms is odd, so we may also assume that $j_i$ is odd.
This suggests   to investigate intersections of the form $D(y_1)\cap D(y_2)$ where we assume that  $y_1, y_2 \in V^{'}$ such that $y_1^2, y_2^2 \in \Zlp\setminus \{0\}$ and $y_1 \perp y_2$ (and we also may assume that $y_1$ is odd).  Note that by the proof of Lemma \ref{Lweg} (which shows that $D(y_1)$ and $D(y_2)$ have no common component), $D(y_1)\cap D(y_2)$ is a divisor in $D(y_1)$. 
\newline

Let $[\Lambda]\in \mathscr{B}$. Let $y_1, y_2 \in V^{'}$ such that $y_1^2, y_2^2 \in \Zlp\setminus \{0\}$ and $y_1 \perp y_2$ and  $\mathbb{P}_{[\Lambda]}\subset D(y_1)\cap D(y_2)$. The divisor $D(y_i)_p$ in ${\MHB}_p$ contains $\mathbb{P}_{[\Lambda]}$ with multiplicity $p^{r_i}$ for a suitable integer $r_i$ given by Corollary \ref{spezf}. Now $\mathbb{P}_{[\Lambda]}$ is a closed irreducible reduced subscheme of codimension $1$ in (the regular) $D(y_i)$, hence a prime divisor. The subsequent two propositions will tell us the multiplicity of  $\mathbb{P}_{[\Lambda]}$ in the divisor $D(y_1) \cap D(y_2)$ in $D(y_1)$.
\begin{Pro} \label{unglmult}
In the situation just described, suppose further that $r_1 \neq r_2 $. Then the divisor $D(y_1)\cap D(y_2)$ in  $D(y_1)$ contains  $\mathbb{P}_{[\Lambda]}$ with multiplicity $p^{\min\{r_1,r_2\}}$.
\end{Pro}
\emph{Proof.} 
We choose $x \in \mathbb{P}_{[\Lambda]}(\F)$ not superspecial.
We choose $y\in V^{'}_p$ such that $ y \perp y_1,y_2$ and such that $x \in Z(y)(\F)$ and $\nu_p(y^2)=1$ (Lemma \ref{drinsenk}). Let $R= \mathcal{O}_{\MHB, x}$ and let $(f_1), \ (f_2), \ (d)$ be the ideals of $D(y_1), \ D(y_2), \  Z(y)$ in $R$. 
It follows from Lemma \ref{diperg} and from \cite{T}, chapter 2 that the ideal of $\mathbb{P}_{[\Lambda]}$ in $x$ is given by $(p,d)=(f_1,d)=(f_2,d)$.

{\bf Claim 1} \emph{We have the following identity of ideals in $R$:} $(d^{p^{r_1}},f_1)=(d^{p^{r_1}},p)=(f_1,p)$. 

For any $z \in R$ we denote by $\overline{z}$ the image of $z$ in $\overline{R}:=R/(p)$ which is a UFD.
First we observe that $(d^{p^{r_1}},p)=(f_1,p)$, since the equation of $\mathbb{P}_{[\Lambda]}$ in $\overline{R}$ is given by $\overline {d}=0$.
After perhaps multiplying $f_1$ by a unit we can therefore write 
$$f_1=d^{p^{r_1}}+p \varrho, $$ for some $\varrho \in R$. By the above descripion of the ideal of $\mathbb{P}_{[\Lambda]}$ in $R$ we can also write
$$f_1=\varepsilon p + \sigma d$$ for some $\varepsilon \in R^{\times}$ and some $\sigma \in R$.
From these equations we get $\overline{f_1}=\overline{d}^{p^{r_1}}= \overline{\sigma}\overline{d}$, hence $\overline{\sigma}=\overline{d}^{p^{r_1}-1}$, hence $\sigma=d^{p^{r_1}-1}+p\sigma^{'}$ for some $\sigma^{'}\in R$.
Hence we have $f_1=\varepsilon p+ d^{p^{r_1}}+pd\sigma^{'}= p(\varepsilon+d\sigma^{'})+d^{p^{r_1}}$. Since $\varepsilon+d\sigma^{'}$ is a unit it follows that $p\in (f_1, d^{p^{r_1}})$.
Hence $(f_1,d^{p^{r_1}})=(f_1,d^{p^{r_1}},p)=(f_1,p)$. This confirms the claim.

Now we distinguish the cases $r_2 < r_1$ and $r_1 < r_2$.

{\bf First case} $r_2 < r_1$.
Since the ideal of $\mathbb{P}_{[\Lambda]}$ in the local ring of $D(y_1)$ in $x$ equals $(d)$, it is enough is enough to show that  $(f_1,f_2) \subset (d^{p^{r_2}},f_1)$ and that  $(f_1,f_2)\not \subset (d^{p^{r_2}+1},f_1)$. This is the content of claims 2 and 3. 

{\bf Claim 2}  $(f_1,f_2) \subset (d^{p^{r_2}},f_1)$.

We have $(f_2)\subset (d^{p^{r_2}},p)$, hence $(f_1,f_2)\subset (d^{p^{r_2}},f_1,p)$. The latter ideal equals by  claim 1 the ideal $(d^{p^{r_2}},d^{p^{r_1}},f_1)= (d^{p^{r_2}},f_1)$ since $r_2 < r_1$. Hence $ (f_1,f_2)\subset (d^{p^{r_2}},f_1)$ as claimed.

{\bf Claim 3} $(f_1,f_2)\not \subset (d^{p^{r_2}+1},f_1)$.

Suppose $(f_1,f_2) \subset (d^{p^{r_2}+1},f_1)$. Then $(f_1,f_2,p) \subset (d^{p^{r_2}+1},f_1,p)$. Since $r_2 < r_1$ we have $(f_1,f_2,p)=(d^{p^{r_2}},p)$ and $(d^{p^{r_2}+1},f_1,p)=  (d^{p^{r_2}+1},p)$. Thus we get 
$(d^{p^{r_2}},p)\subset (d^{p^{r_2}+1},p)$, a contradiction which confirms the claim.

Combining claims 2 and 3 ends the proof in case $r_2 < r_1$.

{\bf Second case} $r_2 > r_1$.
By the first case, the vertical part (the part with support in the special fibre) of $D(y_1)\cap D(y_2)$ is in $R$ given by the ideal $(f_2, d^{p^{r_1}})$. Now $(f_2, d^{p^{r_1}})=(f_2, d^{p^{r_1}}, f_1)= (d^{p^{r_1}}, f_1)$.
(The latter equality results from the fact that by claim 1 both ideals describe the ideal of the special fibre of $D(y_1)$ in $R$.)
\qed

\begin{Pro}\label{glmult}
Suppose we are in the situation described before Proposition \ref{unglmult}. Suppose further that $r_1=r_2=:r$. Let $a_i=\nu_p(y_i^2)$ and write $y_i^2= \varepsilon_ip^{a_i}$ und suppose that $a_1$ and $a_2$ are not both even and that $a_1 \leq a_2$.
Suppose further that in case $a_1=a_2$ we have $\chi(-\varepsilon_1\varepsilon_2)=-1$ (here, $\chi$ is the quadratic residue character of $\Zp^{\times}$). Then $a_1$ is odd and the divisor $D(y_1)\cap D(y_2)$ in $D(y_1)$ contains $\mathbb{P}_{[\Lambda]}$ with multiplicity $(\frac{a_1+1}{2}-r)\cdot p^r$. If $a_2$ is also odd, then  the divisor $D(y_1)\cap D(y_2)$ in $D(y_2)$ also contains $\mathbb{P}_{[\Lambda]}$ with multiplicity $(\frac{a_1+1}{2}-r)\cdot p^r$. 
\end{Pro}
Before we prove this, we show the following lemma.
\begin{Lem}\label{gf0}
Suppose that $y_1, y_2\in V^{'}$  are both odd and that $y_1 \perp y_2$.  Let $a_i=\nu_p(y_i^2)$ and write $y_i^2=\varepsilon_ip^{a_i}$ with $a_1 \leq a_2$ and suppose $\chi(-\varepsilon_1\varepsilon_2)=-1$. (Here, $\chi$ denotes the quadratic residue character of $\Zp^{\times}$). Then the generic fibre of $Z(y_1)\cap Z(y_2)$ is empty (i.e., $Z(y_1)\cap Z(y_2)$ has support in the special fibre).
\end{Lem}
\emph{Proof of Lemma \ref{gf0}.} Since $y_1,y_2 \in V^{'},$ it is enough to show that there are no special cycles $Z=Z(\varepsilon_1p^{a_1}, \omega)$  and $Z^{'}=Z(\varepsilon_2 p^{a_2}, \omega^{'})$ such that there is some point in the generic fibre  of  $Z\times_{\M}Z^{'}$ with fundamental matrix $\diag(\varepsilon_1p^{a_1},\varepsilon_2p^{a_2}).$ (See also the incidence relation below Proposition \ref{locglob}.) Recall our fixed lattice $\Lambda \subset V$ in section 1. By the description of the generic fibre of  special cycles given in \cite{KR2}, section 2, it is enough to show the following

{\bf Claim 1} {\em There are no elements $\lambda_1, \lambda_2\in \omega_p:= \Lambda \otimes \Zp$ such that $\lambda_1 \perp \lambda_2$ and $q(\lambda_i)=\varepsilon_ip^{a_i}$}.

Here, $q$ denotes the extension to $\omega_p$ of the quadratic form $q$ on $\Lambda$. We will write $\langle ,\rangle $ for the induced bilinear form on $\omega_p$. We also write $\lambda^2$ for $q(\lambda),$ and we write $\langle, \rangle$ for the corresponding bilinear form, i.e. $\langle x, y\rangle =  \frac{1}{2}((x+y)^2-x^2-y^2)$. 

Suppose there are $\lambda_1,\lambda_2$ as in the claim. We may suppose that both of  $\lambda_1,\lambda_2$ are not divisible by $p$. 

If  $x_1,...,x_4$ is a orthogonal basis of $\omega_p$ and if $q_i=x_i^2,$ then, since $\Lambda$ is self-dual, we have $q_i \in \Zp^{\times}$ for all $i$. Further we have $\chi(q_1q_2q_3q_4)=-1$ (comp. \cite{KR2}, 7.14).

We choose $z_1\in \omega_p$ such that $z_1^2\in \Zp^{\times}$ and $\langle z_1,\lambda_1 \rangle \in \Zp^{\times}$ (such $z_1$ exists, since $\lambda_1$ is not divisible by $p$).  Let $\sigma_1=z_1^2$ and let $z_2=\lambda_1-\frac{\langle z_1,\lambda_1 \rangle }{\sigma_1}z_1$. Then $0 \neq z_2\perp z_1$ and $z_2^2=:\sigma_2 \in \Zp^{\times}$. Let $b_1=\frac{\langle z_1,\lambda_1 \rangle }{\sigma_1}$. Then  $b_1^2\sigma_1+\sigma_2 \equiv 0 \mod p$, hence $\chi(-\sigma_1\sigma_2)=1$. We extend $z_1,z_2$ to an orthogonal basis $z_1,...,z_4$ of $\omega_p$. Let $\sigma_i=z_i^2$. Then we conclude that $\chi(-\sigma_3\sigma_4)=-1$. We write $\lambda_2= c_1z_1+c_2z_2+c_3z_3+c_4z_4$.
We already know that $\lambda_1=b_1z_1+z_2$.
Since $\lambda_1 \perp \lambda_2$ we have $b_1c_1\sigma_1+c_2\sigma_2=0$.

{\bf Claim 2} $c_2 \in \Zp^{\times}$.

Suppose $c_2 \equiv 0 \mod p$. Since $b_1$ and $\sigma_1$ are units it follows that $c_1 \equiv 0 \mod p$. Hence $c_3^2\sigma_3+c_4^2\sigma_4 \equiv 0 \mod p$. If $c_3$ or $c_4$ is divisible by $p$, it follows that both are divisible by $p$. But then $\lambda_2$ is also divisible by $p$, a contradiction. Hence $c_3$ and $c_4$ are units. Hence $\chi(\sigma_3)=\chi(c_3^2\sigma_3)=\chi(-c_4^2\sigma_4)= \chi(-\sigma_4)$, and hence $\chi(-\sigma_3\sigma_4)=1$. This contradiction confirms the claim that $c_2$ is a unit.

By claim 2 we may assume that $c_2=1$.
From this we get $c_1= - \frac{\sigma_2}{b_1\sigma_1}$.
Next we observe that $\nu_p((c_1z_1+z_2)^2)=\nu_p(\frac{\sigma_2^2}{b_1^2\sigma_1^2}\sigma_1 +\sigma_2) 
=\nu_p(b_1^2\sigma_1+\sigma_2)=\nu_p(\lambda_1^2)=a_1$. Since $a_1 \leq a_2$ it follows that $c_3^2\sigma_3+c_4^2\sigma_4 \equiv 0 \mod p^{a_1}$. Since $\chi(-\sigma_3\sigma_4)=-1$ and $a_1$ is odd it follows that $c_3,c_4 \equiv 0\mod p^{\frac{a_1+1}{2}}$. Hence $c_3^2\sigma_3+c_4^2\sigma_4 \equiv 0 \mod p^{a_1+1}$.

We introduce the following notation. For $x=\eta p^h \in \Zp$ with $\eta \in \Zp^{\times}$ we write $\chi(x)=\chi(\eta)$.

It follows that $ \chi(\varepsilon_2 )=\chi(\lambda_2^2)=\chi(c_1^2\sigma_1+\sigma_2)$. Hence 
\begin{equation*}
\begin{split}
& \chi(-\varepsilon_1\varepsilon_2 ) 
 =\chi(-(b_1^2\sigma_1+\sigma_2)(c_1^2\sigma_1+\sigma_2))  
=\chi(-(b_1^2\sigma_1+\sigma_2)( \frac{\sigma_2^2}{b_1^2\sigma_1^2}\sigma_1+\sigma_2)) \\
&=\chi(-(b_1^2\sigma_1+\sigma_2)( \sigma_2^2\sigma_1+\sigma_2 b_1^2\sigma_1^2)) 
=\chi(-\sigma_1\sigma_2(b_1^2\sigma_1+\sigma_2)(\sigma_2+b_1^2\sigma_1))=\chi(-\sigma_1\sigma_2)=1.
\end{split}
\end{equation*} 
This contradicts the assumption of the lemma.
\qed
\newline

\emph{Proof of Proposition \ref{glmult}.}
Since $y_1\perp y_2$, the intersection $Core(y_1)\cap Core(y_2)$ is not empty. If $y_1$ was even, then $Core(y_1)$ would consist of  a single point and $y_2$ would be odd, hence $a_2 > a_1$. Now using the Description of the special fibres of $D(y_1)$ and $D(y_2)$ given by Corollary \ref{spezf}, it is easy to see that then the case $r_1=r_2$ cannot occur. Hence  $a_1$ is odd.

 To prove the claim on the multiplicity of $ \mathbb{P}_{[\Lambda]}$ in $D(y_1)\cap D(y_2)$ we use induction on $r$, starting with the case $r=0$.

We write $\langle, \rangle$ for the bilinear form corresponding to $Q$ on $V^{'}$, i.e. $\langle x, y\rangle =  \frac{1}{2}(Q(x+y)-Q(x)-Q(y))$.   
We choose an $\F$-valued point $x\in \mathbb{P}_{[\Lambda]}$ which is superspecial and such that 
if  $x$ is the intersection point $\mathbb{P}_{[\Lambda]}$ and $\mathbb{P}_{[\Lambda^{'}]}$, then $\mathbb{P}_{[\Lambda^{'}]} \not \subset Z(y_i), i=1,2$.
Next we choose $y_0\in V^{'}$ such that $\nu_p(y_0^2)=0$ and $x\in Z(y_0)(\F)$ and $\langle y_0,y_1\rangle \neq 0.$ (Such $y_0$ exists since otherwise $\mathbb{P}_{[\Lambda^{'}]} \subset Z(y_1)$.)
We may suppose that $y_2$ is linearly independent of $y_0$ and $y_1$. (Otherwise one chooses  $y_3 \in V^{'}$ such that $ y_3\perp y_0,y_1,y_2$ and $x\in Z(y_3)(\F)$ and replaces $y_0$ by $y_0+py_3$.)
Let $z_1=y_0$ and $z_2= y_1-\frac{\langle y_0,y_1\rangle }{\langle y_0,y_0\rangle }y_0$, let $b_1=\frac{\langle y_0,y_1\rangle }{\langle y_0,y_0\rangle }$. We extend $z_1,z_2$ to an orthogonal basis $z_1,...,z_4$ of the quadratic $\Zlp$-module $N_x$ of special endomorphisms $y \in V^{'}$ such that $x \in Z(y)(\F)$. Let $\sigma_i=z_i^2$. Then $\sigma_1$ and $\sigma_2$ are units and (using the notation introduced in the proof of Lemma \ref{gf0}) we have $\chi(-\sigma_3\sigma_4)=-1$ (see also section 3).

We write $y_1=b_1z_1+z_2$ and $y_2=c_1z_1+c_2z_2+c_3z_3+c_4z_4$. Then $c_1$ and $c_2$ are units since otherwise both would be divisible by $p$ (since $a_2$ is divisible by $p$) and then we would have $\mathbb{P}_{[\Lambda^{'}]}\subset Z(y_2)$. Hence we may assume that $c_2=1$. Hence $b_1c_1\sigma_1+\sigma_2=0$. Arguing as in Lemma \ref{gf0} we see that $\nu_p(c_1^2\sigma_1+\sigma_2)=a_1$ and that $c_3^2\sigma_3+c_4^2\sigma_4$ is divisible by $p^{a_1}$.
Using $\chi(-\sigma_3\sigma_4)=-1$ we see that $c_3$ and $c_4$ are divisible by $p^{\frac{a_1-1}{2}}$. Then $c_3z_3+c_4z_4=p^{\frac{a_1-1}{2}}y_4$, where $y_4 \in N_x$ is an odd special endomorphism. If $a_1<a_2,$ then it follows immediately that $\nu_p(y_4^2)=1$. Again reasoning as in Lemma \ref{gf0} we see that $\chi(-(b_1^2\sigma_1+\sigma_2)(c_1^2\sigma_1+\sigma_2))=1$. Thus also in case $a_1=a_2$ (in which case by assumption $\chi(-\varepsilon_1\varepsilon_2)=-1$ and therefore $\chi(\varepsilon_2)\neq \chi(c_1^2\sigma_1+\sigma_2)$) it follows that $\nu_p(y_4^2)=1$.
It follows that the fundamental matrix of $y_0,y_1,y_2$ can be diagonalized to the matrix
$$
\begin{pmatrix}
\eta_0&& \\
&\eta_1&  \\
&&\eta_2p^{a_1}\\
\end{pmatrix},
$$ where $\eta_i\in \Zlp^{\times}$. Thus by \cite{KR2}, Proposition 6.2,  it follows that $(Z(y_0),Z(y_1),Z(y_2))=\frac{a_1+1}{2}$.

Next we choose $y\in V^{'}_p$ such that $x \in Z(y)(\F)$ and $y \perp y_0,y_1$ and $\nu_p(y^2)=1$ (Lemma \ref{drinsenk}). Then the fundamental matrix of $y,y_0,y_1$ can be diagonalized to the matrix 
$$
\begin{pmatrix}
\eta_0&& \\
&\eta_1&  \\
&&y^2 \\
\end{pmatrix}. 
$$ Thus the intersection multiplicity  $(Z(y_0),Z(y_1),Z(y))$ equals $1$ (Lemma \ref{diperg}). 
Since $\mathbb{P}_{[\Lambda]}\subset Z(y)$ it follows that $(Z(y_1)\cap Z(y_0),\mathbb{P}_{[\Lambda]})=(D(y_1)\cap Z(y_0),\mathbb{P}_{[\Lambda]})=1$ (intersection multiplicity in $D(y_1)$). 
    Hence it follows (writing  $(Z(y_0),Z(y_1),Z(y_2))=((Z(y_1)\cap Z(y_0), (Z(y_1)\cap Z(y_2))=((D(y_1)\cap Z(y_0), (D(y_1)\cap D(y_2))$ as intersection multiplicity in $D(y_1)$) that the multiplicity of $\mathbb{P}_{[\Lambda]}$ in $D(y_1) \cap D(y_2)$ 
as divisor in $D(y_1)$ is $\leq\frac{a_1+1}{2}$ and that it is  $<\frac{a_1+1}{2}$ if and only if there is a horizontal component of $D(y_1)\cap D(y_2)$ passing through $x$. The same reasoning shows that for $a_2$ odd the multiplicity of $\mathbb{P}_{[\Lambda]}$ in $D(y_1) \cap D(y_2)$ as divisor in $D(y_2)$ is $\leq\frac{a_1+1}{2}$ and that it is  $<\frac{a_1+1}{2}$ if and only if there is a horizontal component of $D(y_1)\cap D(y_2)$ passing through $x$. Thus in  case $a_1=a_2$ this ends the induction start  by Lemma \ref{gf0},  which guarantees that there is no such horizontal component, and we also see that it suffices to show the claim in case $a_1\neq a_2$ for the multiplicity of   $\mathbb{P}_{[\Lambda]}$ in $D(y_1) \cap D(y_2)$ as divisor in $D(y_1).$

Next it is easy to see that there is a special endomorphism $y_3$ of the form $c_1z_1+z_2+\rho$, where $\rho$ is in the $\Zlp$-span of $z_3$ and $z_4$ in $V^{'},$ such that $y_3^2=\varepsilon_3p^{a_1}$ (where $\varepsilon_3 \in \Zlp^{\times}$) and $\chi(-\varepsilon_1\varepsilon_3)=-1$.  We fix such a $y_3$. 

Thus $y_2-y_3=p^{\frac{a_1-1}{2}}y_4$, where $y_4$ is some element in the $\Zlp$-span of $z_3$ and $z_4$ in $V^{'}$.

As just shown, we already know that the multiplicity of $\mathbb{P}_{[\Lambda]}$ in $D(y_1) \cap D(y_3)$ (or $D(y_1) \cap Z(y_3)$) as divisor in $D(y_1)$ is $\frac{a_1+1}{2}$. Further, $D(y_1) \cap Z(y_4)$ (as divisor in $D(y_1)$) contains $\mathbb{P}_{[\Lambda]}$ with multiplicity at least $1$. Now Proposition \ref{pmult} shows that $D(y_1) \cap Z(p^{\frac{a_1-1}{2}}y_4)$ (as divisor in $D(y_1)$) contains  $\mathbb{P}_{[\Lambda]}$ with multiplicity at least $1+\frac{a_1-1}{2}=\frac{a_1+1}{2}$.
Hence we see that $D(y_1) \cap Z(y_2)=D(y_1) \cap Z(y_3+p^{\frac{a_1-1}{2}}y_4)$ also contains  $\mathbb{P}_{[\Lambda]}$ with multiplicity at least $\frac{a_1+1}{2}$. Since locally around $x$ we $D(y_2)=Z(y_2)$ and 
since we already know that this multiplicity is at most $\frac{a_1+1}{2}$, this ends the start of the induction.
\newline

Now we come to the induction step from $r-1$ to $r$.

If  $y_1$ and $y_2$ would be replaced by  $y_1/p^r$ and $y_2/p^r,$ then we would be in the situation of the induction start. Now
 we 
choose an $\F$-valued point $x$ of $\mathbb{P}_{[\Lambda]}$ for $y_1/p^r$ and $y_2/p^r$ as we did in the induction start, and also, as before, $x $ is the intersection point of   $\mathbb{P}_{[\Lambda]}$ and  $\mathbb{P}_{[\Lambda^{'}]}$. We also choose $y_0$ for  $y_1/p^r$ and $y_2/p^r$ as before. Then the  fundamental matrix of $y_0,y_1,y_2$ can be diagonalized to the matrix
$$
\begin{pmatrix}
\eta_0&& \\
&\eta_1 p^{2r}&  \\
&&\eta_2p^{a_1}\\
\end{pmatrix},
$$ where $\eta_i\in \Zlp^{\times}$.
Further it follows that the  fundamental matrix of $y_0,y_1/p,y_2$ can be diagonalized to the matrix $\diag(\eta_0,\eta_1 p^{2r-2}, \eta_2p^{a_1})$ and  that the  fundamental matrix of $y_0,y_1/p,y_2/p$ can be diagonalized to the matrix $\diag(\eta_0,\eta_1 p^{2r-2}, \eta_2p^{a_1-2}).$ Analogously one sees that the  fundamental matrix of $y_0,y_1,y_2/p$ can be diagonalized to a matrix of the form $\diag(\eta_0,\tilde{\eta_1} p^{2r-2}, \tilde{\eta_2}p^{a_1}),$ where $\tilde{\eta_i}\in \Zlp^{\times}.$
Using $D(y_1)=Z(y_1)-Z(y_1/p)$ and $D(y_2)=Z(y_2)-Z(y_2/p)$  and \cite{KR2}, Proposition 6.2, one checks that \begin{equation} \label{bez}
(Z(y_0),D(y_1),D(y_2))=(\frac{a_1+1}{2}-r+1)p^{r-1}+(\frac{a_1+1}{2}-r)p^r.
\end{equation}
Using the induction hypothesis (applied to $y_1,y_2$ and the projective line  $\mathbb{P}_{[\Lambda^{'}]}$, note that $\mathbb{P}_{[\Lambda^{'}]}$ has multiplicity $p^{r-1}$ in $D(y_1)_p$ and in $D(y_2)_p$) we see that the multiplicity of  $\mathbb{P}_{[\Lambda^{'}]}$ in $D(y_1) \cap D(y_2)$ as divisor in $D(y_1)$ (and also as divisor in $D(y_2)$ if $a_2$ is odd) is $(\frac{a_1+1}{2}-(r-1))p^{r-1}$. 

As before, we choose $y\in V^{'}_p$ such that $y\perp y_0, y_1$ and  such that $\nu_p(y^2)=1$ and $x \in Z(y)(\F)$.
  Then the fundamental matrix  of $y,y_0,y_1$ can be diagonalized to the matrix 
$$
\begin{pmatrix}
\eta_0&& \\
&y^2&  \\
&&p^{2r}\eta_1 \\
\end{pmatrix}. 
$$ Thus the intersection multiplicity  $(Z(y_0),D(y_1),Z(y))$ equals $2$ (use Lemma \ref{diperg}).

Since
 $\mathbb{P}_{[\Lambda]}, \mathbb{P}_{[\Lambda]^{'}}\subset Z(y)\cap D(y_1)$ it follows that  $(Z(y_0)\cap D(y_1),\mathbb{P}_{[\Lambda]})=(Z(y_0)\cap D(y_1),\mathbb{P}_{[\Lambda^{'}]})=1$ (intersection multiplicities in $D(y_1)$ ). Using (\ref{bez}) and  the above expression for the multiplicity of  $\mathbb{P}_{[\Lambda^{'}]}$ in $D(y_1) \cap D(y_2),$ 
 we see that the multiplicity of  $\mathbb{P}_{[\Lambda]}$ in $D(y_1) \cap D(y_2)$ as divisor in $D(y_1)$ is $ \leq (\frac{a_1+1}{2}-r)p^r$, with equality if and only if there is no horizontal component of $D(y_1)\cap D(y_2)$ passing through $x$. Thus, again by Lemma \ref{gf0}, we are done in case $a_1=a_2$. The same reasonig shows in case that $a_2$ is odd the corresponding statement in $D(y_2).$ Thus for the rest we may restrict ourselves to the case that we consider $D(y_1) \cap D(y_2)$ as divisor in $D(y_1).$ 
 
 The remaining reasoning (for the  case $a_1 \neq a_2$ ) is analogous to the reasoning in the induction start: One chooses $y_3 \in V^{'}$ such that $y_2-y_3$ is a linear combination of $z_3$ and $z_4$  (in particular $y_3 \perp y_1$) and  such that $x \in Z(y_3)(\F)$ and  $\nu_p(y_3^2)=a_1$ and writing $y_3^2=\varepsilon_3 p^{a_1}$ we have 
 $\chi(-\varepsilon_1\varepsilon_3)=-1.$ Then $y_2-y_3=p^{\frac{a_1-1}{2}}y_4$ for some $y_4$ such that $x \in Z(y_4)(\F).$  Then as just shown, $D(y_1)\cap D(y_3)$ as divisor in $D(y_1)$ contains  $\mathbb{P}_{[\Lambda]}$ with multiplicity  $ (\frac{a_1+1}{2}-r)p^r.$ Thus (using Proposition \ref{unglmult}) we see that  $D(y_1)\cap Z(y_3)$ as divisor in $D(y_1)$ contains  $\mathbb{P}_{[\Lambda]}$ with multiplicity  $ 1+ ... + p^{r-1}+ (\frac{a_1+1}{2}-r)p^r.$ Now Proposition \ref{unglmult} and claim 1 in its proof show that $Z(p^r y_4)\cap D(y_1)$ as divisor in $D(y_1)$ contains $\mathbb{P}_{[\Lambda]}$ with multiplicity at least $1+ ... + p^r.$ Now using Proposition \ref{pmult} we see that  $Z(p^{\frac{a_1-1}{2}}y_4) \cap D(y_1)$  as divisor in $D(y_1)$ contains $\mathbb{P}_{[\Lambda]}$ with multiplicity at least $1+ ... + p^r+ (\frac{a_1-1}{2}-r)p^r=1+ ... + p^{r-1}+ (\frac{a_1+1}{2}-r)p^r.$  Thus $D(y_1) \cap Z(y_2)=D(y_1) \cap Z(y_3+p^{\frac{a_1-1}{2}}y_4)$ also contains  $\mathbb{P}_{[\Lambda]}$ with multiplicity at 
 least $1+ ... + p^{r-1}+ (\frac{a_1+1}{2}-r)p^r $. Thus we see (again using Proposition \ref{unglmult}) that the multiplicity of  $\mathbb{P}_{[\Lambda]}$ in $D(y_1)\cap D(y_2)$ (as divisor in $D(y_1)$) is at least $(\frac{a_1+1}{2}-r)p^r. $
Since we already know that this multiplicity is at most $ (\frac{a_1+1}{2}-r)p^r,$ this is the precise multiplicity. \qed
\newline

Propositions \ref{unglmult} and \ref{glmult} give us the necessary informations about the vertical part (i.e. the part having support in the special fibre) of the divisors $D(j_i)\cap D(j_l)$ in say $D(j_i)$. Next we investigate the horizontal components (the components which do not have support in the special fibre). We consider two special endomorphisms  $y_1, y_2 \in V^{'}$ such that $y_1 \perp y_2$, and we suppose that $y_1$ is odd. By the results of of \cite{KR2} and \cite{T} we may assume that $\nu_p(y_i)\geq 2$. (In case that the three special endomorphisms $j_i$ are pairwise perpendicular to each other and at least one of the exponents $\nu_p(j_i^2)$ is $\leq 1$, the intersection multiplicity $(Z(j_1),Z(j_2),Z(j_3))$ is known by \cite{KR2}, resp. \cite{T}.) We start with the following lemma, which is given a geometric interpretation in the subsequent proposition.
\begin{Lem}\label{horizl}
Let $y_1, y_2 \in V^{'}$ such that $y_1\perp y_2$ and suppose that $y_1$ is odd and that $y_2$ is even. Write $a_i=\nu_p(y_i^2)$. Suppose that $a_1>a_2\geq 2$. Let $x$ be the unique  $\F$-valued point in $Core(y_2)$. Let $R=\mathcal{O}_{\MHB, x}$ and let $f_1$ resp. $f_2$ be generators of the ideals of $D(y_1)$ resp. $D(y_2)$ in $R$. Let $d\in R$ such that the image of $d$ in $R/(f_1)$ is a generator of the ideal of the underlying reduced subscheme of $D(y_1)$. By Proposition \ref{unglmult}, we can write $f_2=d^{\alpha}h+\rho f_1$, for some $h, \rho \in R$ such that the image of $h$ in $R/(f_1)$ is coprime to $p$, and where $\alpha=p^{a_2/2-1}$.
Further let $t \in R$ such that the image of $t$ in $R/(p)$ is a generator of the ideal of the divisor $s$ in Corollary \ref{spezf} (applied in case $j=y_2$). Let $\beta=p^{a_2/2-1}(p-1)$.
Then the ideals $I_1:=(f_1,h,d)$ and $I_2:=(t^{\beta},d,p)$ are equal and the length (over $W$) of $R/I_1$ is $2\beta$.
\end{Lem}
\emph{Proof.}
We choose a special endomorphism $y $ such that $y\perp y_1,y_2$ such that $\nu_p(y)=1$ and $x\in Z(y)(\F)$ (Lemma \ref{drinsenk}). By Lemma \ref{diperg} (and using Lemma \ref{hillem}), we may assume that the ideal of $Z(y)$ in $R$ is $(d)$. Further by Lemma \ref{diperg} we know that  $p \in (f_1,d)$ and $f_1 \in (d,p)$.
Thus $I_1=(f_1,h,d,p)$ and $I_2=(f_1,t^{\beta},d,p)$.
For any $r \in R$ we denote by $\overline{r}$ its image in $R/(p)$ (which is a UFD). Further let $\overline{I}_i$ be the image of ${I}_i$ in $R/(p)$. Then, since $p \in I_1 \cap I_2$, it is enough to show that $\overline{I}_1=\overline{I}_2$.
We write 
$$
f_2=d^{\alpha} t^{\beta}+\sigma p
$$
for some $ \sigma \in R$. Then 
$$
\overline{f_2}=\overline{d}^{\alpha}\overline{h}+\overline{\rho} \overline{ f_1}=\overline{d}^{\alpha}\overline{t}^{\beta}
$$
and 
$$
\overline{f_1}=\overline{d}^{\gamma},
$$
where $\gamma=p^{(a_1-1)/2}$. Combining both equations we get
$$
\overline{t}^{\beta}=\overline{h}+\overline{\rho}\overline{d}^{\gamma-\alpha}. 
$$
Note that $\gamma-\alpha>0$. Thus we see that $\overline{t}^{\beta} \in \overline{I}_1$ and $\overline{h}\in \overline{I}_2$. 
This confirms that $\overline{I}_1=\overline{I}_2$ and hence $I_1=I_2$.

Further $\lg_W(R/I_1)=\lg_W(R/I_2)=\beta \cdot \lg_W(R/(p,d,t))$.  By Lemma \ref{diperg},1 (together with Lemma \ref{hillem} below), the equation for $Z(y_1/p^{(a_1-1)/2})$ in $R/(d)$ is given by $p=0$, hence we see that $\lg(R/(p,d,t))= (Z(y_1/p^{(a_1-1)/2}), Z(y_2/p^{a_2/2}), Z(y))$. The latter equals $2$, as follows again from Lemma \ref{diperg}. Thus $\lg_W(R/I_1) =2\beta $ as claimed. \qed

\begin{Pro}\label{horiz}
In the situation of Lemma \ref{horizl} there is a horizontal part $h$ of the divisor $D(y_1)\cap D(y_2)$ in $D(y_1)$ which passes through the unique point in $Core(y_2)$. Its intersection multiplicity (in $D(y_1)$) with each of the two projective lines which contain the point in $Core(y_2)$ is equal to $\beta=p^{a_2/2-1}(p-1)$. 
\end{Pro}
\emph{Proof.} This follows immediately from Lemma \ref{horizl}.
\qed

\begin{Pro} \label{horiz0}
Let $y_1, y_2 \in V^{'}$ such that $y_1\perp y_2$ and suppose that $y_1$ is odd.  Let $a_i=\nu_p(y_i^2)$ and write $y_i^2= \varepsilon_i p^{a_i}$. Suppose that $a_1,a_2\geq 2$ and that in case $a_1=a_2$ we have $\chi(-\varepsilon_1\varepsilon_2)=-1$. Then in case that $a_1 > a_2$ and that $a_2$ is even,  the horizontal component of $D(y_1)\cap D(y_2)$ passing through  the point in $Core(y_2)$ (see Proposition \ref{horiz}) is the only horizontal component of $D(y_1)\cap D(y_2)$. In all other cases  there does not exist any horizontal component of $D(y_1)\cap D(y_2)$.
\end{Pro}
The proof will be given in the next section.


\section{Calculation of intersection multiplicities}

In view of Proposition \ref{multlin} (and its proof) it is enough to calculate  (if at least one of $j_1,j_2,j_3$ is odd) $(D(j_i),D(j_l),D(j_k))$ or (if at least one of $j_i,j_k$ is odd) $(D(j_i),D(j_k),Z(j_l))$, where $\{i,k,l\}=\{1,2,3\}$. Using Lemma \ref{Lweg} we see (as before) that 
$$(D(j_i),D(j_l),D(j_k))=\chi(\mathcal{O}_{D(j_i) \cap D(j_l)}\otimes^{\mathbb{L}}_{\mathcal{O}_{D(j_i)}}\mathcal{O}_{D(j_i) \cap D(j_k)}),$$ which we write as $((D(j_i)\cap D(j_l)),(D(j_i)\cap D(j_k)))$ (intersection multiplicity in $D(j_i)$). The corresponding analogous formula holds if one replaces $D(j_l)$ by $Z(j_l)$.
We may and will assume that $j_i$ is odd (see below). For any odd $j \in V^{'}$ with $\nu_p(j^2) \geq 1$ by Corollary \ref{spezf} the special fibre of the regular formal scheme $D(j)$ is a scheme and its underlying reduced subscheme is a union of copies of $\mathbb{P}^1_{\F}$. The same reasoning as in \cite{D}  now shows that the intersection number of two divisors $E=E_1+E_2$ and $F$ in $D(j)$ (defined as usual as the Euler-Poincaré characteristic of their structure sheaves) is bilinear i.e. satisfies $(E,F)=(E_1,F)+(E_2,F)$ provided that the support of  $E \cap F$ is contained in $D(j)_p$ and is proper over $\F$.

Next we compute the intersection multiplicity of two of the projective lines in $D(j)_p$.

\begin{Lem}\label{selsch}
Let $j \in V^{'}$ be  odd with $\nu_p(j^2)\geq 1$ and 
let $\mathbb{P}_{[\Lambda]}, \mathbb{P}_{[\Lambda^{'}]} \subset D(j)_p$. Then their intersection multiplicity in $D(j)$ has the value:
$$
(\mathbb{P}_{[\Lambda]}, \mathbb{P}_{[\Lambda^{'}]})=
\begin{cases}
0, &   \ \text{if } \mathbb{P}_{[\Lambda]} \text{ and } \mathbb{P}_{[\Lambda^{'}]} \text{ do not intersect} \\
1, &   \ \text{if } \mathbb{P}_{[\Lambda]} \text{ and } \mathbb{P}_{[\Lambda^{'}]} \text{ intersect in precisely one point} \\
-(p+1), &  \text{if } [\Lambda]= [\Lambda^{'}] \text{ and } \nu_p(j^2)=1 \\
-2p, &  \text{if } [\Lambda]= [\Lambda^{'}] \text{ and } \nu_p(j^2)> 1\text{ and } \mathbb{P}_{[\Lambda]} \subset Z(j/p)\\ 
-p, &  \text{if } [\Lambda]= [\Lambda^{'}] \text{ and } \nu_p(j^2)> 1\text{ and } \mathbb{P}_{[\Lambda]} \not\subset Z(j/p).
\end{cases}
$$
\end{Lem}
\emph{Proof.}
The first case is obvious. In the second case we have to compute the length of the local ring of the intersection point of $\mathbb{P}_{[\Lambda]} \text{ and } \mathbb{P}_{[\Lambda^{'}]}$. We find a special endomorphism $y \in V^{'}_p$ with $\nu_p(y^2)=1$ and such that $Z(y)$ contains $\mathbb{P}_{[\Lambda]} \text{ and } \mathbb{P}_{[\Lambda^{'}]}$. But then the claim follows from \cite{KR1}, Lemma 4.7. The same lemma also shows the third case. 
To compute the intersection multiplicity in the fourth  and in the fifth case, we use the fact that $(\mathbb{P}_{[\Lambda]}, D(j)_p)=0$. This follows as in the proof of Proposition \ref{multlin} from the exact sequence
$$
\begin{CD}
0 @ >>> \mathcal{O}_{D(j)} @> p \cdot>> \mathcal{O}_{D(j)} @>>>\mathcal{O}_{D(j)}/(p)@>>> 0.
\end{CD}
$$
Let $p^a$ be the multiplicity of $\mathbb{P}_{[\Lambda]}$ in $D(j)_p$ (see Corollary \ref{spezf}).
We write $j^2 = p^ b j_0$, where $\nu_p(j_0^2)=1$. Suppose first that $\mathbb{P}_{[\Lambda]}$ does not belong to $Z(j_0)$.
Then using the second case and Corollary \ref{spezf} we therefore get in the fourth case
$$
0=(\mathbb{P}_{[\Lambda]}, D(j)_p)=p^a\cdot (\mathbb{P}_{[\Lambda]}, \mathbb{P}_{[\Lambda]})+ p^2\cdot p^{a-1}+p^{a+1}.
$$
Thus $(\mathbb{P}_{[\Lambda]}, \mathbb{P}_{[\Lambda]})=-2p$.

If  $\mathbb{P}_{[\Lambda]} \subset Z(j_0),$ then we calculate (still in the fourth case)
$$0=p^a\cdot (\mathbb{P}_{[\Lambda]}, \mathbb{P}_{[\Lambda]})+ (p^2-p)\cdot p^{a-1}+(p+1)\cdot p^{a}. $$ (Note that there are $p+1$ projective lines intersecting $\mathbb{P}_{[\Lambda]}$ which also belong to $Z(j_0)$ and there are $p^2-p$ projective lines intersecting $\mathbb{P}_{[\Lambda]}$ which do not  belong to $Z(j_0)$.)
Thus again $(\mathbb{P}_{[\Lambda]}, \mathbb{P}_{[\Lambda]})=-2p$.

In the fifth case we calculate (using $a=0$):
$$
0= (\mathbb{P}_{[\Lambda]}, \mathbb{P}_{[\Lambda]})+ p.
$$
Thus $(\mathbb{P}_{[\Lambda]}, \mathbb{P}_{[\Lambda]})=-p$.
\qed
\newline

In the sequel we will use the following terminology. If $j$ is an odd special endomorphism and $j=p^b j_0$, where $\nu_p(j_0^2)=1$, then {\em the center of $Z(j)$}, denoted by $Cent(j)$, is the set of all $\mathbb{P}_{[\Lambda]}$ which are contained in $Z(j_0)$. 

If $j$ is an arbitrary special endomorphism such that $j^2 \neq 0$, then {\em the boundary of $Z(j)$} is the set of all  $\mathbb{P}_{[\Lambda]}$ which belong to $Z(j)$ but not to $Z(j/p)$ and which do not contain any points of $Core(j)$. It will be denoted by $B(j)$. 

If $E$ is a divisor in $D(j),$ then we denote by $E_v$ the vertical part of $E$, i.e. the part having support in the special fibre.

For any $\mathbb{P}_{[\Lambda]}$ and any $r \in \{ 1,2,3\}$ we denote by $\mult_r(\mathbb{P}_{[\Lambda]})$ the multiplicity of $\mathbb{P}_{[\Lambda]}$ in $D(j_r)_p$ given by Corollary  \ref{spezf}.

Further we denote again by $\chi$ the quadratic residue character of $\Zp^{\times}$ resp. $\Zlp^{\times}.$ 

\begin{Lem} \label{hillem} 

Let $y_1, y_2 \in V^{'}_p$ such that  $y_1\perp y_2$ and $y_1^2,y_2^2 \neq 0$.  Suppose that $a_i=\nu_p(y_i^2) \geq 0$ and write $y_i^2= \varepsilon_i p^{a_i}$. 
\begin{enumerate}[(i)]
\item If $y_1$ and $y_2$ are both odd  and $\chi(-\varepsilon_1\varepsilon_2)=1,$ then $Cent(y_1)\cap Cent(y_2)$ consists of precisely one projective line. 

\item If $y_1$ and $y_2$ are both odd  and $\chi(-\varepsilon_1\varepsilon_2)=-1,$ then $Cent(y_1)\cap Cent(y_2)$ consists of an appartment of projective lines, i.e. each $\mathbb{P}_{[\Lambda]} \in Cent(y_1) \cap Cent(y_2)$ is intersected by precisely two others.

\item If $y_1$ is odd and $y_2$ is even, then  the point in  $Core(y_2)$  is the intersection point of two projective lines belonging to $Cent(y_1)$.

\item If $y_3 \in V^{'}_p$ is a third special endomorphism with $y_3^2=\varepsilon_3 p^{a_3}$ (where $a_3:=\nu_p(y_3^2) \geq 0$) and $y_3 \perp y_1,y_2,$ then $Z(y_1)\cap Z(y_2) \cap Z(y_3)$ is not empty, and if $y_1$ and $y_2$ are odd and $y_3$ is even, then $\chi(-\varepsilon_1\varepsilon_2)=-1$. Further, if $y_1,y_2, y_3$ are all odd, then $Cent(y_1)\cap Cent(y_2)\cap Cent(y_3)$ consists of precisely one projective line. 
\end{enumerate}
\end{Lem}
\emph{Proof.} The first three points follow by combining \cite{T}, Proposition 4.2, together with loc.cit Proposition 2.4 and \cite{KR1}, Corollary 2.6, or can be proved directly like \cite{KR1}, Corollary 2.6.
The first assertion of  point (iv) is clear, the second follows from the first three points. The last assertion follows from \cite{KR2}, Proposition 8.13 or alternatively from [KR1], Proposition 2.12.
\qed

\begin{Lem} \label{mitte0}
Suppose that $y_1, y_2 \in V^{'}$ such that $y_1 \perp y_2$ and that $a_1:=\nu_p(y_1^2)\geq 1$ is odd and that $a_2:=\nu_p(y_2^2)\geq 2$. Write $y_i^2=\varepsilon_ip^{a_i}$. Suppose that in case that 
 $a_1=a_2$ we have $\chi(-\varepsilon_1\varepsilon_2)=-1$.
Let $\mathbb{P}_{[\Lambda]}\subset D(y_1)\cap D(y_2)$  and suppose that $\mathbb{P}_{[\Lambda]} \notin  B(y_2)$. In case that $a_2$ is even and  that $a_2 < a_1$ or $a_1=1$ suppose further that $\mathbb{P}_{[\Lambda]}$ does not contain  the  point in $Core(y_2).$  Then $(\mathbb{P}_{[\Lambda]}, (D(y_2)\cap D(y_1))_v)=0$ (intersection multiplicity in $D(y_1)$). 
\end{Lem}
\emph{Proof.}
Suppose first that $a_1 \geq 3.$
 We distinguish (a priori) 4 times 3 cases. The four cases are given by the distiguishing whether $\mathbb{P}_{[\Lambda]} \in Cent(y_1)$ or not and whether $\mathbb{P}_{[\Lambda]} \in Cent(y_2)$ or not. Let $p^{m_i}$ be the multiplicity of $\mathbb{P}_{[\Lambda]}$ in $D(y_i)_p$ given by Corollary \ref{spezf}. For each of the four cases one  distuishes whether $m_1<m_2$ or $m_1=m_2$ or $m_1 > m_2$.   

For example, if  $\mathbb{P}_{[\Lambda]} \in Cent(y_1)$ and $m_1 < m_2,$ then (using Propositions \ref{unglmult}, \ref{glmult} and Lemma \ref{selsch}) $$(\mathbb{P}_{[\Lambda]},( D(y_2)\cap D(y_1))_v)= (-2p)p^{m_1}+(p+1)p^{m_1}+(p^2-p)p^{m_1-1}=0. $$

Let us consider the case $\mathbb{P}_{[\Lambda]} \notin Cent(y_1) \cup Cent(y_2)$ and $m_1=m_2$.  Then there are $p^2$ projective lines intersecting  $\mathbb{P}_{[\Lambda]}$ such that each is contained in $D(y_1)_p$ and $D(y_2)_p$ with multiplicity $p^{m_1-1}$. There is one projective line intersecting  $\mathbb{P}_{[\Lambda]}$ which is contained in $D(y_1)_p$ and $D(y_2)_p$ with multiplicity $p^{m_1+1}$. Let $a= \min\{a_1,a_2\}.$ Now using Lemma \ref{selsch} and Proposition \ref{glmult} we compute
\begin{equation*}
\begin{split}
& (\mathbb{P}_{[\Lambda]}, (D(y_2)\cap D(y_1))_v)= \\
& (\frac{a+1}{2}-m_1)(-2p)p^{m_1}+(\frac{a+1}{2}-m_1-1)p^{m_1+1}+p^2(\frac{a+1}{2}-m_1+1)p^{m_1-1}=0. 
\end{split}
\end{equation*} 

The other cases are proved analogously, as well as the case $a_1=1$ which is even simpler (use lemma \ref{diperg}).
\qed
\newline

In the sequel we will say that two projective lines have distance $d$, if the corresponding vertices in the building have  distance $d$. 
\newline

\emph{Proof of Proposition \ref{horiz0}.} 
Denote by $(D(y_1) \cap D(y_2))_{vh}$ the divisor in $D(y_1)$ which is given by the part of  $(D(y_1) \cap D(y_2))$ having support in the special fibre and,  in case  that $a_1 > a_2$ and that $a_2$ is even, includes additionally the horizontal part of $D(y_1)\cap D(y_2)$ passing through the point in $Core(y_2)$ (see Proposition \ref{horiz}).

Suppose that $x\in (D(y_1)\cap D(y_2))(\F)$ such that in case that $a_1 > a_2$ and that $a_2$ is even $\{x\}$ is not equal to $Core(y_2)$. Then we must show that there is no horizontal component of $D(y_1)\cap D(y_2)$ passing through $x$. We choose $y \in V_p^{'}$ such that $\nu_p(y^2)=1$ and $y \perp y_1, y_2$ and $x \in Z(y)(\F)$. 
By the results of \cite{T} we know the value of $(D(y_1),D(y_2),Z(y))$. We also know by loc. cit. the structure of $D(y_1)\cap Z(y)$ (Lemma \ref{diperg}). Using this  we  calculate (below)  
$((D(y_1) \cap D(y_2))_{vh}, (D(y_1)\cap Z(y)))$ (intersection multiplicity in $D(y_1)$) and observe that it equals $(D(y_1),D(y_2),Z(y))$. If there was a horizontal component of $D(y_1)\cap D(y_2)$ passing though $x$, then it would cause a strictly positive contribution to the intersection multiplicity which could not be compensated by any other contributions of horizontal components. Hence $(D(y_1) \cap D(y_2))_{vh}=(D(y_1) \cap D(y_2))$. 

Now we come to the calculation of $((D(y_1) \cap D(y_2))_{vh}, (D(y_1)\cap Z(y)))$. 

We distinguish cases according as $a_2$ is even or odd and if $a_1<a_2$ or  $a_1=a_2$ or $a_1>a_2$. 
First we handle the case that $a_2$ even and $a_2 < a_1$. 

For each $\mathbb{P}_{[\Lambda]}\subset D(y_1)\cap Z(y)$ we compute the contribution to  $((D(y_1) \cap D(y_2))_{v}, (D(y_1)\cap Z(y))_v)$ coming from the part of $(D(y_1) \cap Z(y))_v$ which has support in 
$ \mathbb{P}_{[\Lambda]}$. Then we add all these contributions. (More precisely, this means the following. If  $\mathbb{P}_{[\Lambda]}$ is contained in the divisor $D(y_1)\cap Z(y)$ in $D(y_1)$ (then automatically with  multiplicity is $m=1$, see Lemma \ref{diperg}), then we compute $((D(y_1) \cap D(y_2))_{v},\mathbb{P}_{[\Lambda]} )$ as intersection multiplicity in $D(y_1).$)  By Lemma \ref{mitte0} the contribution of $\mathbb{P}_{[\Lambda]}$ is $0$ if $\mathbb{P}_{[\Lambda]}\subset D(y_2)$ and  $\mathbb{P}_{[\Lambda]} \notin  B(y_2)$ and $ Core(y_2) \not\subset \mathbb{P}_{[\Lambda]}$.
If $\mathbb{P}_{[\Lambda]}\subset D(y_1)\cap Z(y)$ and $\mathbb{P}_{[\Lambda]}\in B(y_2),$ then the contribution to  $((D(y_1) \cap D(y_2))_{v}, (D(y_1)\cap Z(y)))$ coming from the part of $(D(y_1) \cap D(y))$ which has support in 
$ \mathbb{P}_{[\Lambda]}$ is $-2p+p=-p$. (where $-2p$ comes from the self-intersection and $p$ comes from the intersection with the unique projective line $ \mathbb{P}_{[\Lambda^{'}]}$ in $D(y_1)\cap D(y_2)$ which intersects $ \mathbb{P}_{[\Lambda]}$. By Proposition \ref{unglmult} it has multiplicity $p$ in $D(y_2)\cap D(y_1)$.) Now any such $ \mathbb{P}_{[\Lambda]}$ is intersected in $D(y_1)$ by $p$ additional projective lines ($\neq  \mathbb{P}_{[\Lambda^{'}]}$) which belong to $Z(y)$ but not to $D(y_2)$. Each of these $p$ projective lines contributes with the value $1$ to the intersection multiplicity. Therefore together with $ \mathbb{P}_{[\Lambda]}$ they contribute with the value $0$ to the intersection multiplicity. Next we come to the contribution of the two projective lines passing through the point in $Core(y_2)$. If $a_2 \geq 4$  each contributes with the value 
$$(-2p)\cdot  p^{a_2/2-1}+p^2\cdot p^{a_2/2-2}+p^{a_2/2-1}=-p^{a_2/2}+p^{a_2/2-1}.$$ In case $a_2=2$ the contribution is $-2p+1.$ In this case we have to add the contributions of the $p$ projective lines which do not belong to $D(y_2)$ but to $Z(y)$ and intersect $ \mathbb{P}_{[\Lambda]}$.  Thus we have in both cases $$((D(y_1) \cap D(y_2))_{v}, (D(y_1)\cap Z(y))_v)=-2p^{a_2/2}+2p^{a_2/2-1}.$$ 
By 
Lemma \ref{diperg}
we have $ (D(y_1)\cap Z(y))_v= (D(y_1)\cap Z(y))$. Using Proposition \ref{horiz} we finally get
$$((D(y_1) \cap D(y_2))_{vh}, (D(y_1)\cap Z(y)))=0.$$ 
By  \cite{T}, Theorem 5.1 (and its proof) we see that this is already the value of $((D(y_1) \cap D(y_2)), (D(y_1)\cap Z(y)))$. Thus $((D(y_1) \cap D(y_2))_{vh}=((D(y_1) \cap D(y_2))$.

Next we treat the case that $a_2$ is odd. Then by Lemma \ref{gf0} and our assumption we only need to treat the case $a_2>a_1$ and $\chi(-\varepsilon_1\varepsilon_2)=1$.
Then again  $ (D(y_1)\cap Z(y))_v= (D(y_1)\cap Z(y))$ (Lemma \ref{diperg}) and we only have to calculate  $((D(y_1) \cap D(y_2))_{v}, (D(y_1)\cap Z(y)_v))$. 
We write $y^2= \varepsilon p$.

{\bf Claim}    \emph{$((D(y_1) \cap D(y_2))_{v}, (D(y_1)\cap Z(y)_v))$ equals $0$ if $\chi(-\varepsilon \varepsilon_1)=1$, and 
  it equals 
 $-2(p-1)p^{\frac{a_1-1}{2}}$ if $\chi(-\varepsilon \varepsilon_1)=-1$}.

The first claim follows from Lemma \ref{mitte0}.

In the second case we compute for each projective line in   $(D(y_1)\cap Z(y))_v$  the contribution to the intersection multiplicity coming  from the part of  $(D(y_1)\cap Z(y))_v$ having support in this projective line. By Lemma \ref{mitte0} we only have to consider such projective lines which are in $B(y_2)$ or which intersect a projective line in $B(y_2)$ and are not in $D(y_2)$. If such a projective line lies in $B(y_1)\cap B(y_2)$ (and lies in $Z(y)$), then its contribution is $(-p)(\frac{a_1+1}{2})+p(\frac{a_1+1}{2}-1)=-p$. There are precisely $2(p-1)p^{\frac{a_1-1}{2}-1}$ such projective lines. (There are precisely two projective lines in the apartment of projective lines belonging to $Cent(y_1)\cap Cent(y_2)$ for which the distance to $B(y_1)$ (defined as the corresponding distance in the building) is $\frac{a_1-1}{2}$.  For each of these two, there are $(p-1)p^{\frac{a_1-1}{2}-1}$ projective lines with the described properties which have distance $\frac{a_1-1}{2}$ to the projective line in the apartment.) 

For all other projective lines in $B(y_2)$, their respective contributions add to zero with the contributions of the $p$ projective lines in $Z(y)$ which intersect it but are not contained in $D(y_2)$ (as explained in the case that $a_2$ is even and $a_2<a_1$).
This yields the claim.

Again, in both cases of the claim, this is also equal to the intersection  multiplicity $((D(y_1) \cap D(y_2)), (D(y_1)\cap Z(y)))$ given by  \cite{T}, Theorem 5.1.

Finally, the case that $a_1 < a_2$ and $a_2$ is even is treated in the same way  as the second case of the above claim (the intersection multiplicity is also again $-2(p-1)p^{\frac{a_1-1}{2}}$).

\qed
\newline

Now we are ready to start the intersection calculus. We fix the notation $a_r=\nu_p(j_r^2)$ and write $j_r^2=\varepsilon_r p^{a_r}$. 
 By Proposition \ref{diag} we may assume that $j_1,j_2,j_3$ are pairwise perpendicular to each other. From now on we make this assumption (unless otherwise mentioned). Further we assume that $2 \leq a_1 \leq a_2 \leq a_3$. 
To calculate $(D(j_i),D(j_l),D(j_k))=((D(j_l)\cap D(j_i)), (D(j_k)\cap D(j_i)))$ (intersection multiplicity in $D(j_i)$) we want to use Propositions \ref{unglmult} - \ref{horiz0} and Lemmas \ref{selsch} - \ref{mitte0}. 
Note that we always find $i \in \{1,2,3\}$ such that $j_i$ is odd. This follows from \cite{Ku1}, 1.16, see also (\ref{iden}) below.

\begin{Pro}
We have the following identities of intersection multiplicities.
\begin{enumerate}[(i)]
\item If $a_1$ is even and $a_2$ is odd, then $(D(j_1),D(j_2),D(j_3))=2p^{(a_1+a_2-3)/2}(p - 1)$.
\item If $a_1$ and $a_3$ are odd and $a_2$ is even, then $(D(j_1),D(j_2),D(j_3))=0$.
\item If $a_1$ and $a_2$ are odd and $a_3$ is even, then 

$(D(j_1),D(j_2),D(j_3))=
-2p^{(a_1+a_2-4)/2}(\frac{a_1+1}{2}p-\frac{a_1-1}{2})(p - 1)$.

\item If $a_1$ is odd and $a_2$ and $a_3$ are even and $a_2<a_3,$ then  $(D(j_1),D(j_2),D(j_3))=0$.
\item If $a_1$ is odd and $a_2$ and $a_3$ are even and $a_2=a_3,$ then 

 $(D(j_1),D(j_2),D(j_3))=-2p^{(a_1+a_2-3)/2}(\frac{a_1+1}{2}p - \frac{a_1-1}{2})$.

\end{enumerate}
\end{Pro}

The case that $a_1$ and $a_2$ are even and the case  that $a_1$ and $a_2$ and $a_3$ are odd will be treated later.

\emph{Proof.} Let us consider the first case. In view of Propositions \ref{horiz} and \ref{horiz0} we have to show that $(D(j_1)\cap D(j_2)_v,D(j_3)\cap D(j_2)_v)=0$ (intersection multiplicity in $D(j_2)$). But this follows from Lemma \ref{mitte0}. The same reasoning shows the second and the fourth case.

Let us consider the third case. For each $\mathbb{P}_{[\Lambda]}$ blonging to $D(j_3)$ and to $D(j_1)$  we compute the contribution to the intersection multiplicity which comes from the part of $D(j_3)\cap D(j_1)$ which has support in $\mathbb{P}_{[\Lambda]}$ and we add the contributions for the several projective lines.  (More precisely, if $D(j_3)\cap D(j_1)$ contains $\mathbb{P}_{[\Lambda]}$ with multiplicity $m$ (Propositions \ref{unglmult}, \ref{glmult}), we compute  the intersection  $(m\cdot \mathbb{P}_{[\Lambda]}, D(j_2)\cap D(j_1))$ in $D(j_1)$ and add the contributions of the several $\mathbb{P}_{[\Lambda]}$.)  
By Lemma \ref{mitte0} we only need to consider such $\mathbb{P}_{[\Lambda]}$ which belong to $B(j_2)$ or which do not belong to $D(j_2)$.
First we assume  that $a_1 < a_2$.

First we  consider the contributions of such $\mathbb{P}_{[\Lambda]}$, which belong to $B(j_2) \cap B(j_3)$. This contribution is $-2p +p=-p$ if $\mathbb{P}_{[\Lambda]}$ does not belong to $B(j_1)$ and it is $(-p)\cdot (\frac{a_1+1}{2})^2+ (\frac{a_1+1}{2}-1)\cdot \frac{a_1+1}{2}\cdot p$ if it belongs to $B(j_1)$.
Next it is easy to see that there are precisely $2(p-1)p^{\frac{a_2-a_1}{2}-1}\cdot p\cdot {(p^2)}^{\frac{a_1-1}{2}-1}$ projective lines belonging to $B(j_2)\cap B(j_3)$ but not to $B(j_1)$: We start  with one of the two projective lines in the apartment of projective lines belonging to $Cent(j_1)$ and to $Cent(j_2)$ such that $D(j_2)_p$ and $D(j_3)_p$ contain this projective line with the same multplicity (given by Corollary \ref{spezf}).  There are $p-1$ neighboring projective lines belonging to $Cent(j_1)$ but not to $Cent(j_2)$. For each such choice there are $p$ new neighboring projective lines belonging to $Cent(j_1)$.  We continue walking in $Cent(j_1)$ and enlarging our distance to the mentioned apartement. After $\frac{a_2-a_1}{2}-1$ such steps, the distances to $B(j_1)$ and to $B(j_2)$ are equal. Now there are $p\cdot {(p^2)}^{\frac{a_1-1}{2}-1}$ possible ways to go on until $B(j_2)$ but not to $B(j_1)$ (in each step enlarging the distance to the apartment by one). 
Similarly there are precisely 
$2(p-1)p^{\frac{a_2-a_1}{2}-1}\cdot (p^2-p)\cdot ({p^2)}^{\frac{a_1-1}{2}-1}$ projective lines belonging to $B(j_1) \cap B(j_2)\cap B(j_3)$.

Next is easy to see that the contribution of all $\mathbb{P}_{[\Lambda]}$ which do not belong to $B(j_3)$ sum up to zero. (More precisely, one checks the following fact for any  projective line $\mathbb{P}_{[\Lambda^{'}]}$ which lies in the 
 apartment of projective lines belonging to
  $Cent(j_1) \cap Cent(j_2)$ and for which
 the  multiplicities of  $\mathbb{P}_{[\Lambda^{'}]}$  in $D(j_2)_p$ and in  $D(j_3)_p$ are different:  
 The contributions coming from all projective lines $\mathbb{P}_{[\Lambda]}$ in       
 $B(j_2)$ (or intersecting it) with the property that $\mathbb{P}_{[\Lambda^{'}]}$ is the projective line in the apartment such that its distance to $\mathbb{P}_{[\Lambda]}$ is minimal, sum up to $0$.)

 Thus 
\begin{equation*}
\begin{split}
(D(j_1),D(j_2),D(j_3))= 2(p-1)p^{\frac{a_2-a_1}{2}-1}\cdot p\cdot {(p^2)}^{\frac{a_1-1}{2}-1}\cdot (-p) &\\ +2(p-1)p^{\frac{a_2-a_1}{2}-1}\cdot (p^2-p)\cdot {(p^2)}^{\frac{a_1-1}{2}-1}  &\\ \times ((-p)\cdot (\frac{a_1+1}{2})^2+ (\frac{a_1+1}{2}-1)\cdot \frac{a_1+1}{2}\cdot p)  &\\=-2p^{(a_1+a_2-4)/2}(\frac{a_1+1}{2}p-\frac{a_1-1}{2})(p - 1).
\end{split}
\end{equation*}
In case $a_1=a_2$ the calculation is almost the same, we have only to replace  the above formulas  for the number of projective lines which belong to $B(j_2)\cap B(j_3)$ but not to $B(j_1)$ resp. which belong to $B(j_1)\cap B(j_2) \cap B(j_3)$ by the  expressions $2(p-1)\cdot {(p^2)}^{\frac{a_1-1}{2}-1}$ resp. $2(p-1)^2\cdot {(p^2)}^{\frac{a_1-1}{2}-1}$. The resulting expression for $(D(j_1),D(j_2),D(j_3))$ does not change.

The fifth case is done analogously.
\qed
\newline

Next we come to the case that $a_1$ and $a_2$ are even and (hence) $a_3$ is odd. We intersect in $D(j_3)$. Now there is a horizontal component in $D(j_1)\cap D(j_3)$ passing through the unique point in  $Core(j_1)=Core(j_2)$ and also one of 
 $D(j_2)\cap D(j_3)$ passing through the point in  $Core(j_1)=Core(j_2)$. If we proceed as before, then we are faced with the problem that we have to compute the intersection multiplicity of these two horizontal components in $D(j_3)$. To avoid this problem, we proceed as follows. First we choose $\gamma \in \Zlp$ such that $\nu_p(j_2^2+\gamma^2j_1^2)=a_2+2$ (which is possible, as follows from the fact that $\chi(\varepsilon_1\varepsilon_2)=\chi(-1)$, which in turn follows from the form of the matrix $S^{'}$ given in section 3). Now let  $j_2^{'}=j_2+\gamma j_1$. By Proposition \ref{diag}, we get $(Z(j_1),Z(j_2),Z(j_3))=(Z(j_1),Z(j_2^{'}),Z(j_3))$.
 Using also Proposition \ref{multlin} we also see  $(D(j_1),Z(j_2),D(j_3))=(D(j_1),Z(j_2^{'}),D(j_3))$. It follows from the construction of $j_2^{'}$ that $Core(j_1) \neq Core(j_2^{'})$ and that $ Core(j_2^{'})$ consists of a supersingular  point lying on one of the projective lines passing through the point in $ Core(j_1).$

\begin{Pro}
Suppose that $a_1$ and $a_2$ are even. Then $$ (D(j_1),Z(j_2),D(j_3))=p^{\frac{a_1+a_2}{2}}+p^{\frac{a_1+a_2}{2}-1}-2p^{a_1-1}.$$
\end{Pro}
\emph{Proof.} We compute the value of $(D(j_1),Z(j_2^{'}),D(j_3))$. We write  $(D(j_1),Z(j_2^{'}),D(j_3))=((D(j_1)\cap D(j_3)), (Z(j_2^{'})\cap D(j_3)))$ as intersection multiplicity in $D(j_3).$ Assume first that $a_2>2 $ and $a_2<a_3-1$. We write $(D(j_1)\cap D(j_3))=v_1+h_1$, where $v_1$ denotes the vertical part (i.e. the part with support in the special fibre) of $D(j_1)\cap D(j_3)$ as divisor in $D(j_3)$ and $h_1$ denotes the horizontal part. Similarly we write $Z(j_2^{'}) \cap D(j_3)=v_2+h_2$. Then $(D(j_1),Z(j_2^{'}),D(j_3))=(v_1,h_2)+(h_1,v_2)+(v_1,v_2)$. (Note that by the above construction of $j_2^{'},$ the horizontal part $h_1$ and $h_2$ do not intersect.)
Using Propositions \ref{horiz} and \ref{unglmult} we compute $$(v_1,h_2)= (p^{\frac{a_1}{2}-1}+p^{\frac{a_1}{2}-2})(1+(p-1)+...+(p-1)p^{\frac{a_2}{2}})=(p^{\frac{a_1}{2}-1}+p^{\frac{a_1}{2}-2})p^{\frac{a_2}{2}+1}.$$ (Here the first summand 1 in the second parenthesis claims that the horizontal component of $Z(j_2^{'}/p^{a_2/2+1}) \cap D(j_3)$ intersects the two projective lines which contain $Cent(j_2^{'})$ each with multiplicity $1$. This follows by combining Lemma \ref{drinsenk} and point 3 of Lemma \ref{diperg}.)

Similarly, we compute $$(h_1,v_2)=(p-1)p^{\frac{a_1}{2}-1}(1+(1+p)+...+(p^{\frac{a_2}{2}-1}+p^{\frac{a_2}{2}})).$$ 
Next one checks that the the contribution to $(v_1,v_2)$ coming from the part of $(Z(j_2^{'})\cap D(j_3))_v$ which has support in the two projective lines passing through the point of $Core(j_1)$ is exactly $-( h_1,v_2)$.

To calculate the rest of $(v_1,v_2)  $  we write $ v_2 =(Z(j_2^{'})\cap D(j_3))_v= \sum_r  (D(j_2^{'}/p^r) \cap D(j_3))_v$  and see (using Lemma \ref{mitte0}) that for each $r$ we only have to consider the contributions coming from the projective lines in $(D(j_2^{'}/p^r) \cap D(j_3))_v$ which are in $B(j_1)$ or which intersect a projective line in $B(j_1)$. Further these contributions add to zero  if $r \neq \frac{a_2-a_1}{2}, \frac{a_2-a_1}{2}-1$. In each of these two cases these contributions add to $(-p)(p^2)^{\frac{a_1}{2}-1}$.   
 Adding everything we get the desired result.
 The cases that $a_1=2$ and/or $a_2=a_3-1$ are computed analogously.
  \qed
\newline

Next we treat the case that $a_1a_2a_3$ is odd. Here we are faced with the problem that in case $a_i=a_k$ and $\chi(-\varepsilon_i\varepsilon_k)=1$ we cannot apply Propositions \ref{glmult}  and \ref{horiz0} to get information about the structure of $D(j_i) \cap D(j_k)$.
But if in this case $a_l\neq a_i,a_k$ we can compute $((D(j_i)\cap D(j_l)),( D(j_k)\cap D(j_l)))$ as intersection multiplicity in $D(j_l)$. Thus the only case where a problem arises, is the case that $a_1=a_2=a_3=:a$ and at least two of the three expressions  $\chi_{rs}:= \chi(-\varepsilon_r \varepsilon_s)$ are equal to $1$. Suppose  $\chi_{ij}=\chi_{ik}=1$. But then the matrix $T= \diag(\varepsilon_1 p^{a}, \varepsilon_2 p^{a}, \varepsilon_3 p^{a})$ is $\GL_3(\Zlp)$-equivalent to a matrix of the form $\diag(\varepsilon_i p^a, \eta_k\varepsilon_k p^a, \eta_l\varepsilon_l p^a)$, where $\eta_l, \eta_k \in \Zlp^{\times}$ such that $\chi(\eta_l)= \chi(\eta_k)=-1$. Then $\chi(-\varepsilon_i \eta_k\varepsilon_k)=\chi(-\varepsilon_i \eta_l \varepsilon_l)=-1$. Since we are finally interested in the value of $(Z(j_1),Z(j_2),Z(j_3))$ and since for this  we may replace $T$ by  $\diag(\varepsilon_i p^a, \eta_k\varepsilon_k p^a, \eta_l\varepsilon_l p^a)$ (Proposition \ref{diag}), we are reduced to the case that $\chi_{ij} =\chi_{ik}=-1$.

\begin{Pro}
Suppose that $a_1a_2a_3$ is odd. Then we have the following identities of intersection multiplicities.
\begin{enumerate}[(i)]
\item If $a_1=a_2=a_3=:a$ and $\chi_{12}=\chi_{13}=-1$ and $\chi_{23}=1,$ then $$(D(j_1),D(j_2),D(j_3))=-\frac{a+1}{2}p^{a}+3\frac{a+1}{2}p^{a-1}-(a-1)p^{a-2}.$$
\item If $a_1=a_2=a_3=:a$ and $\chi_{12}=\chi_{13}=\chi_{23}=-1,$ then $$(D(j_1),D(j_2),D(j_3))=-\frac{a+1}{2}p^{a}+3\frac{a+1}{2}p^{a-1}-2(a-1)p^{a-2}.$$
\item If $a_1\leq a_2<a_3$ and $\chi_{12}=-1,$ then $$(D(j_1),D(j_2),D(j_3))=-2p^{\frac{a_1+a_2-4}{2}}(p-1)(\frac{a_1+1}{2}p-\frac{a_1-1}{2}).$$
\item If $a_1=a_2<a_3$ and $\chi_{12}=1,$ then $$(D(j_1),D(j_2),D(j_3))=2p^{a_1-1}.$$
\item If $a_1<a_2=a_3$ and $\chi_{12}=\chi_{13}=1,$ then $$(D(j_1),D(j_2),D(j_3))=-p^{\frac{a_1+a_2-4}{2}}(p+1)(\frac{a_1+1}{2}p-\frac{a_1-1}{2}).$$
\item If $a_1<a_2=a_3$ and $\chi_{12}=-1$  and $\chi_{13}=1,$ then $$(D(j_1),D(j_2),D(j_3))=-p^{\frac{a_1+a_2-4}{2}}(p-1)(\frac{a_1+1}{2}p-\frac{a_1-1}{2}).$$
\item If $a_1<a_2<a_3$ and $\chi_{12}=1,$   then $$(D(j_1),D(j_2),D(j_3))=0.$$

\end{enumerate}
\end{Pro}
 Note that the calculation of the intersection multiplicity $(Z(j_1),Z(j_2),Z(j_3))$ in case  $a_1<a_2=a_3$ and $\chi_{13}=-1$ can be reduced to the case case  $a_1<a_2=a_3$ and $\chi_{13}=1$ by the same reasoning as the one before the proposition.
 \newline
 
\emph {Proof.} 
In case (i) we write $(D(j_1),D(j_2),D(j_3))=((D(j_1)\cap D(j_2)),(D(j_1)\cap D(j_3))),$ and we compute for each projective line in $D(j_1)\cap D(j_3)$ the contribution to the intersection multiplicity $ ((D(j_1)\cap D(j_2)),(D(j_1)\cap D(j_3)))$ which comes from the part of  $D(j_1)\cap D(j_3)$ which has support in that projective line. By Lemma \ref{mitte0} we only have to consider such projective lines, which are in $B(j_2),$ or which are not in $D(j_2)$ but intersect an element in $B(j_2)$.

First we consider projective lines in $B(j_1) \cap B(j_2) \cap B(j_3)$. Then each contributes with the value $-p\cdot \frac{a+1}{2}$ (Proposition \ref{glmult}) and there are precisely $(p^2-3p+2)(p^2)^{\frac{a-1}{2}-1}$ such.
Next, there are $(p-3)(p^2)^{\frac{a-1}{2}-1}$ projective lines which are in $B(j_2) \cap B(j_3)$ but not in $B(j_1)$. Each contributes the value $-p$.
Further, there are $(p-1)(p^2)^{\frac{a-1}{2}-1}$ projective lines which are in $B(j_1) \cap B(j_2)$ and $D(j_3)$ but not in $B(j_3)$. Each contributes with value $-p$.

Next we count the contributions coming from the projective lines in $B(j_2)$ which are not in $B(j_1)\cup B(j_3)$.
For each such we also add the contributions coming from the projective lines in $D(j_1)\cap D(j_3)$ which intersect it but are not in $D(j_2)$. Now, if $\mathbb{P}_{[\Lambda]} \in B(j_2)$ but not in $B(j_1)\cup B(j_3)$ and if $\mult_1(\mathbb{P}_{[\Lambda]})\neq \mult_3(\mathbb{P}_{[\Lambda]}),$ then the contribution coming from $\mathbb{P}_{[\Lambda]}$ (including the contributions coming from the projective lines in $D(j_1)\cap D(j_3)$ which intersect 
$\mathbb{P}_{[\Lambda]}$ but are not in $D(j_2)$) is easily seen to be $0$. 
Next for each $m$ between $1$ and $\frac{a-1}{2}$ we count the contributions coming from such $\mathbb{P}_{[\Lambda]} \in B(j_2)$ but not in $B(j_1)\cup B(j_3)$ for which $\mult_1(\mathbb{P}_{[\Lambda]})=\mult_3(\mathbb{P}_{[\Lambda]})=p^m$. In case $m=\frac{a-1}{2}$ there are precisely two such (lying in the apartment of projective lines lying in $Cent(j_1)\cap Cent(j_3)$).
The contribution for each of the two is then $(-2p)p^{\frac{a-1}{2}}+p\cdot p^{\frac{a-1}{2}}+p^{\frac{a-1}{2}}+(2p-2)p^{\frac{a-1}{2}-1}+2(p^2-2p+1)p^{\frac{a-1}{2}-1}$. (The first two summands come from the projective line itself, the rest comes from the projective lines in $D(j_1)\cap D(j_3)$ which intersect it but do not belong to $D(j_2)$.)
Thus the contribution coming from the case $m=\frac{a-1}{2}$ is $2((-2p)p^{\frac{a-1}{2}}+p\cdot p^{\frac{a-1}{2}}+p^{\frac{a-1}{2}}+(2p-2)p^{\frac{a-1}{2}-1}+2(p^2-2p+1)p^{\frac{a-1}{2}-1})=2p^{\frac{a-1}{2}}(p-1)$. Similarly, for $m < \frac{a-1}{2}$, the contribution becomes $2(p^2-2p+1)(p^2)^{\frac{a-1}{2}-m-1}((-2p)p^m(\frac{a-1}{2}-m)+p\cdot p^{m}(\frac{a-1}{2}-m)+p^2p^{m-1}(\frac{a-1}{2}-m+1))=2(p^2-2p+1)p^{a-m-2}$.
Now summing everything up, we get the desired result.

Case (ii) is done analogously.

In case (iv)

we  write $(D(j_1),D(j_2),D(j_3))=((D(j_1)\cap D(j_3)),(D(j_2)\cap D(j_3)))$ (intersection multiplicity in $D(j_3)$). 
As usual, for each projective line in $(D(j_1)\cap D(j_3))$ we compute the contribution to intersection multiplicity 
 coming from the part of $(D(j_1)\cap D(j_3))$ which has support in that  projective line. We only have to consider projective lines in $B(j_2)$ or   intersecting a projective line in $B(j_2)$.
There are no projective lines in $B(j_2)\cap B(j_3)$ which  belong to $D(j_1)$.
We count the contribution coming from the projective lines in $B(j_1)\cap B(j_2)$. Since $\chi_{12}=1$ there precisely $(p^2-2p-1)(p^2)^{\frac{a_1-1}{2}-1}$ such and each contributes with the value $-p$ to the intersection multiplicity. Similarly one checks that the projective lines in $B(j_2)$ which are not in $B(j_1)$ together with the projective lines which intersect such altogether contribute with the value $(p^2-1)p^{a_1-2}$. Summing everything up gives the claimed result.

In case (v) we write $(D(j_1),D(j_2),D(j_3))=((D(j_1)\cap D(j_2)),(D(j_1)\cap D(j_3)))$ (intersection multiplicity in $D(j_1)$). Again we have to compute the contribution coming from the part of $(D(j_1)\cap D(j_3))$ which has support in a projective line in $B(j_2)$ or a projective line intersecting one in $B(j_2)$. Since $\chi_{13}=1$ it follows that, if $\mathbb{P}_{[\Lambda]} \in B(j_2)$ and   $\mathbb{P}_{[\Lambda]}\subset D(j_1)\cap D(j_3),$ then $\mathbb{P}_{[\Lambda]}\in B(j_3)$. Now it is easy to see that if in this situation $\mathbb{P}_{[\Lambda]}\notin B(j_1),$ then the contribution coming from $\mathbb{P}_{[\Lambda]}$ is $-p$ and that there are precisely 
$(p+1)\cdot p^{\frac{a_2-a_1}{2}-1}\cdot p \cdot (p^2)^{\frac{a_1-1}{2}-1}$ such projective lines. Similarly,
if $\mathbb{P}_{[\Lambda]}\in B(j_1),$ then the contribution coming from $\mathbb{P}_{[\Lambda]}$ is $\frac{a_1+1}{2}(-p)$ and there are precisely $(p+1)\cdot p^{\frac{a_2-a_1}{2}-1}\cdot (p^2-p) \cdot (p^2)^{\frac{a_1-1}{2}-1}$ such projective lines. Summing up all contributions, we obtain the claim.

In case (iii) we proceed similarly. Writing  $(D(j_1),D(j_2),D(j_3))=((D(j_1)\cap D(j_2)),(D(j_1)\cap D(j_3)))$ (intersection multiplicity in $D(j_1)$), 
we have to compute the contribution coming from the part of $(D(j_1)\cap D(j_3))$ which has support in a projective line in $B(j_2)$ or a     
projective line intersecting one in $B(j_2)$.  Now we consider the apartment of projective lines corresponding to the intersection $Cent(j_1) \cap Cent(j_2)$. For any  projective line in
 $B(j_2)$ (or intersecting it) there is a unique projective line in the apartment, such that the distance of the two projective lines (i.e. the distance of the corresponding vertices in the building) is minimal. Now we consider the two projective lines  $\mathbb{P}_{[\Lambda_0]}, \mathbb{P}_{[\Lambda_1]}$ in the apartment such that $\mult_3(\mathbb{P}_{[\Lambda_0]})=\mult_2(\mathbb{P}_{[\Lambda_0]})$ and $\mult_3(\mathbb{P}_{[\Lambda_1]})=\mult_2(\mathbb{P}_{[\Lambda_1]})$. Then by the same reasoning as in case (v) the contribution coming from all projective lines in $B(j_2)$ whose distance to the apartment is the distance to  $\mathbb{P}_{[\Lambda_0]}$ or $\mathbb{P}_{[\Lambda_1]}$ gives the claimed intersection multiplicity. One easily checks that for any  other projective line $\mathbb{P}_{[\Lambda]}$ in the apartment the contributions coming from all projective lines $\mathbb{P}_{[\Lambda^{'}]}$ in $B(j_2)$ (or intersecting it) with the property that $\mathbb{P}_{[\Lambda]}$ is the projective line in the apartment such that its distance to $\mathbb{P}_{[\Lambda^{'}]}$ is minimal, sum up to $0$.

Case (vi) is proved analogously to case (v). 

Case (vii) follows from Lemma \ref{mitte0}. 
\qed
\newline

We now drop the assumption that $T$ is diagonal.
We want  to  obtain an expression for $(Z(j_1),Z(j_2),Z(j_3))$. 
Suppose that $T$ is
$\GL_3(\Zlp)$-equivalent to $\diag(\varepsilon_1p^{a_1},\varepsilon_2p^{a_2},\varepsilon_3p^{a_3})$,  where  $\varepsilon_i \in \Zlp^{\times} $ for all $i$ and $a_1\leq a_2 \leq a_3.$ 
 (Now the restriction that $2 \leq a_1, a_2, a_3$ is dropped.) To state the result,  we  introduce the following 
 invariants of
 $T$, comp. \cite{W}.  
  Let
 \[
 \tilde{\xi}=
 \begin{cases}
\chi(-\varepsilon_1\varepsilon_2)\ & \text{if $a_1 \equiv a_2 \mod 2,$} \\
0 \ & \text{if $a_1 \not\equiv a_2 \mod 2$,} 
 \end{cases}
 \]
and let 
\[
\sigma=
 \begin{cases}
2\ & \text{if $a_1 \equiv a_2 \mod 2,$ } \\
1 \ & \text{if $a_1 \not\equiv a_2 \mod 2$.} 
 \end{cases}
 \]
Further, let 
\[
\eta=
 \begin{cases}
+1\ & \text{if $T$ is isotropic over $\Qp$,} \\
-1 \ & \text{if $T$ is anisotropic over $\Qp$.} 
 \end{cases}
 \]
To distinguish whether $T$ is isotropic or anisotropic over $\Qp$, we recall the following fact (see \cite{W}, p. 189). Let $i,j \in \{1,2,3 \}$
with $i \neq j$ and $a_i \equiv a_j \mod 2$, and define $k\in\{1,2,3 \}$ by $\{i,j,k\}=\{1,2,3\}$. Then
$T$ is isotropic over $\Qp$ if and only if $\chi(-\varepsilon_i\varepsilon_j)=1$ or $a_k \equiv a_j \mod 2$. On
the other hand, since $T$ is represented by $V^{'}$ we have
\begin{equation}\label{iden}
-1=(-1)^{a_1+a_2+a_3}\chi(-1)^{a_1+a_2+a_3+a_1 a_2 + a_2 a_3 + a_1 a_3}
\chi(\varepsilon_1)^{a_2+a_3}\chi(\varepsilon_2)^{a_1+a_3}\chi(\varepsilon_3)^{a_1+a_2},
\end{equation}
see \cite{KR2}, section 7 or \cite{Ku1}, (1.16).
Using this we can decide, whether $T$ is isotropic or anisotropic  over $\Qp$, and we can also determine $\tilde{\xi}$ from $a_1,a_2,a_3$, except in case that $a_1a_2a_3$ is odd (then $\tilde{\xi}=-1$ and $\tilde{\xi}=1$ is possible).

\begin{The}\label{expform} 
Suppose that $T$ is $\GL_3(\Zlp)$-equivalent to $\diag(\varepsilon_1p^{a_1},\varepsilon_2p^{a_2},\varepsilon_3p^{a_3})$, where  $\varepsilon_i \in \Zlp^{\times} $ for all $i$ and $a_1\leq a_2 \leq a_3.$ Then
there is the following explicit expression for the intersection multiplicity  $(Z(j_1),Z(j_2),Z(j_3))$.
\[
\begin{split}
(Z(j_1),Z(j_2),Z(j_3))&= -\sum_{i=0}^{a_1}\sum_{j=0}^{\frac{a_1+a_2-\sigma}{2}-i}p^{i+j}(-1)^i (i+2j) \\
& - \eta \sum_{i=0}^{a_1}\sum_{j=0}^{\frac{a_1+a_2-\sigma}{2}-i}p^{\frac{a_1+a_2-\sigma}{2}-j}(-1)^{a_3+\sigma+i}(a_3+\sigma+i+2j) \\
& -  \tilde{\xi}^2 p^{\frac{a_1+a_2-\sigma}{2}+1}\sum_{i=0}^{a_1}\sum_{j=0}^{a_3-a_2+2\sigma-4}\tilde{\xi}^j
(-1)^{a_2-\sigma+i+j} (a_2-\sigma+2+i+j).
\end{split}
\]
\end{The}
\emph{Proof.}
We already know by Proposition \ref{diag} that $(Z(j_1),Z(j_2),Z(j_3))$ only depends on the $\GL_3(\Zlp)$-equivalence class of $T$. Thus we may assume that $T=\diag(\varepsilon_1p^{a_1},\varepsilon_2p^{a_2},\varepsilon_3p^{a_3})$. Next by comparing with \cite{KR2}, Proposition 6.2 resp. \cite{T}, Theorem 5.1 we get the claim in case $a_1=0$ resp. $a_1=1$. Now by induction on $a_1+a_2+a_3$ we see (using Proposition \ref{multlin}) that it is enough to show that the above formula predicts the same values of $(D(j_1),D(j_2),D(j_3))$ resp. $(D(j_1),Z(j_2),D(j_3))$ as  the propositions of this section. 
This is in all cases checked by a straightforward calculation.
\qed
 

\section{The connection to representation densities and Eisenstein series}

In this section we want to express the local intersection multiplicity $(Z(j_1),Z(j_2),Z(j_3))$ in terms of certain representation densities and the global intersection multiplicity $\chi_T(Z_1,Z_2,Z_3)$ in terms of the derivative of the $T$-th Fourier coefficient  a certain Eisenstein series.
\newline

First we recall that, for $S\in \begin{rm}  Sym \end{rm}_m(\mathbb{Z}_p)$ and 
$U\in \begin{rm}  Sym \end{rm}_n(\mathbb{Z}_p)$
with $\begin{rm}  det \end{rm}(S)\neq 0$ and $\begin{rm}  det \end{rm}(U)\neq 0$,
 the representation density is defined as 
\[
\alpha_p(S,U)= 
\operatorname*{lim}_{t\rightarrow\infty} p^{-tn(2m-n-1)/2} \mid \{x \in M_{m,n}(\mathbb{Z}/p^t\mathbb{Z});
\ S[x]-U \in p^t\begin{rm}  Sym \end{rm}_n(\mathbb{Z}_p)\}\mid.
\]
Given $S$ as above, let 
\[
S_r=
\begin{pmatrix}
 S \\ 
 & 1_r \\
 & & -1_r
\end{pmatrix}.
\]
Then there is a rational function $A_{S,U}(X)\in \mathbb{Q}(X)$ of $X$ such that 
\[
\alpha_p(S_r,U)= A_{S,U}(p^{-r}).
\]
One defines 
\[
\alpha_p^{'}(S,U)=\frac{\partial}{\partial X}(A_{S,U}(X)) \arrowvert_{X=1}.
\]
(Comp. \cite{KR1}.) 
 Let $S=\begin{rm}  diag \end{rm}(1,-1,1,-\Delta)$. (Since $T$ is represented by the space $V^{'}(\Qp)$, it
 is not represented  by $V(\Qp)$ resp.  $S$, see \cite{Ku1}, Proposition 1.3. )

\begin{The}\label{theo2}
There is the following relation between intersection multiplicities and representation densities:
\[
(Z(j_1),Z(j_2),Z(j_3))= -\frac{p^4}{(p^2+1)(p^2-1)}\alpha_p^{'}(S,T).
\]
\end{The}
\emph{Proof.} 
Since both sides only depend on the $\GL_3(\Zlp)$-equivalence class of $T$, we may assume that $T$ is diagonal and equals $\diag(\varepsilon_1p^{a_1},\varepsilon_2p^{a_2},\varepsilon_3p^{a_3}),$  where  $\varepsilon_i \in \Zlp^{\times} $ for all $i$ and $a_1\leq a_2 \leq a_3.$ 
We compute the expression on the right hand side and compare with Theorem \ref{expform}. The right hand side can be expressed explicitly, combining a result of Katsurada and a result of Shimura. This was already done in \cite{T}  in case $a_1=1$, but the reasoning carries over. More precisely, by 
the same reasoning as in the proof of \cite{T}, Theorem 5.1, the function $A_{S,T}(X)$ is given as follows.
\[
A_{S,T}(X)=(1+p^{-2}X)(1-p^{-2}X^2)\tilde{F}_p(T;-X),
\]
where 
$\tilde{F}_p(T;X)$ is given by the following expression.
\[
\begin{split}
\tilde{F}_p(T;X) &= \sum_{i=0}^{a_1}\sum_{j=0}^{(a_1+a_2-\sigma)/2-i}p^{i+j}X^{i+2j} \\
& + \eta \sum_{i=0}^{a_1}\sum_{j=0}^{(a_1+a_2-\sigma)/2-i}p^{(a_1+a_2-\sigma)/2-j}X^{a_3+\sigma+i+2j} \\
& +  \tilde{\xi}^2p^{(a_1+a_2-\sigma+2)/2}\sum_{i=0}^{a_1}\sum_{j=0}^{a_3-a_2+2\sigma-4}\tilde{\xi}^j
X^{a_2-\sigma+2+i+j}.
\end{split}
\]
Here, the invariants $\eta, \sigma, \tilde{\xi} $ are defined as in the last section.
This yields an explicit expression for the right hand side of the formula in the statement of the theorem which is straighforward to calculate, and one checks that it is the same expression as the one given in  Theorem \ref{expform}.
\qed
\newline

Next we come to the comparison of the global intersection multiplicity $\chi_T(Z_1,Z_2,Z_3)$ and the  derivative of the $T$-th Fourier coefficient of  a certain Eisenstein series.

First we shortly recall some notations and facts introduced in \cite{KR2}. See \cite{Ku1}  and \cite{KR2}, section 7 for details.

 Let $K^{'}_p$ be the stabilizer in $G^{'}(\Qp)$ of a superspecial lattice $L$ 
in $X^{HB}$ (see section 2). Define $V^{'}(\Zp)=\End(L,F)\cap V^{'}_p$ and $V(\Zp)= \Lambda \otimes \Zp,$ where $\Lambda$ is the fixed self-dual lattice of section 1. 
 Let $K^{'}=K^{'}_p K^p \subset  G^{'}(\mathbb{A}_f).$ 
Let  $\varphi_f^p=\text{char}(\omega_1 \times \omega_2 \times \omega_3) \in S(V(\A_f^p)^3)$ and 
$\varphi_p =\text{char}(V(\Zp)^3) \in S(V(\Qp)^3)$ and  $\varphi_p ^{'}=\text{char}(V^{'}(\Zp)^3) \in S(V^{'}(\Qp)^3)$. 

 Let $\chi=\chi_V$ be the quadratic character of $\A^{\times}/\Q^{\times}$ associated to $V$, i.e., $\chi(x)=(x, \det(V))_{\A},$ where $(\ , \ )_{\A}$ denotes the global Hilbert symbol.
To $\varphi_p$ one associates $\Phi_p \in I_3(0,\chi_{V_p})$ (see \cite{Ku1} for an explanation of the  induced representation $I_3(0,\chi_{V_p})$ of $\Sp_6(\Qp)$, and more generally of  the  induced representation $I_3(s,\chi)$ of $\Sp_6(\A)$).  
Analogously  to  
 $\varphi_f^p$ one associates $\Phi_f^p$ and 
to $\varphi_p ^{'}$ one associates $\Phi_p ^{'}.$

 Then $\Phi_p$ is completed to an incoherent standard section $\Phi(s)=\Phi_{\infty}^2(s)\cdot \Phi_f^p(s)\cdot \Phi_p(s)\in I_3(s,\chi)$ (see \cite{Ku1}), where $\Phi^2_{\infty}(s)$ is associated to the Gaussian in $S(V^{'}(\R)^3)$. Denote by $pr(K^{'})$ the image of $K^{'}$ under the projection $pr: G^{'}(\A_f) \rightarrow SO(V^{'})(\A_f)$. Analogously one defines $pr(K).$ Let $\underline{x}=(j_1,j_2,j_3)\in V^{' 3}$ be the triple of our fixed three special endomorphisms (section 2) so that the matrix of the quadratic form with respect to $\underline{x}$ is our fixed $T.$ Finally let (having chosen Haar measures on $G^{'}(\A_f)$ and $Z^{'}(\A_f)$)
 $$
 O_T(\varphi_f^{'})=\int_{Z^{'}(\A_f)\setminus G^{'}(\A_f)} \varphi_f^{'}(g^{-1}\underline{x}g)dg.
  $$
 Analogously one defines  $O_T(\varphi_f^{p})$ and $O_T(\varphi_p^{'})$.
Let $P\subset \Sp_6$ be the Siegel parabolic subgroup (see introduction).
  We consider the following Eisenstein series on  $\Sp_6 (\A)$, 
 $$
E(g,s,\Phi)=\sum_{\gamma \in P(\Q) \setminus \Sp_6 (\Q) } \Phi(\gamma g,s).
$$
 Suppose that $\omega_1 \times \omega_2 \times \omega_3$ is locally centrally symmetric (i.e. invariant under the action of $\mu_2(\A_f^p)$). 
 Then for $h \in \Sp_{6}( \R)$  the derivative in $s=0$ of $E_T(h,s,\Phi)$ (see also introduction) can be expressed as 
 \[
 \begin{split}
 E^{'}_T(h,0,\Phi)&=\vol(SO(V^{'})(\R))\cdot pr(K^{'}))\cdot W_T^2(h) \\
 & \times \frac{W^{'}_{T,p}(e,0,\Phi_p)}{W_{T,p}(e,0,\Phi_p^{'})}\cdot \vol(K^{'})^{-1}\vol(Z(\Q)\setminus Z(\A_f)) \cdot O_T(\varphi_f^{'}).
 \end{split}
 \]
 See \cite{KR2} and \cite{Ku1} for an explanation of the Whittaker functions $W_T^2$ and $W_{T,p}$.
   With these preparations we are ready to state the final result.
 \begin{The}\label{Ziel}
 Let $h \in \Sp_6(\R)$ and suppose that $\omega_1 \times \omega_2 \times \omega_3$ is locally centrally symmetric. Then there is the following relation between Eisenstein series and intersection multiplicities.
 \[
 \begin{split}
  E^{'}_{T}(h,0,\Phi)&= -\frac{1}{2}\log(p) \cdot \kappa \cdot \chi_T(Z_1,Z_2,Z_3) \cdot W_T^2(h),
 \end{split}
 \]
 where $\kappa$ is the volume constant with value $\kappa=\vol(SO(V^{'})(\R)) \vol(pr(K)).$
 \end{The}
 
 \emph{Proof.} 
 By the above formula for $E^{'}_{T}(h,0,\Phi)$ we can write 
 $$
  E^{'}_{T}(h,0,\Phi)=c_1\cdot W^2_T(h)  \frac{O_T(\varphi_p^{'})}{W_{T,p}(e,0,\Phi_p^{'})} W^{'}_{T,p}(e,0,\Phi_p) O_T(\varphi_f^{p})
 $$
 for some constant $c_1$ (meaning that it is independent of $T$ and $h$ and $\omega_1 \times \omega_2 \times \omega_3$). Further we can write
 $
 W^{'}_{T,p}(e,0,\Phi_p)=c_2 \cdot  \alpha_p^{'}(S,T)
 $
 for some constant $c_2$ (see the proof  of \cite{KR2}, Proposition 7.2). Using Theorem \ref{theo2} we see that $
 W^{'}_{T,p}(e,0,\Phi_p)=c_3 \cdot  (Z(j_1),Z(j_2),Z(j_3))
 $ for some constant $c_3$. Combining this with Proposition \ref{locglob}, we obtain 
 $$  W^{'}_{T,p}(e,0,\Phi_p) O_T(\varphi_f^{p})= c_4 \cdot \chi_T(Z_1,Z_2,Z_3)$$  for some constant $c_4$. Next we claim 
 that the quotient 
$\frac{O_T(\varphi_p^{'})}{W_{T,p}(e,0,\Phi_p^{'})}$ is a constant. The proof for this is the same as the proof of formula (5.3.33) in the proof of Proposition 5.3.3 in \cite{KRY}.
Thus we obtain
$$
E^{'}_{T}(h,0,\Phi)=c \cdot \chi_T(Z_1,Z_2,Z_3)\cdot W_T^2(h)
$$
for some constant $c$. It remains to prove that $c=-\frac{1}{2}\vol(SO(V^{'})(\R)) \vol(pr(K)) \log(p)$. But by \cite{KR2}, Theorem 7.3, this is true provided that $T$ is not divisible by $p$. Since this claim is independent of $T$, it is always true.
This ends the proof.
\qed

\end{document}